\documentclass{article}
\usepackage{arxiv}

\usepackage[utf8]{inputenc}
\usepackage{url}
\usepackage{eurosym}
\usepackage{bm}
\usepackage{subfigure}
\usepackage{amsmath}
\usepackage{multirow}
\usepackage{booktabs}
\usepackage{caption}
\usepackage{graphicx}
\usepackage{float}
\usepackage[version=4]{mhchem}
\usepackage{siunitx}
\usepackage{bbold}
\usepackage{longtable,tabularx}
\usepackage{algorithm}
\usepackage{algorithmic}
\usepackage{multirow}
\usepackage[colorinlistoftodos]{todonotes}
\usepackage[nopostdot,acronym,nogroupskip,nonumberlist]{glossaries}

\newcommand{\bbb}[1]{\boldsymbol{#1}}

\newcommand{\bvec}[1]{\boldsymbol{#1}}

\makeglossaries
\newacronym{sa}{SA}{short arc}
\newacronym{cdm}{CDM}{Conjunction Data Message}
\newacronym{ar}{AR}{Admissible Region}
\newacronym{ads}{ADS}{Automatic Domain Splitting}
\newacronym{nlp}{NLP}{Nonlinear Programming Problem}
\newacronym{sqp}{SQP}{Sequential Quadratic Programming}
\newacronym{cam}{CAM}{Collision Avoidance Maneuver(s)}
\newacronym{da}{DA}{Differential Algebra}
\newacronym{soc}{SOC}{Second Order Cone}
\newacronym{socp}{SOCP}{Second-Order Cone Programming}
\newacronym{rso}{RSO}{Resident Space Object(s)}
\newacronym{dlr}{DLR}{German Aerospace Center}
\newacronym{esa}{ESA}{European Space Agency}
\newacronym{ca}{CA}{Closest Approach}
\newacronym{rtn}{RTN}{Radial, Transverse, and Normal}
\newacronym{raan}{RAAN}{Right Ascension of the Ascending Node}
\newacronym{leo}{LEO}{Low Earth Orbit}

\begin{document}

\title{Collision Avoidance Maneuver Optimization with a Multiple-Impulse Convex Formulation}

\author{Roberto Armellin \thanks{email: roberto.armellin@auckland.ac.nz} \\
  Te P\=unaha \=Atea – Auckland Space Institute,\\
  The University of Auckland \\
  20 Symonds Street, 1010 Auckland, New Zealand\\
  }

\maketitle

\begin{abstract}
 A method to compute optimal collision avoidance maneuvers for short-term encounters is presented. The maneuvers are modeled as multiple-impulses to handle impulsive cases and to approximate finite burn arcs associated either with short alert times or the use of low-thrust propulsion. The maneuver design is formulated as a sequence of convex optimization problems solved in polynomial time by state-of-the-art primal-dual interior-point algorithms. The proposed approach calculates optimal solutions without assumptions about the thrust arc structure and thrust direction. The execution time is fraction of a second for an optimization problem with hundreds of variables and constraints, making it suitable for autonomous calculations.   
\end{abstract}

\section{Introduction}
The number of catalogued \gls{rso}s is  growing due to spacecraft miniaturization (e.g., CubeSats) and the launch of mega-constellations (e.g., Starlink \cite{McDowell_2020}). In parallel, the number of tracked space debris is expected to increase due to the improvement in tracking systems (e.g., the space fence system \cite{haimerl2015space}), resulting in a larger number of conjunctions to process and more \gls{cam}s to design and execute. 

Conjunction analysis and collision avoidance are currently performed by agencies and operators on the ground using several tools and processes that were developed over the last twenty years \cite{Schiemenz2019}. These tools support operators' activities; however, the decision process still requires human intervention. This approach will not be practical in the future when conjunction screening, collision avoidance decision processes, and \gls{cam} design and execution need to be automated. Furthermore, more accurate conjunction services (i.e., limiting unnecessary maneuvers) and optimal \gls{cam}s will be required to reduce the propellant budget allocated for conjunction management. In this context, this work aims to propose a method for \gls{cam} optimization that is suitable for autonomous use with either high- or low-thrust capabilities.           

A \gls{cam} is performed when, at the time of closest approach, a threshold on the miss distance, or the collision probability, is exceeded \cite{Klinkrad2006a}. In \cite{Patera2003} a method to optimize an impulsive \gls{cam} based on the assumption of small maneuvers was introduced. This simplification allowed  the problem to be decoupled into maneuver direction and magnitude determination. The direction was determined by collision probability gradient while the magnitude with an iterative process. In \cite{Alfano2005} Alfano describes a tool for \gls{cam} analysis that can perform parametric studies of single-axis and dual-axis maneuvers. Collision probability contours for single-axis maneuvering are calculated based on an upper bound on the impulse magnitude and a range of permissible maneuver times. By selecting a specific time, contours are produced for dual-axis maneuvering. Agencies and operators use similar tools \cite{Schiemenz2019}. For example, a simple approach employing tangential maneuvers, thus sacrificing optimality, is adopted by the German Aerospace Center \cite{Aida2016}. The \gls{esa}'s tool CAMOS can instead deal with single and multiple impulses with arbitrary direction, different objective functions (e.g., collision probability, miss distance, total $\Delta v$), bounds on the maneuvers, and constraints on the post-maneuver orbital elements \cite{Cobo-Pulido2015}. The drawback of this flexibility is that the obtained solutions are only locally optimal, and therefore the analyst must critically analyze the results. In \cite{Morselli2014} a multi-objective approach for \gls{cam} design was presented that enabled an exhaustive analysis of the problem building Pareto optimal solutions according to multiple criteria. This single-impulse approach also allows for merging station-keeping with \gls{cam} and assessing collision risk for a one-week window after the maneuver. However, this approach is numerically intensive, as it requires multiple evaluations of complex objective functions. Bombardelli and Hernando-Ayuso \cite{Bombardelli2015} developed an analytical and semi-analytical method to find the impulse that minimizes either the miss distance or the collision probability for a given $\Delta v$ magnitude. The proposed methods have proven convergence as the problem is reduced either to an eigenvalue problem or a convex optimization one. 

The research on low-thrust debris avoidance is not very developed. A similar problem has been studied in a great deal of detail for formation flying collision-free reconfiguration. Both direct and indirect optimal control approaches were also proposed and compatible with onboard use (see the introduction in \cite{DiMauro2019} and the references therein for an exhaustive overview). Restricting the analysis to low-thrust \gls{cam} design, in \cite{Gonzalo2019} an approach based on averaged dynamics and Gauss variational equation is proposed with the underlying assumption of continuous tangential thrust. Four methods based on an indirect formulation of an optimal control problem are presented in \cite{Salemme2020}. It is concluded that a semi-analytical method based on the linearization of the dynamics offers the best compromise between accuracy and computational time. However, formulating an energy optimal control problem without bounds on thrust magnitude is a significant limitation for this approach. Recently, Hernando-Hayuso and Bombardelli \cite{Hernando-Ayuso2021} proposed a solution to the low-thrust \gls{cam} design problem by applying a maximum thrust for a fixed time span, while optimizing its orientation to yield minimum collision probability. The proposed method is valid for circular orbits and does not directly account for fuel optimality in the formulation.

In this work, we present a methodology for optimal \gls{cam} design with convergence properties and execution time suitable for autonomous use, potentially onboard spacecraft. The approach is suitable for both impulsive and low-thrust \gls{cam} design and can handle short-term encounters with warning times from few minutes to multiple orbital revolutions. Furthermore, no a priori assumption on the direction of the impulses is made, and an arbitrary dynamical model can be used. Constraints either on miss distance, maximum collision probability, or an approximated value of the collision probability can be enforced while minimizing the total $\Delta v$. The approach is based on framing a multiple-impulse \gls{cam} optimization problem as a convex optimization one \cite{Boyd2004}. Thanks to the proven existence and uniqueness of the solution and computational advantages ensured by polynomial complexity, convex optimization has found many applications in aerospace engineering over the last 15 years \cite{Liu2017}, including long-duration low-thrust transfers \cite{Tang2018}, \cite{Wang2018}. The multiple-impulse \gls{cam} optimization problem is a \gls{nlp}. Formally, three steps are required to transform the \gls{nlp} into a \gls{socp} problem \cite{Boyd2004}, due to the following three issues:
\begin{enumerate}
    \item the objective function (Eq.\eqref {eq:NLPobj}) is a nonlinear function of the optimization variables; 
    \item the maps describing the dynamics (Eq. \eqref{eq:SST}) and for the calculation of closest approach quantities (Eq. \eqref{eq:conj1}--\eqref{eq:conj4}) are nonlinear;
    \item the collision avoidance constraints are (Eq. \eqref{eq:NLPcon}) non-convex.
\end{enumerate}

The introduction of slack variables and lossless convexification \cite{Ackmese2011},\cite{Wang2018} solve the first issue. Only the linear terms in the dynamics and b-plane (refer to the next section for its definition) maps are retained, using \gls{da} implemented in DACE\footnote{The open-source software package is available at https://github.com/dacelib/dace} as a first-order automatic differentiation technique \cite{Rasotto2016}. Due to the small deviations introduced by CAMs, linearized maps accurately describe the conjunction dynamics \cite{Bombardelli2015,Salemme2020}. As a result, the problems associated with linearization (i.e., artificial infeasibility and unboundedness \cite{mao2018successive}) are not relevant here, and successive convexification is introduced for refinement purposes only, with no need for virtual controls and complex trust-region strategies \cite{Mao2018}. The third issue is tackled by working on the squared Mahalanobis distance, whose contour lines describe an ellipse on the b-plane. The squared Mahalanobis distance allows us to set a constraint either on minimum distance, maximum collision probability, or an approximation of the collision probability \cite{Alfriend1999}. However, the constraint on this quantity results in a non-convex problem, as the admissible region on the b-plane is non-convex (as is the keep-out zone in rendezvous dynamics \cite{Liu2014}). The projection and linearization technique proposed in \cite{Mao2017} is adopted to tackle this issue, resulting in a second iterative procedure. Depending on the chosen initial guess, the iterations can converge to different optima. However, at least for short alert times, two suitably selected guesses allow for the identification of the global minimum. 

The final \gls{cam} optimization approach is solved with the state-of-the-art primal-dual interior-point algorithm implemented in the software MOSEK \cite{mosek}. Solutions with hundreds of impulses are obtained robustly and efficiently, representing both high-thrust and low-thrust maneuvers. We test our algorithm on 2,170 real conjunctions derived from the ESA Collision Avoidance Challenge \url{https://kelvins.esa.int/collision-avoidance-challenge}    
The paper is organized as follows. In Sec. \ref{conjDin} a brief overview of the short-term conjunction dynamics is provided together with a nonlinear approach to study the effect of a \gls{cam} on the conjunction geometry. Section \ref{Sec2} contains a description of the methodology developed in this work. We first state the \gls{cam} design as \gls{nlp} problem, followed by the description of the steps required for its convexification. The algorithm's application to a large set of cases is described in \ref{NumRef}, leading to the conclusive section. 

\section{Conjunction Dynamics} \label{conjDin}
A brief introduction of the key quantities of a short-term encounter is provided. Only relevant concepts for the design of \gls{cam}s are summarized. Afterward, we introduce the \gls{da}-based approach to study the effect of maneuvers on the conjunction geometry used to define the optimization problem in Sec. \ref{Sec2}. 

\subsection{Collision Probability Computation}
\begin{figure}[h!]
\centering
\includegraphics[width=0.5\textwidth]{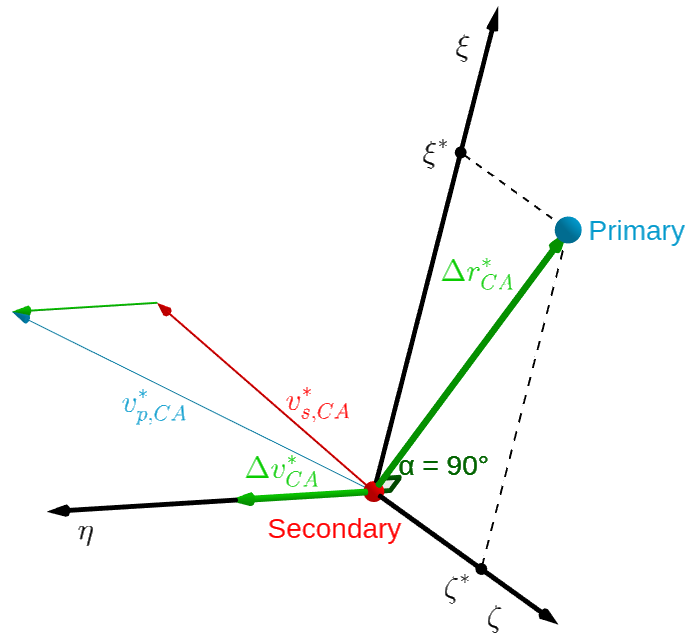}
\caption{Conjunction geometry and b-plane definition}
\label{Bplane}
\end{figure}

We consider the conjunction between a primary (subscript $p$) and a secondary (subscript $s$). The primary is the spacecraft we control, whereas the secondary is assumed to be passive. We indicate the relative position and velocity vectors at the \gls{ca} as
 \begin{align}
 \boldsymbol{\Delta r}^*_{CA} &= \bvec{r}^*_{p,CA} - \bvec{r}^*_{s,CA}\, ,\\
 \boldsymbol{\Delta v}^*_{CA} &= \bvec{v}^*_{p,CA} - \bvec{v}^*_{s,CA}\,
\end{align}
in which $\bvec{r}$ and $\bvec{v}$ are the absolute position and velocity vectors. In the remainder of the paper the asterisk is used to highlight quantities when the objects are not maneuvered.
At the \gls{ca} the distance between the two objects is minimum, and as such it follows 
\begin{equation}\label{eq:conj}
\boldsymbol{\Delta r}^*_{CA} \cdot \boldsymbol{\Delta v}^*_{CA} = 0.
\end{equation}
To compute the collision probability it is useful to introduce a coordinate system referred to as the b-plane. The origin of the axes of this frame lies at the centre of the secondary object at the time of \gls{ca}; the $\eta$-axis is defined along the direction of the relative velocity of the primary with respect to the secondary object; the $\xi\zeta$ plane is perpendicular to that $\eta$-axis 
 \begin{align}
 \hat{\boldsymbol{u}}_{\xi} &= \frac{\boldsymbol{v}^*_{s,CA}\times\boldsymbol{v}^*_{p,CA}}{\|\boldsymbol{v}^*_{s,CA}\times\boldsymbol{v}^*_{p,CA}\|}\, ,\\
 \hat{\boldsymbol{u}}_{\eta} &= \frac{\boldsymbol{v}^*_{p,CA}-\boldsymbol{v}^*_{s,CA}}{\| \boldsymbol{v}^*_{p,CA}-\boldsymbol{v}^*_{s,CA}\|}\, ,\\
 \hat{\boldsymbol{u}}_{\zeta} &=\hat{\boldsymbol{u}}_{\xi}\times\hat{\boldsymbol{u}}_{\eta}.
\end{align}
All the introduced quantities are shown in Fig. \ref{Bplane}. 

The unit vectors define the rotation matrix from the inertial reference frame to the b-plane
\begin{equation}
    R_{3D} = [\hat{\boldsymbol{u}}_{\xi}\   \hat{\boldsymbol{u}}_{\eta}\  \hat{\boldsymbol{u}}_{\zeta}]^T,
\end{equation}
while the projection in the $\eta$-axis is achieved by 
\begin{equation}\label{eq:project}
    R_{2D} = [\hat{\boldsymbol{u}}_{\xi}\     \hat{\boldsymbol{u}}_{\zeta}]^T.
\end{equation}

The nominal position of the primary on the b-plane at the time of closest approach $t^*_{CA}$ is  $\boldsymbol{\Delta r}^*_{CA} = (\xi^*, \zeta^*)$. The particular case in which $\Delta r^*_{CA} = 0$  is referred to as a direct impact. 

Under the short encounter approximation and a Gaussian distribution of objects state vectors, the collision probability is 
 \begin{equation}\label{eq:Pc_2D}
 P^*_C = \dfrac{1}{2\pi \left(\det{C^*_{CA}}\right)^{1/2}} \iint_{A} e^{-\frac{1}{2} \left(\bvec{\Delta r}-\bvec{\Delta r}^*_{CA}\right)^T ({C^*_{CA}})^{-1}\left(\bvec{\Delta r}-\bvec{\Delta r}^*_{CA}\right)} \textrm{d}\xi \textrm{d}\zeta  
\end{equation}
in which $C_{CA}^*$ is the sum of the positional covariance of the two objects referred to a common reference frame and projected onto the b-plane via \eqref{eq:project}. $A$ collision cross sectional area, a circle of radius $R = R_p + R_s$ where $R_{p/s}$ is the radius of the sphere enclosing the primary/secondary respectively. 
Several methods have been developed over the years to calculate $P_C$. Among them, \cite{Chan2008}, \cite{Serra2016}, and \cite{Garcia-Pelayo2016} are worth particular merit as they provide analytical solutions. As explained in \cite{Alfriend1999}, the simplest approach consists of approximating the value of $P_C$ by assuming the probability density is constant over the collision circle. The approximate value is    
 \begin{equation}\label{eq:approxPC}
 P^*_C = \dfrac{R^2}{2 \left(\det{{C^*_{CA}}}\right)^{1/2}} e^{-\frac{1}{2} {\left(d^*_{CA}\right)^{2}}},  
\end{equation}
in which $d^*_{CA} = \sqrt{{\bvec{\Delta r}^*_{CA}}^T ({C^*_{CA}})^{-1}\bvec{\Delta {r}}^*_{CA}}$ is the Mahalanobis distance. Additionally, the maximum collision probability can be obtained by optimally scaling the combined covariance (\cite{Alfriend1999}) resulting in  
\begin{equation}\label{eq:maxPc}
    P_{C,\max}^* = \dfrac{R^2}{\left({d^*_{CA}}^{2}\right) \left(\det{{C^*_{CA}}}\right)^{1/2} e}.
\end{equation}
These approximations allow us to use the squared Mahalanobis distance to set constraints on an approximate value of the collision probability through Eq. \eqref{eq:approxPC}, the maximum collision probability through Eq. \eqref{eq:maxPc}, or the miss distance (by setting $C^*_{CA} = \mathcal{I}$, the identity matrix). Moreover, for a constant covariance matrix, the contour lines of the squared Mahalanobis distance describe ellipses on the b-plane, a type of constraint that can be dealt with efficiently by successive convexifications, as is shown later.

\subsection{Effect of Maneuvers on the Conjunction}
When the primary is maneuvered, the conjunction geometry changes. This aspect is analyzed here with the use of arbitrary order Taylor expansions enabled by \gls{da} (the notation adopted in \cite{Armellin2017} is used). These effects are generally small for short-term encounters and small maneuvers. Nevertheless, a general treatment is provided as the proposed approach could potentially be applied to longer encounters with minimal changes (provided that the distance function is convex in time).  

The first step is to use \gls{da} to introduce perturbations to the closest encounter time $t_{CA}^* + \delta t$, the primary position $\bvec{r}^*_{p,CA} + \bvec{\delta} \bvec{r}_p$, and the primary velocity $\delta \bvec{v}^*_{p,CA} +  \bvec{\delta} \bvec{v}_p$. Taylor expansions of both the primary and secondary states are obtained by \gls{da}-based numerical integrations (expansion in time and state, see \cite{Armellin2010} for details) of the orbital dynamics, delivering 
\begin{equation}\label{eq:perturbedState}
\begin{array}{l}
 \bvec{r}_{p} = \mathcal{T}_{\bvec{r}_{p}}(\delta t, \bvec{\delta} \bvec{r}_p, \bvec{\delta}\bvec{v}_p)\\
  \bvec{v}_{p} = \mathcal{T}_{\bvec{v}_{p}}(\delta t,\bvec{\delta}\bvec{r}_p, \bvec{\delta} \bvec{v}_p)\\
   \bvec{r}_{s} = \mathcal{T}_{ \bvec{r}_{s}}(\delta t)\\ 
   \bvec{v}_{s} = \mathcal{T}_{\bvec{v}_{s}}(\delta t).\\
\end{array}
\end{equation}
Note that in Eq. \eqref{eq:perturbedState}, due to the effect of the introduced perturbations, we have dropped the \gls{ca} subscript. Additionally, the secondary's state is only affected by the time perturbation.

From Eq. \eqref{eq:perturbedState} the relative quantities can be calculated: 
\begin{equation}\label{eq:deltaPerturbed}
\begin{array}{l}
 \bvec{\Delta r} = \mathcal{T}_{\bvec{\Delta} \bvec{r}}(\delta t, \bvec{\delta} \bvec{r}_p, \bvec{\delta} \bvec{v}_p)\\
  \bvec{\Delta v} = \mathcal{T}_{\bvec{\Delta} \bvec{v}}(\delta t, \bvec{\delta} \bvec{r}_p, \bvec{\delta} \bvec{v}_p).\\
\end{array}
\end{equation}
The closest encounter condition, Eq. \eqref{eq:conj}, in DA formalism, reads   
\begin{equation}
\begin{array}{l}\label{eq:PIE}
 \bvec{\Delta r}\cdot \bvec{\Delta v} = \mathcal{T}_{\bvec{\Delta r}\cdot \bvec{\Delta v}}(\delta t, \bvec{\delta} \bvec{r}_p, \bvec{\delta} \bvec{v}_p) = 0.\\
\end{array}
\end{equation}
This constraint is a parametric implicit equation that can be solved for $\delta t$ using polynomial partial inversion techniques (see \cite{Armellin2010} and \cite{He2018} for more details). The polynomial partial inversion provides   
\begin{equation}
\begin{array}{l}
 \delta t  = \mathcal{T}_{\delta t} (\bvec{\Delta r}\cdot \bvec{\Delta v}, \bvec{\delta} \bvec{r}_p, \bvec{\delta} \bvec{v}_p)\\
\end{array}
\end{equation}
and, by substitution of the closest encounter constraint $\bvec{\Delta r}\cdot \bvec{\Delta v} = 0$, we obtain 
\begin{equation}\label{eq:tcaManeuvered}
\begin{array}{l}
 \delta t_{CA}  = \mathcal{T}_{\delta t_CA} ( \bvec{\delta} \bvec{r}_p, \bvec{\delta} \bvec{v}_p)\\
\end{array}
\end{equation}
This map gives a Taylor approximation for the variation in the closest encounter time due to change in the primary state (due to a maneuver), an aspect that is commonly ignored in CAM design. This polynomial can be inserted back in Eq. \eqref{eq:perturbedState}, obtaining 
\begin{equation}\label{eq:perturbedStateCA}
\begin{array}{l}
 \bvec{r}_{p,CA} = \mathcal{T}_{\bvec{r}_{p,CA}}(\bvec{\delta} \bvec{r}_p, \bvec{\delta} \bvec{v}_p)\\
  \bvec{v}_{p,CA} = \mathcal{T}_{\bvec{v}_{p,CA}}(\bvec{\delta} \bvec{r}_p, \bvec{\delta} \bvec{v}_p)\\
   \bvec{r}_{s,CA} = \mathcal{T}_{ \bvec{r}_{s,CA}}(\bvec{\delta} \bvec{r}_p, \bvec{\delta} \bvec{v}_p)\\ 
   \bvec{v}_{s,CA} = \mathcal{T}_{\bvec{v}_{s,CA}}(\bvec{\delta} \bvec{r}_p, \bvec{\delta} \bvec{v}_p).\\
\end{array}
\end{equation}
These polynomial maps approximate the states of both objects at the different times of \gls{ca} as a result of a \gls{cam} that produces a variation in position and velocity of the primary object at $t_{CA}^*$. Note that in Eq. \eqref{eq:perturbedStateCA} we re-introduce the \gls{ca} subscript, but we remove the asterisk, as now each perturbed solution has a different time of closest approach, determined by Eq. \eqref{eq:tcaManeuvered}. Similarly, the projection matrix in Eq. \eqref{eq:project} is expanded by using Eq. \eqref{eq:perturbedStateCA} for the calculation of the unit vectors, resulting in
\begin{equation}\label{eq:R2DExp}
    R_{2D} = \mathcal{T}_{R_{2D}}(\bvec{\delta} \bvec{r}_p, \bvec{\delta} \bvec{v}_p)
\end{equation}

Equation \eqref{eq:R2DExp} is used both to project the relative state and the combined covariance matrix on the b-plane, thus allowing for the calculation of the expansion of all the relevant conjunction quantities: 

\begin{equation}\label{eq:conj1}
\begin{array}{l}
 \bvec{\Delta r}_{CA} = \mathcal{T}_{\bvec{\Delta} \bvec{r}_{CA}}(\bvec{\delta} \bvec{r}_p, \bvec{\delta} \bvec{v}_p),
\end{array}
\end{equation}
\begin{equation}\label{eq:conj2}
\begin{array}{l}
 {d^2_{CA}} = \mathcal{T}_{{d^2_{CA}}}(\bvec{\delta} \bvec{r}_p, \bvec{\delta} \bvec{v}_p),\\
\end{array}
\end{equation}
\begin{equation}\label{eq:conj3}
\begin{array}{l}
 P_C = \mathcal{T}_{P_C}(\bvec{\delta} \bvec{r}_p, \bvec{\delta} \bvec{v}_p),\\
\end{array}
\end{equation}
and
\begin{equation}\label{eq:conj4}
\begin{array}{l}
 P_{C,\max} = \mathcal{T}_{P_{C,\max}}(\bvec{\delta} \bvec{r}_p, \bvec{\delta} \bvec{v}_p).
\end{array}
\end{equation}

Note that in Eqs. \eqref{eq:conj2}-\eqref{eq:conj4}  the covariance  $C_{CA}$ is function of $\bvec{\delta} \bvec{r}_p$ and $\bvec{\delta} \bvec{v}_p$ as a result of the projection matrix's dependency on the perturbed conjunction state (see Eq. \eqref{eq:R2DExp}). On the other hand, it is always assumed that the state covariances provided in the \gls{cdm} are not directly altered by the implementation of the maneuver.

\section{Collision Avoidance Maneuver Design}\label{Sec2}
The details of the \gls{cam} design algorithm for multiple-impulse maneuvers is presented. Before framing it as a successive convexification problem, a general \gls{nlp} formulation is described to provide the general setting.

\subsection{Nonlinear Programming Formulation}
A uniform $N$-point time grid is  constructed by selecting a discretization time step $\Delta t$ and starting from the earliest maneuvering time  $t_0$. Note that more refined discretization schemes are possible without adding complexity to the algorithm (e.g., using a suitable angular variable to better deal with eccentric orbits). 
Here $N = \min\left(\textrm{floor}\left(\dfrac{t^*_{CA}-t_0}{\Delta t}\right),N_{\max}\right)$, where $N_{\max}$ accounts for the maximum time span in which a maneuver can be implemented. At every discretization point a maneuver can be added in the form of an instantaneous change in velocity $\boldsymbol{\Delta v}_i$ with $i = 0,\dots, N-1$.
The nominal trajectory at discretization points is given by $\bvec{r}^*_i, \bvec{v}^*_i$ with $i = 0, \cdots, N$, with $\bvec{r}^*_N = \bvec{r}^*_{CA}, \bvec{v}^*_N = \bvec{v}^*_{CA}$. 
For each $i=0,\dots,N-1$ we calculate a $o$-th order Taylor approximation of the mapping between deviations in the initial state and deviations of the final state
\begin{equation}\label{eq:SST}
\begin{array}{l}
     \bvec{\delta r}_{i+1}^- = \mathcal{T}_{\bvec{\delta r}_{i+1}^-}(\bvec{\delta r}_{i}^+, \bvec{\delta v}_{i}^+ )  \\
     \bvec{\delta v}_{i+1}^- = \mathcal{T}_{\bvec{\delta r}_{i+1}^-}(\bvec{\delta r}_{i}^+, \bvec{\delta v}_{i}^+ ).  \\
\end{array}
\end{equation}
In Eq. \eqref{eq:SST} a $+$ superscript indicates chosen quantities at the beginning of an interval, whereas a $-$ propagated quantities from the previous interval. These maps are built with $N$ \gls{da} integrations of a dynamical model of choice using the unperturbed trajectory as a reference (e.g., obtained by backward propagation from the \gls{cdm} with the same dynamical model). As the \gls{cam} $\Delta v$s are small, a low order (in most cases order 2)  allows for sufficiently accurate approximations without the need of iterations.

The optimization variables are the set of $N$ $\bvec{\Delta v}_i$ applied at nodes $i = 0,\dots,N-1$. The effect of these maneuvers on the trajectory of the primary are obtained by the use of maps from Eq. \eqref{eq:SST}. In particular, for $i=0$ we can set $\bvec{\delta r}_{0}^+ = 0$ and $\bvec{\delta v}_{0}^+ = \bvec{\Delta v}_{0}$, and, by applying \eqref{eq:SST}, we obtain the mapped perturbations $\bvec{\delta r}^-_{1}$ and $\bvec{\delta v}^-_{1}$. For $i = 1,\dots,N-1$ we proceed by defining the new perturbations $\bvec{\delta r}_{i}^+ = \bvec{\delta r}^-_{i-1}$ and $\bvec{\delta v}_{i}^+ = \bvec{\delta v}^-_{i-1} + \bvec{\Delta v}_{i} $ and use the $i$-th set of maps of Eq. \eqref{eq:SST} to map these perturbations to the end of the segment. 
At $i=N-1$ we learn how the set $\mathbb{x} = [\bvec{\Delta v}_0; \dots; \bvec{\Delta v}_{N-1}]$ is mapped into the final perturbations $\bvec{\delta r}_p$ and $\bvec{\delta v}_p$, which then in turn allows us to compute the relevant quantities through Eqs. \eqref{eq:conj1}-\eqref{eq:conj4}. 
The multiple-impulse \gls{cam} optimization problem, can be stated as follows: 
minimize the sum of magnitudes of the impulses 
\begin{equation}\label{eq:NLPobj}
    \min_{\mathbb{x}} \sum_{i=0}^{N-1} \Delta v_i
\end{equation}
subject to nonlinear inequality constraints
\begin{equation}\label{eq:NLPcon}
\begin{array}{l}
    \Delta v_i \le \Delta \bar{v}\quad \textrm{for}\quad i = 0, \dots, N-1, \\
    P_{C} \le \bar{P}_{C} \quad \textrm{or}\quad
    P_{C,\max} \le \bar{P}_{C,\max} \quad \textrm{or}\quad
    \Delta r_{CA} \ge  \bar{d}_{\min},
\end{array}
\end{equation}
where the overline indicates assigned values, and $\Delta v_i$ the impulse magnitude. This optimization problem is a \gls{nlp} that can be solved with dedicated solvers. In this work we use the \gls{sqp} algorithm implemented within the MATLAB {\it fmincon} function providing analytical gradient of the objective function and Jacobian of the constraints. For the latter the derivatives included in the polynomials maps are used together with the chain rule to calculate the sensitivity of the constraints with respect to the optimisation vector. The \gls{nlp} formulation is summarized in Algorithm \ref{algo:NLP}. As \gls{nlp} problems, this formulation is non-deterministic polynomial-time hard (NP-hard), meaning that the computation time may be very long if the problem is solved at all \cite{Liu2014}.

\begin{algorithm}[!ht]
  \caption{Nonlinear programming formulation}\label{algo:NLP}
  \begin{algorithmic}[1]
     \STATE Get inputs from \gls{cdm}: $R$, $\boldsymbol{r}^*_{p/s,CA}$ $\boldsymbol{v}^*_{p/s,CA}$, $t^*_{CA}$, $C^*_{p/s, CA}$;
     \STATE Assign $t_0$, $\Delta t$, $\Delta \bar{v}$, $N$,$\bar{P}_{C}$ or $\bar{P}_{C,\max}$ or $\bar{d}_{\min}$, and $\mathbb{x}_0$;
     \STATE Back propagate the trajectories from $t^*_{CA}$ to $t_0$ and save $\boldsymbol{r}_{p/s,t_0}$, and $\boldsymbol{v}_{p/s,t_0}$; 
     \STATE Define the time grid $(t_0:\Delta t:t_0 + N\Delta t)$;  
     \STATE  Build maps Eq. \eqref{eq:SST} by $N$ $o$-th order DA forward propagations;
     \STATE Solve the \gls{nlp} problem defined by \eqref{eq:NLPobj} and \eqref{eq:NLPcon};
   \end{algorithmic}
\end{algorithm}

\subsection{Convex Problem Formulation}\label{3A}
As described in the introduction, three main steps are required to formulate the \gls{cam} design as a \gls{socp} problem. Firstly, in Sec. \ref{losslessConv} the objective function and the constraints on the velocity magnitude are reformulated by introducing slack variables and lossless convexification. Afterward, the dynamics are linearized in Sec. \ref{dynLinear}, followed by constraints linearization in Sec. \ref{consLin}.  Only the constraint on the squared Mahalanobis distance is taken into account as this type of constraint is handled by a projection and linearization approach \cite{Mao2017}. The details of the  algorithm are then presented in Sec. \ref{scvx}.   

\subsubsection{Lossless Convexification}\label{losslessConv} The $\Delta v$ magnitudes are introduced as slack variables in the optimization problem. As a result, each impulse is described by four independent variables $\bvec{\Delta \tilde{v}}_i = [\bvec{\Delta v}_i; {\Delta v}_i]$, and the  optimization vector becomes 
\begin{equation} \mathbb{x} = [\bvec{\Delta {v}}_0; \dots; \bvec{\Delta {v}}_{N-1}; \Delta {v}_0; \dots; {\Delta {v}}_{N-1}].
\end{equation}
Introducing the slack variables renders objective function linear 
\begin{equation}\label{eq:obj}
    \min_{\mathbb{x}} \sum_{i=0}^{N-1} \Delta v_i
\end{equation}
and transforms the $N$ constraints on the impulse magnitudes into \gls{soc} constraints
\begin{equation} \label{eq:SOCC}
\begin{array}{l}
    \sqrt{\Delta v_{i,x}^2 + \Delta v_{i,y}^2 + \Delta v_{i,z}^2}  \le \Delta {v}_i\quad \textrm{for}\quad i = 0, \dots, N-1. 
\end{array}
\end{equation} 
Lastly, the bounds on the slack variables
\begin{equation}
\begin{array}{l}\label{eq:Bounds}
    0 \le \Delta {v}_i \le \Delta {\bar{v}}  \quad \textrm{for}\quad i = 0, \dots, N-1 
\end{array}
\end{equation}
are added to the problem. This convexification step is referred to as lossless because it can be proved that the optimal solution of the convexified problem is also the optimal solution of the original one \cite{Liu2014}.  

\subsubsection{Linearization of the Dynamics}\label{dynLinear}
The introduction of slack variables is not sufficient to make the problem convex due to nonlinearities in Eq. \eqref{eq:SST}. These are dealt with by successive linearizations, requiring an iterative process. We will refer to these iterations as \textit{major iterations} associated with index $j$ in the remainder of this section.

Assume a solution $\mathbb{x}^{j-1}$ is available providing a reference trajectory (only at the first major iteration the reference trajectory is ballistic) about which the dynamics are linearized. At the $j$-th iteration, the linear part of Eq. \eqref{eq:SST} can be extracted, resulting in 
\begin{gather}\label{eq:Amatrix}
 \begin{bmatrix} \bvec{\delta r}_{i+1}^-  \\ \bvec{\delta v}_{i+1}^-  \end{bmatrix}^j
 =
  \begin{bmatrix}
   A_{\bvec{\delta r},i} &
   A_{\bvec{\delta v},i} 
   \end{bmatrix}^j
    \begin{bmatrix} \bvec{\delta r}_{i}^+  \\  \bvec{\delta v}_{i}^+  \end{bmatrix}^j = A_i^j     \begin{bmatrix} \bvec{\delta r}_{i}^+  \\  \bvec{\delta v}_{i}^+  \end{bmatrix}^j 
\end{gather}

The composition of all these linear maps results in
\begin{gather}\label{eq:Abbmatrix}
\begin{split}
 & \begin{bmatrix} \bvec{\delta r}_{N}  \\ \bvec{\delta v}_{N}  \end{bmatrix}^j = 
 \begin{bmatrix} \bvec{\delta r}_{p}  \\ \bvec{\delta v}_{p}  \end{bmatrix}^j
 =\\
 & \begin{bmatrix}
   A_{N-1}A_{N-2}\dots A_{\bvec{\delta v},0},  A_{N-1}A_{N-2}\dots A_{\bvec{\delta v},1},   \dots,  A_{\bvec{\delta v},N-1},  \boldsymbol{0}_{6\times N}
   \end{bmatrix}^j  \left(\mathbb{x}^j-\mathbb{x}^{j-1}\right) =\\ 
   & = \mathbb{A}^j \left(\mathbb{x}^j-\mathbb{x}^{j-1}\right),
   \end{split}
\end{gather}
where $(\mathbb{x}^j-\mathbb{x}^{j-1})$ result from the absolute impulses in the optimization vector $\mathbb{x}^j$ contains absolute impulses, whilst the linearizations are about the optimal impulses of the $(j-1)$-th iteration.  
The perturbed time of closest approach, $t^j_{CA}$, is   
\begin{equation}\label{eq:deltatca}
\begin{array}{l}
 t^j_{CA} = t^{j-1}_{CA} + \mathbb{B}^j\mathbb{A}^j\left(\mathbb{x}^j-\mathbb{x}^{j-1}\right) = \mathbb{B}^j\mathbb{A}^j\mathbb{x}^j + \left( t^{j-1}_{CA} - \mathbb{B}^j\mathbb{A}^j\mathbb{x}^{j-1}\right),\\
\end{array}
\end{equation}
where $\mathbb{B}$ is the $1\times6$ linear part of Eq. \eqref{eq:tcaManeuvered}. Similarly, the perturbed relative position vector on the b-plane, $\bvec{\Delta r}^j_{CA}$, is given by  
\begin{equation}\label{eq:deltarrb}
\begin{array}{l}
 \bvec{\Delta r}^j_{CA} = \bvec{\Delta r}^{j-1}_{CA} + \mathbb{C}^j\mathbb{A}^j\left(\mathbb{x}^j-\mathbb{x}^{j-1}\right) = \mathbb{C}^j\mathbb{A}^j\mathbb{x}^j + \left(\bvec{\Delta r}^{j-1}_{CA} - \mathbb{C}^j\mathbb{A}^j\mathbb{x}^{j-1}\right),\\
\end{array}
\end{equation}
in which $\mathbb{C}^j$ is the $2\times6$ linear part of Eq. \eqref{eq:conj1} at the $j$-th iteration. When $j=1$, all the $(j-1)$ quantities in the equations above are relative to the un-maneuvered case, i.e. $\bvec{\Delta r}^{0}_{CA} = \bvec{\Delta r}^{*}_{CA}$ and $\mathbb{x}^{0} = \mathbb{0}$.
Note that as CAMs result in a variation of the position vector of only few kilometers with respect to the ballistic trajectory (i.e., a relative variation of $0.1\%$), the convexification of the dynamics does not come with issues like artificial infeasibility and unboundedness. For this reason we did not introduce any artificial control or sophisticated  trust-region algorithms, but only problem-driven bounds on the impulse magnitudes and a maximum final state deviation on the b-plane. 

\subsubsection{Linearization of the Squared Mahalanobis Distance} \label{consLin}

The last step consists of dealing with the squared Mahalanobis distance constraint that defines an elliptically shaped avoidance region (or keep-out region using rendezvous terminology) on the b-plane. This non-convex constraint is dealt with by a second iterative process, nested in each major iteration, consisting of a projection and a linearization. We refer to these iterations as \textit{minor iterations}, with index $k$.

Assume the solution $\mathbb{x}^{j,k-1}$ of the $j$-th major iteration and $(k-1)$-th minor iteration is available. This defines the relative position vector on the b-plane $\bvec{\Delta r}^{j,k-1}_{CA}$. The $k$-th iteration starts by the projection algorithm that finds the point $\bvec{z}^{k}$ on the ellipse $\left({d_{CA}^2}\right)^j = \bar{d}_{CA}^2$ closest to $\bvec{\Delta r}^{j,k-1}_{CA}$. This is a convex optimization sub-problem with objective function
\begin{equation}\label{eq:minDistance}
    \min_{\bvec{z}}  ||\bvec{\Delta r}^{j,k-1}_{CA} - \bvec{z}|| 
\end{equation}
subject to the inequality constraint
\begin{equation}\label{eq:ellipseDistance}
    {\bvec{z}}^T (C_{CA}^j)^{-1} {\bvec{z}}^T \le  \bar{d}_{CA}^2,
\end{equation}
that can be solved efficiently using convex optimization algorithms (the interested reader can refer to \cite{He2018} for an exhaustive analysis of this sub-problem). $C_{CA}^j$ does not depend on $k$ as this quantity is assumed constant within each minor iteration.

Once $\bvec{z}^{k}$  is computed, the squared Mahalanobis distance constraint is linearized. The linear constraint ensures that $\bvec{\Delta r}^{j,k}_{CA}$ belongs to the half-plane tangent to the constraint in $\bvec{z}^{k}$, i.e.    
\begin{equation}
    \nabla {{{d}^2_{CA}}^j}(\bvec{z}^{k}) (\bvec{\Delta r}^{j,k}_{CA}-\bvec{z}^{k}) \ge 0.
\end{equation}
By substituting Eq. \eqref{eq:deltarrb}, the final expression is obtained 
\begin{equation}\label{eq:linMaha}
    -\nabla {{{d}^2_{CA}}^j}(\bvec{z}^{k})\mathbb{C}^j\mathbb{A}^j\mathbb{x}^{j,k} \le\nabla {{{d}^2_{CA}}^j}(\bvec{z}^{k}) (\bvec{\Delta r}^{j,k-1}_{CA}-\bvec{z}^{k} -\mathbb{C}^j\mathbb{A}^j\mathbb{x}^{j,k-1}),
\end{equation}
in which $\bvec{\Delta r}^{j,0}_{CA} = \bvec{\Delta r}^{j-1}_{CA}$ (i.e., the value at the end of the previous major loop). 

The selection of the first starting point determines the algorithm's convergence to a local optimum of the original, non-convex, problem. However, as will be illustrated in Sec. \ref{sec:MandmIterations}, for alert times on the order of two orbits, the problem appears to have only two local minima on opposite sides of the elliptical boundary of the avoidance region. Thus, two initial starting points determined by $\bvec{\Delta r}^{1,0}_{CA} = \pm \bvec{\Delta r}^{*}_{CA}$ are sufficient to automatically identify the global minimum. More guesses are needed for longer alert times, but this is not a practical issue thanks to the efficiency of the method and because \gls{cam}s are preferably executed close to the conjunctions \cite{Patera2003}. 

\subsubsection{Successive Convex Optimization Algorithm}\label{scvx}

With the introduction of slack variables, the linearization of the dynamics, and the linearization of the constraints, the multiple-impulse CAM design problem can be solved by primal-dual interior-point methods. We use MOSEK \cite{mosek} through its MATLAB interface. 

\begin{algorithm}
  \caption{Successive convexification optimization algorithm}\label{algo:SCVX}
  \begin{algorithmic}[1]
     \STATE Get inputs from \gls{cdm}: $R$, $\boldsymbol{r}^*_{p/s,CA}$ $\boldsymbol{v}^*_{p/s,CA}$, $t^*_{CA}$, $C^*_{p/s, CA}$;
     \STATE Assign $t_0$, $\Delta t$, $\Delta \bar{v}_i$ $N$, $\bar{P}_{C}$ or $\bar{P}_{C,\max}$ or $\bar{d}_{\min}$, $\textrm{tol}_{M}$, $\textrm{tol}_{m}$;
     \IF{$\bar{P}_{C,\max}$ is defined}
     \STATE Calculate $C^*_{CA}$ and $\bar{d}_{CA}^2$ by Eq. \eqref{eq:approxPC};
     \ELSIF{{$P_{C,\max}$ is defined}}
     \STATE Calculate $C^*_{CA}$ and $\bar{d}_{CA}^2$ by Eq. \eqref{eq:maxPc};
     \ELSE 
     \STATE Set $C^*_{CA} \gets \mathcal{I}$ and $\bar{d}_{CA}^2 \gets \bar{d}_{\min}^2$ 
     \ENDIF
     \STATE Back propagate the trajectories from $t^*_{CA}$ to $t_0$ and save $\boldsymbol{r}_{p/s,t_0}$, and $\boldsymbol{v}_{p/s,t_0}$; 
     \STATE Define the time grid $(t_0:\Delta t:t_0 + N\Delta t)$;  
     \STATE $j \gets 0$, $\mathbb{x}^0 \gets \mathbb{0}$, $t^0_{CA} \gets t^*_{CA}$, $\bvec{\Delta r}^{1,0}_{CA} \gets \pm \bvec{\Delta r}^{*}_{CA}$, $C^0_{CA} \gets C^*_{CA}$;
     \WHILE{$(j = 0)$ \OR $ ||\mathbb{x}^j-\mathbb{x}^{j-1}||_{\infty}\ge \textrm{tol}_{M}$}
      \STATE $j \gets j+1$;
      \STATE Perform a $1$-st order DA propagation of the trajectories at $(t_0:\Delta t:t_0 + N\Delta t)$ and $t^{j-1}_{CA}$ with impulses extracted from $\mathbb{x}^{j-1}$;
     \IF{$P_{C}$ is defined}
     \STATE Calculate $C_{CA}^j$ and then $\bar{d}_{CA}^2$ from Eq. \eqref{eq:approxPC};
     \ELSIF{{$P_{C,\max}$ is defined}}
     \STATE  Calculate $C_{CA}^j$ and then $\bar{d}_{CA}^2$ from Eq. \eqref{eq:maxPc};
     \ELSE
     \STATE Set $C^j_{CA} \gets \mathcal{I}$ and $\bar{d}_{CA}^2 \gets d_{\min}^2$;
     \ENDIF
      \STATE Assemble the matrices $\mathbb{A}^j$, $\mathbb{B}^j$, $\mathbb{C}^j$ as in Eq. \eqref{eq:Amatrix}-\eqref{eq:deltarrb};
      \STATE $k \gets 0$ and $\bvec{\Delta r}^{j,0}_{CA} \gets \bvec{\Delta r}^{j-1}_{CA}$;
          \WHILE{$(k=0)$ \OR $ ||\bvec{\Delta r}^{j,k}_{CA} - \bvec{\Delta r}^{j,k-1}_{CA} ||_2\ge \textrm{tol}_{m}$}
          \STATE $k \gets k+1$;
          \STATE Calculate $z_k$ by solving the convex optimization sub-problem \eqref{eq:minDistance}-\eqref{eq:ellipseDistance};
          \STATE Calculate $\mathbb{x}^{j,k}$ by solving the convex optimization problem defined by \eqref{eq:obj}-\eqref{eq:Bounds} and \eqref{eq:linMaha}; 
          \STATE  $t^{j,k}_{CA} \gets \mathbb{B}^j\mathbb{A}^j\mathbb{x}^{j,k} + \left( t^{j-1}_{CA} - \mathbb{B}^j\mathbb{A}^j\mathbb{x}^{j-1}\right)$;
          \STATE $\bvec{\Delta r}^{j,k}_{CA} \gets \mathbb{C}^j\mathbb{A}^j\mathbb{x}^{j,k} + \left(\bvec{\Delta r}^{j-1}_{CA} - \mathbb{C}^j\mathbb{A}^j\mathbb{x}^{j-1}\right)$;
     \ENDWHILE
     \STATE $\mathbb{x}^{j} \gets \mathbb{x}^{j,k}$, $t_{CA}^j \gets t_{CA}^{j,k} $, and $\bvec{\Delta r}^{j}_{CA} \gets \bvec{\Delta r}^{j,k}_{CA}$;
     \ENDWHILE
   \end{algorithmic}
\end{algorithm}

The Algorithm \ref{algo:SCVX} provides a complete overview of the solution process. As previously described, two iterations are needed: a major $j$-th iteration for the linearization of the dynamics and a minor $k$-th iteration for linearization of the constraint. The two iterations stop when two conditions are met: the minor one ends when $||\bvec{\Delta r}^{j,k}_{CA} - \bvec{\Delta r}^{j,k-1}_{CA} ||_2 \le \textrm{tol}_m$, while the major one when $ ||\mathbb{x}^j-\mathbb{x}^{j-1}||_{\infty}\le \textrm{tol}_M$. The algorithm minimizes the total $\Delta v$ while constraining either the risk, the maximum risk, or the closest approach's distance. The latter case is dealt with by setting $C_{CA}^* = \mathcal{I}$ in all iterations.

\section{Test Cases}\label{NumRef}
\begin{figure}[h!]
\centering
\includegraphics[width=0.5\textwidth]{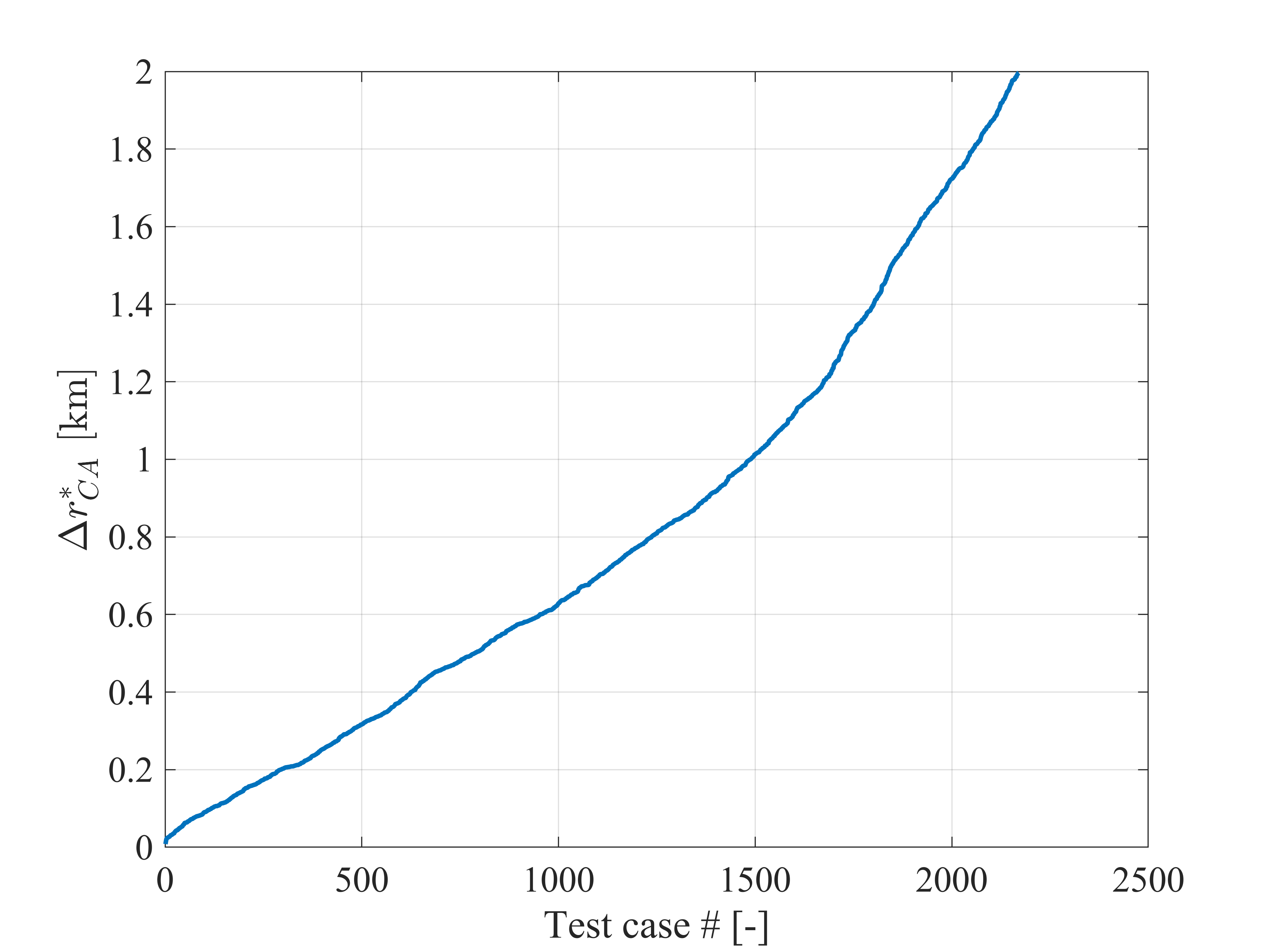}  
\caption{Distribution of minimum distance at closest approach}
\label{miss_distanceminDist}
\end{figure}

\begin{figure}[h!]
\centering
\includegraphics[width=0.5\textwidth]{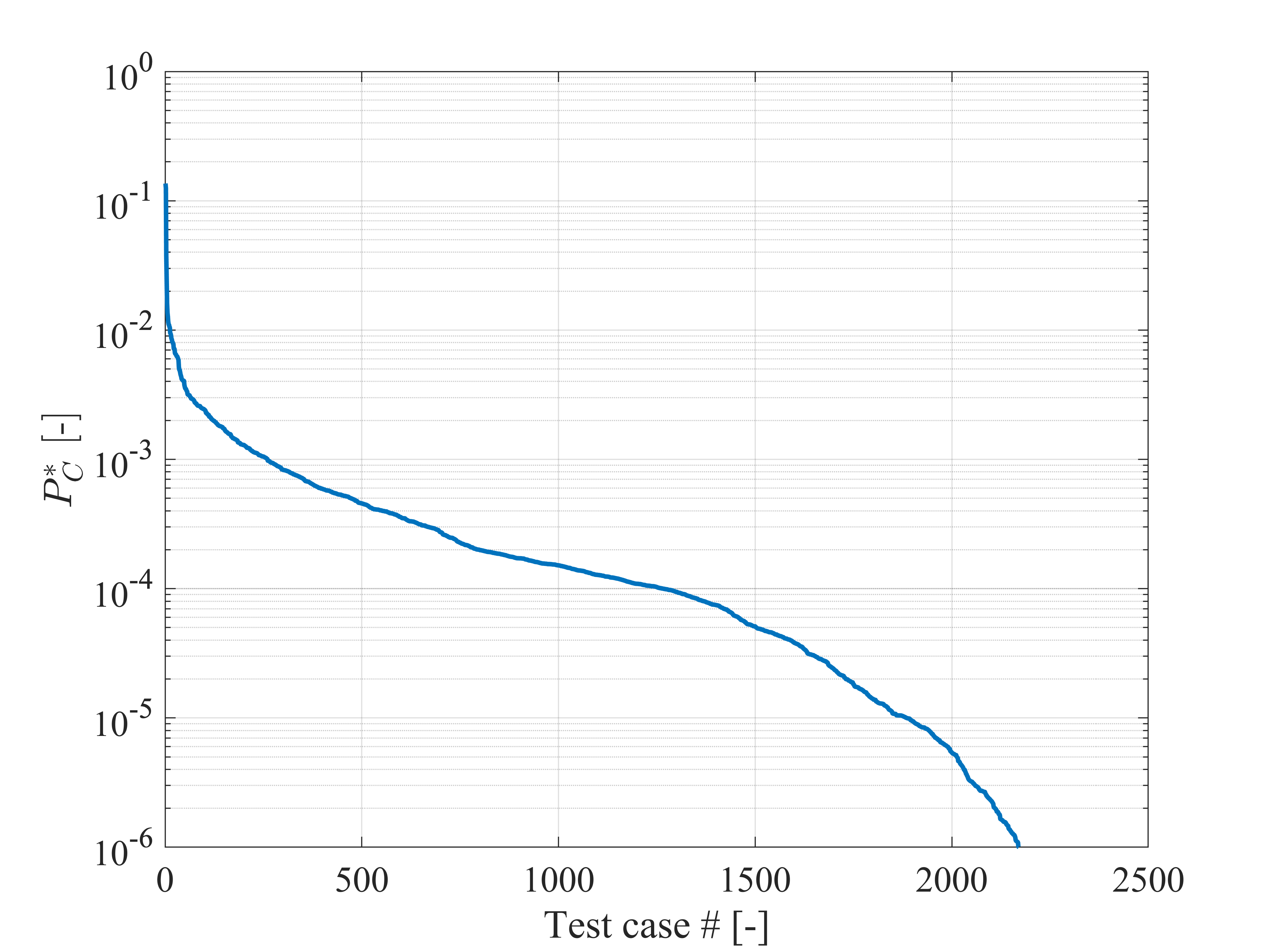} 
\caption{Distribution of collision probabilities}
 \label{miss_distanceminPcData}
\end{figure}

\begin{figure}[h!]
\centering
\includegraphics[width=0.5\textwidth]{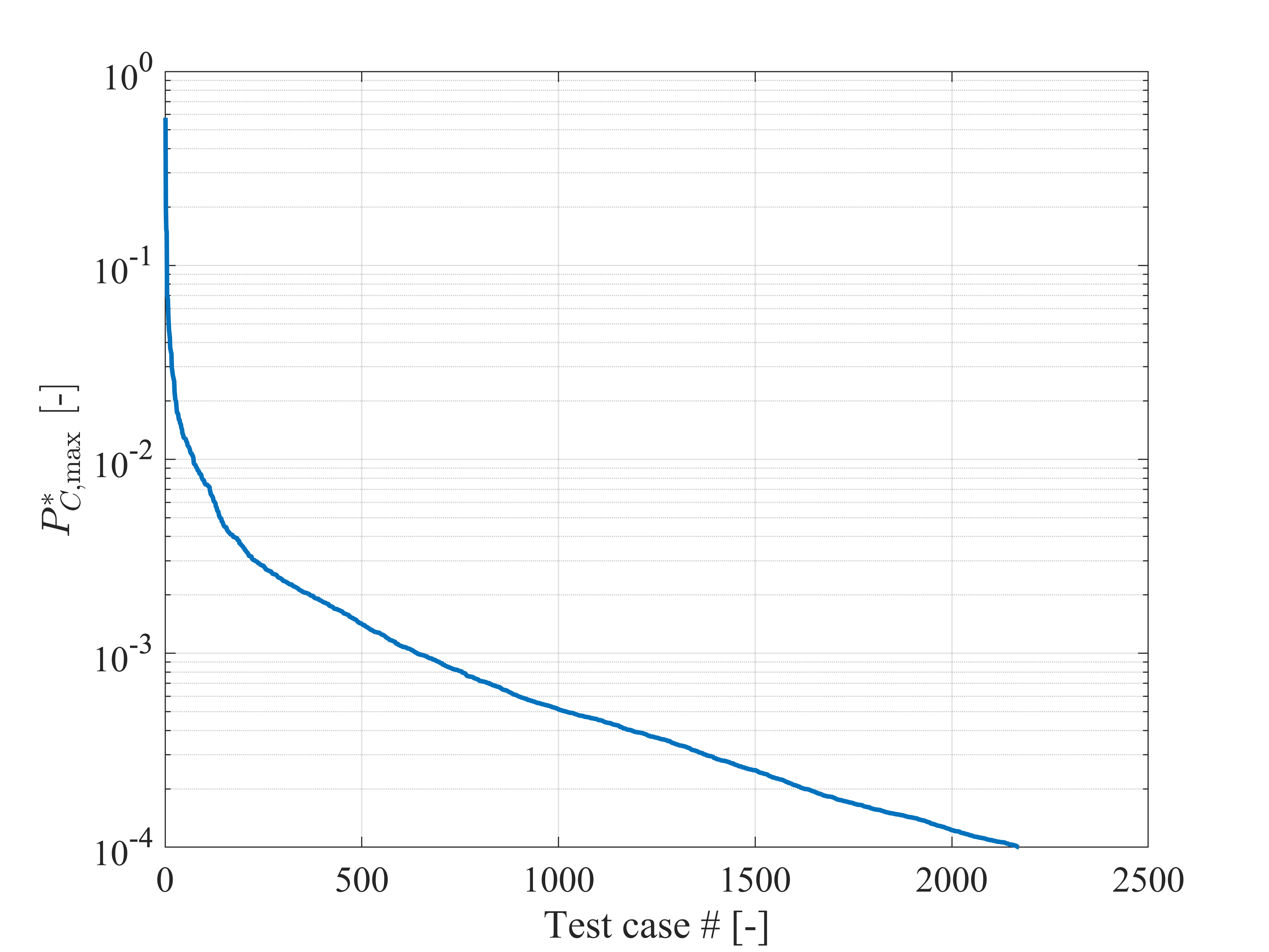}
\caption{Distribution of maximum collision probabilities}
\label{miss_distancemaxPcData}
\end{figure}

The methodology described in the previous sections is applied to test cases derived from the ESA Collision Avoidance Challenge. For this competition, ESA provided the teams with real conjunction data extracted from 162,634 \gls{cdm}s, corresponding to 13,154 unique events. However, ESA did not distribute the full orbital elements set and provided the positional covariances in the primary \gls{rtn}  reference frame only. With a procedure omitted here, we managed to reconstruct the full orbital data of the objects except for the \gls{raan}, for which only a relative estimate was computed (i.e., we set the \gls{raan} of one object to zero and solved for the other to obtain the prescribed conjunction). These data were filtered to consider conjunctions with $\Delta r_{CA}^* \le 2$ km, $P_C^* > 10^{-6}$, and $P_{C,\max}^* > 10^{-4}$, resulting in a new data file with 2,170 conjunctions, available for download at \url{github.com/arma1978/conjunction}. The distribution of the minimum distance at \gls{ca}, the collision probability and maximum collision probability are shown in Fig. \ref{miss_distanceminDist}--\ref{miss_distancemaxPcData}. All the conjunctions are relative to objects in \gls{leo} with 90$\%$ of cases having relative conjunction speeds $\in [1.80, 14.98]$ km/s.

All the simulations presented in the next sections were obtained with a dynamical model including $J_2-J_4$ zonal harmonics
\begin{equation}\label{EoM}
\left\{
	\begin{array}{l}
		\dot{x} = v_x \\ 
		\dot{y} = v_y \\ 
		\dot{z} = v_z\\
		\dot{v}_x = -\frac{\mu x}{r^3} + \frac{3\mu J_2 R_e^2}{2 r^5}(\frac{5z^2}{r^2}-1)x + \frac{5\mu J_3 R_e^3 x z}{2r^7}(\frac{7z^2}{r^2}-3) +\frac{15\mu J_4 R_c^4 x}{8r^7}(1-\frac{14z^2}{r^2}+\frac{21z^4}{r^4}) \\
        \dot{v}_y = -\frac{\mu y}{r^3} + \frac{3\mu J_2 R_e^2}{2 r^5}(\frac{5z^2}{r^2}-1)y + \frac{5\mu J_3 R_e^3 y z}{2r^7}(\frac{7z^2}{r^2}-3) +\frac{15\mu J_4 R_e^4 y}{8r^7}(1-\frac{14z^2}{r^2}+\frac{21z^4}{r^4})\\
\dot{v}_z = -\frac{\mu z}{r^3} + \frac{3\mu J_2 R_e^2}{2 r^5}(\frac{5z^2}{r^2}-3)z + \frac{5\mu J_3 R_e^3}{2r^5}(\frac{3}{5}-\frac{6z^2}{r^2} +\frac{7z^4}{r^4}) +\frac{15\mu J_4 R_e^4 z}{8r^7}(5-\frac{70z^2}{3r^2}+\frac{21z^4}{r^4})
	\end{array}
	\right.
\end{equation}
in which $\bbb{r} = [x,y,z]^T$ and $\bbb{v} = [v_x,v_y,v_z]^T$ are the spacecraft position and velocity vectors; $\mu$, $R_e$, and $J_i$ are the gravitational parameter, the mean equatorial radius, and the $i$-th zonal harmonic coefficient of the Earth. 
Any dynamical model can be selected without affecting the algorithm complexity as the required state transition matrices are automatically obtained with \gls{da} without the need to derive and integrate the variational equations.   
The model in Eq. \eqref{EoM} was used to demonstrate that the method is agnostic to the dynamical model used, as the state transition matrices are automatically obtained with \gls{da} without the need to derive and integrate the variational equations. In Sec. \ref{sec:extensiveSim} it will be shown that the perturbation effect is negligible for the vast majority of the cases presented in this work, this feature is relevant when extending to strongly non Keplerian cases.   

In Sec. \ref{sec:MandmIterations}--\ref{sec:ecc} we offer a detailed analysis of the reference scenario, which is the case with the highest collision probability in the dataset (details are found in the Appendix). In Sec. \ref{sec:extensiveSim} we provide a summary of the method's performance when applied to the entire dataset. Unless specified, all simulations assume $\Delta t = 1$ min, a maximum impulse of 6 mm/s (corresponding to a constant thrust  of 30 mN for a reference spacecraft of 300 kg), tol$_M = 1$ mm/s, and tol$_m = 1$ m. The one-minute discretization is deemed appropriate to approximate finite burns with a sequence of impulses, even in the case of eccentric orbits. Although not implemented here, a uniform discretization in an angular variable, such as the true anomaly, to avoid oversampling would, in general, be preferred. 

The simulations are run on a MacBook Pro with a 2,3 GHz Quad-Core Intel Core i7 and 16 GB Memory.

\subsection{Major and Minor Iterations}\label{sec:MandmIterations}
\begin{figure}[h!]
\subfigure[Convergence summary]{
\includegraphics[width=0.47\textwidth]{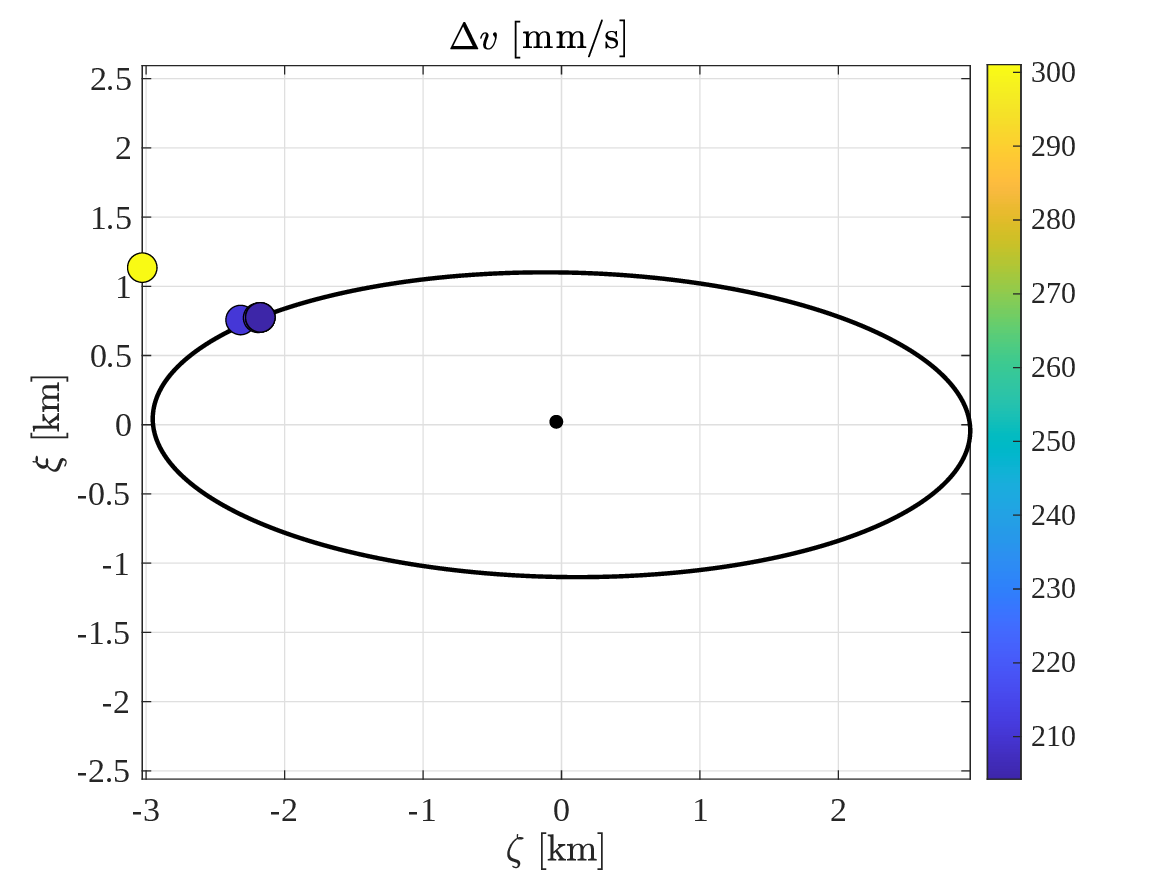}
\label{exBplane}}\hfill
\subfigure[First and second iterations]{
\includegraphics[width=0.47\textwidth]{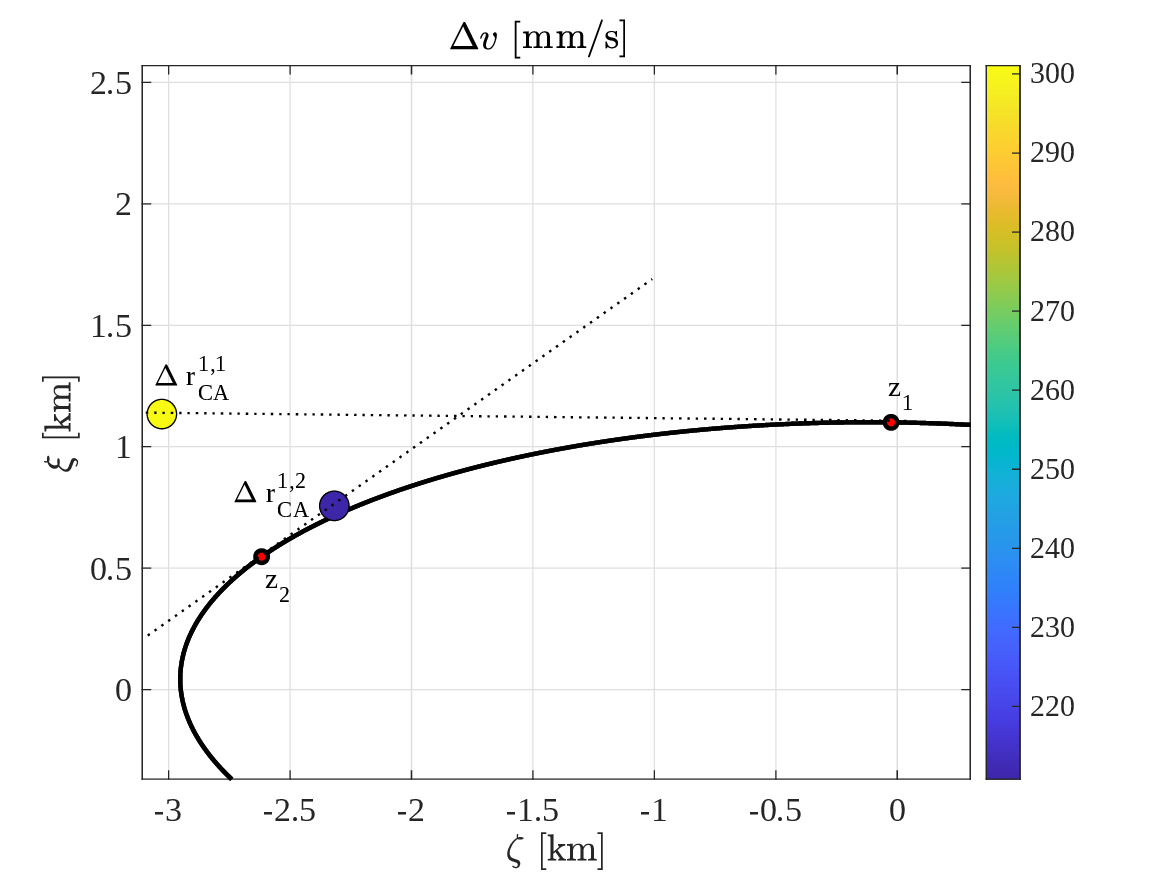}
\label{exConv12}} \hfill
\subfigure[Third and fourth iterations]{
\includegraphics[width=0.47\textwidth]{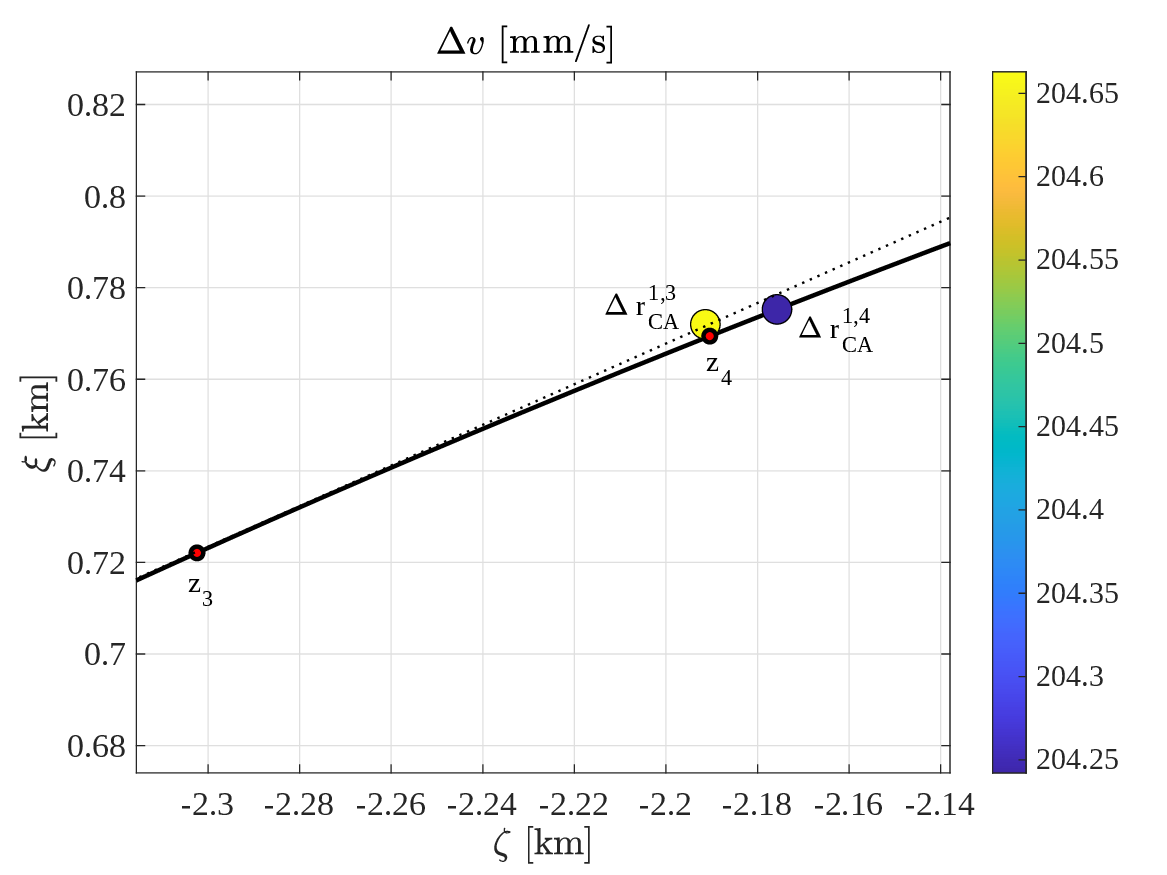}
\label{exConv34}}\hfill
\subfigure[Fifth iteration and first iteration of second major iteration]{
\includegraphics[width=0.47\textwidth]{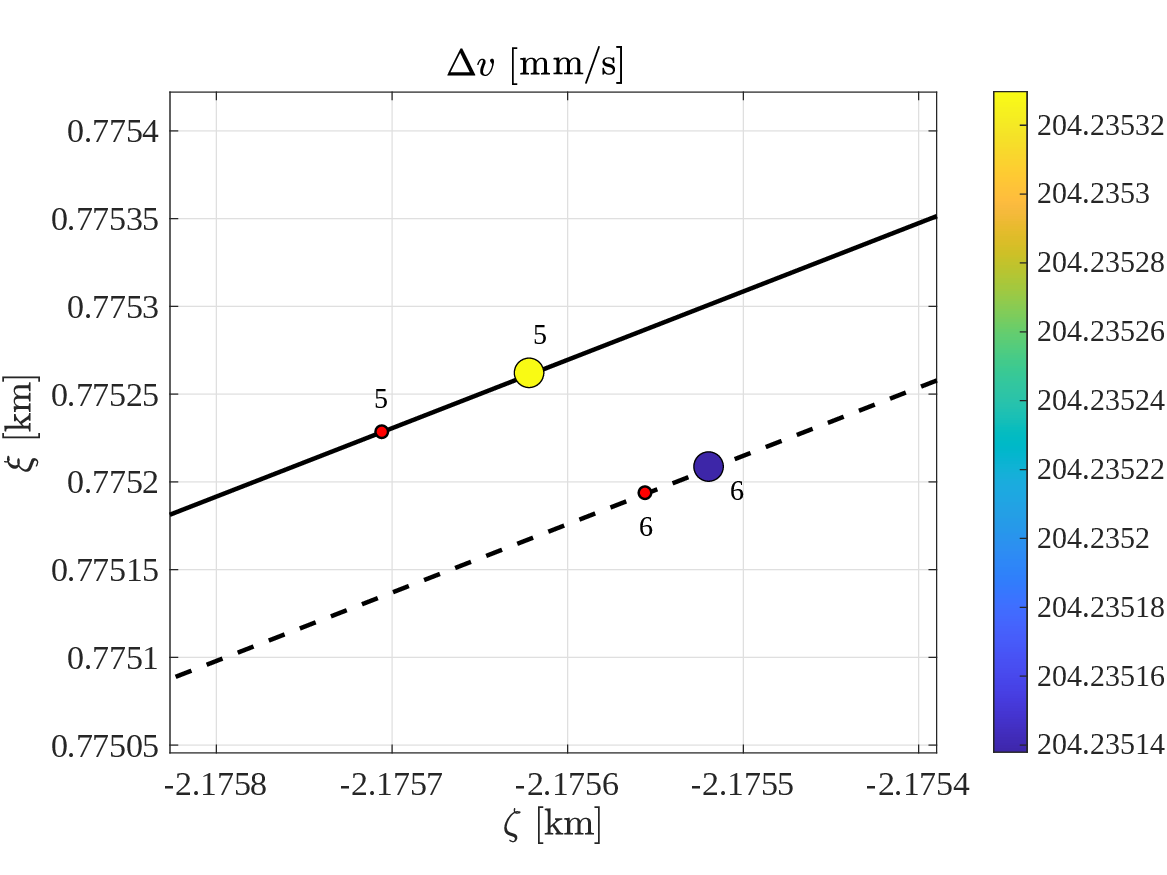}
\label{exConv56}}\hfill
\caption{B-plane convergence analysis}
\label{exConv}
\end{figure}

Figure \ref{exConv} provides details of the convergence of Algorithm \ref{algo:SCVX} when the \gls{cam} can be applied between 8 and 6 orbits before the \gls{ca}, with a maximum of 200 impulses and a target $\bar{P}_{C,\max} = {10}^{-4}$. In Fig. \ref{exBplane} the contour line $P_{C,\max} = \bar{P}_{C,\max}$ is plotted as a black solid line, the dot close to the origin is the unperturbed primary position with $\Delta r^*_{CA} = 43$ m, while the circles represent the perturbed primary position during the algorithm iterations colored according to the maneuver $\Delta v$. The optimization requires two major iterations, with 5 and 1 minor iterations, respectively. In Fig. \ref{exConv12} the first minor iteration starts from the red circle labeled $1$ representing $\bvec{z}_1$, which is the point on the constraint line closest to the unperturbed solution. The output of the first optimizer run delivers the solution $\bvec{\Delta r}_{CA}^{1,1}$ indicated as a yellow filled circle, which belongs to the ellipse tangent in $\bvec{z}_1$. From $\bvec{\Delta r}_{CA}^{1,1}$ the new $\bvec{z}_2$ is calculated allowing for the updated solution $\bvec{\Delta r}_{CA}^{1,2}$ with a significant reduction of maneuver $\Delta v$. Figure  \ref{exConv34} shows similar information for the third and fourth iterations and \ref{exConv56} for the fifth iteration of the first major loop and the first and only iteration of the second major loop. It is worth highlighting that the contour line of equal maximum collision probability is slightly changed (dashed line) because, between the first and second major iteration, the b-plane and the combined 2D covariance are updated. However, note that between the first and the second major loop the $\Delta v$ changes only by $3\times 10^{-4}$ mm/s, showing that linearized dynamics and transformation are particularly accurate for thanks also to a relative conjunction velocity of almost 15 km/s.

\begin{figure}[h!]
\subfigure[Thrust arcs on primary trajectory]{
\includegraphics[width=0.47\textwidth]{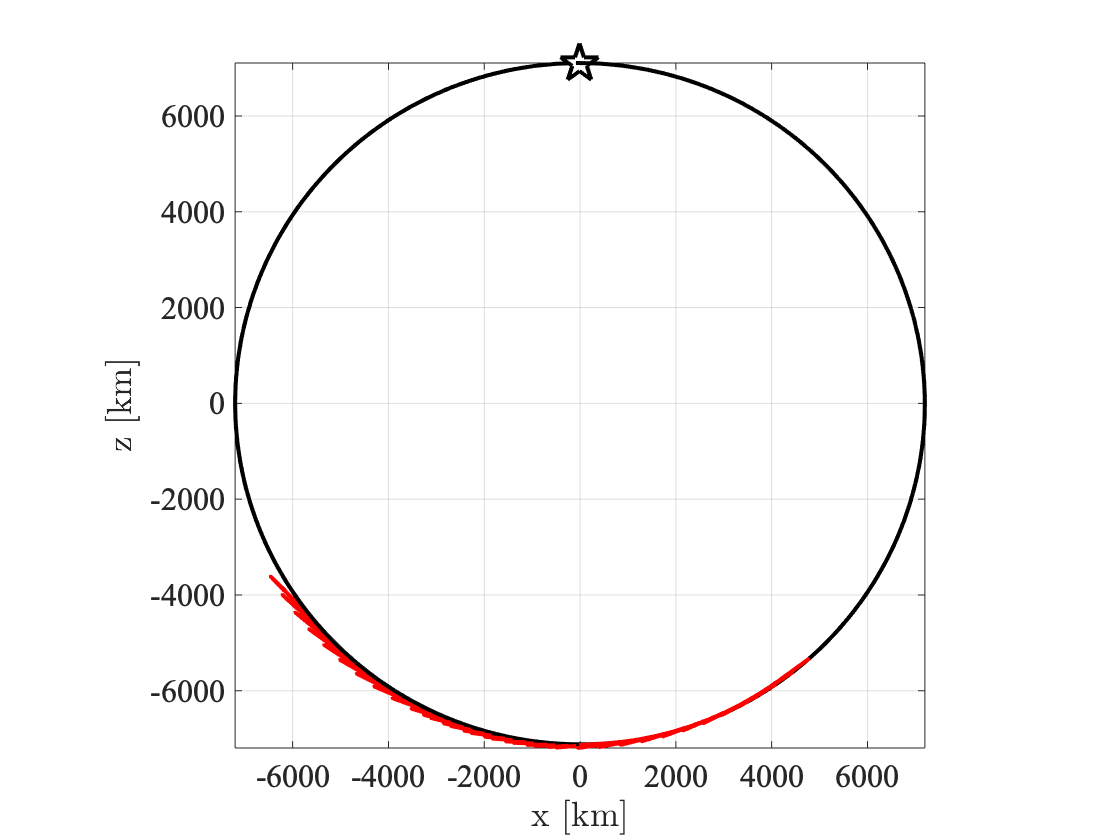}
\label{exTrj}}\hfill
\subfigure[Impulses timing and components]{
\includegraphics[width=0.47\textwidth]{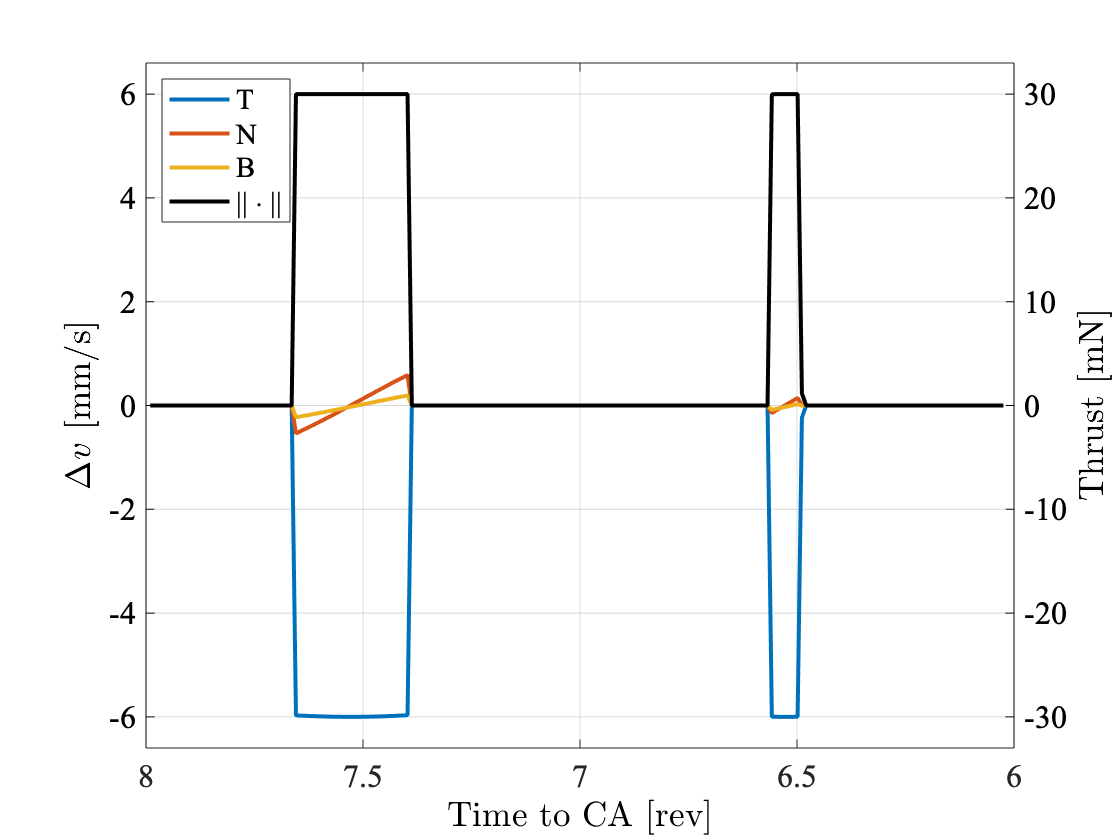}
\label{exDvComp}}
\caption{Example of trajectory and impulses profile}
\label{exFig}
\end{figure}

Figure \ref{exTrj} highlights the thrust arcs on the primary trajectory (the pentagram indicates the conjunction point) and Fig. \ref{exDvComp} shows impulses in the Frenet frame; where T, N and B indicate the tangential, normal and binormal directions respectively. Note that in the test cases, with the exception of Sec. \ref{sec:ecc}, the normal direction is almost aligned with the negative radial direction due to the limited eccentricity of the primary orbits. The thrust arcs' optimal location is spread along the portion of the orbit opposite to the conjunction. The total $\Delta v$ is $204.2$ mm/s, requiring 34 impulses.

Figure \ref{exDvComp} shows that the optimal maneuver has components in the radial and, to a lesser extent, out of plane directions. Remarkably, each minor iteration of the problem is solved in approximately ~12 ms and 14 iterations, allowing for a complete solution, in around 0.15 s. This computational time requires solving a large optimization problem entailing 800 optimization variables, 200 second-order cone constraints, one inequality constraint on the squared Mahalanobis distance, and simple bounds. By contrast, the \gls{nlp} method described in Sec. \ref{algo:NLP} requires, when it reaches convergence, thousands of iterations and a few minutes to achieve a comparable solution starting from a random feasible guess.

\begin{figure}[h!]
\centering
\includegraphics[width=0.5\textwidth]{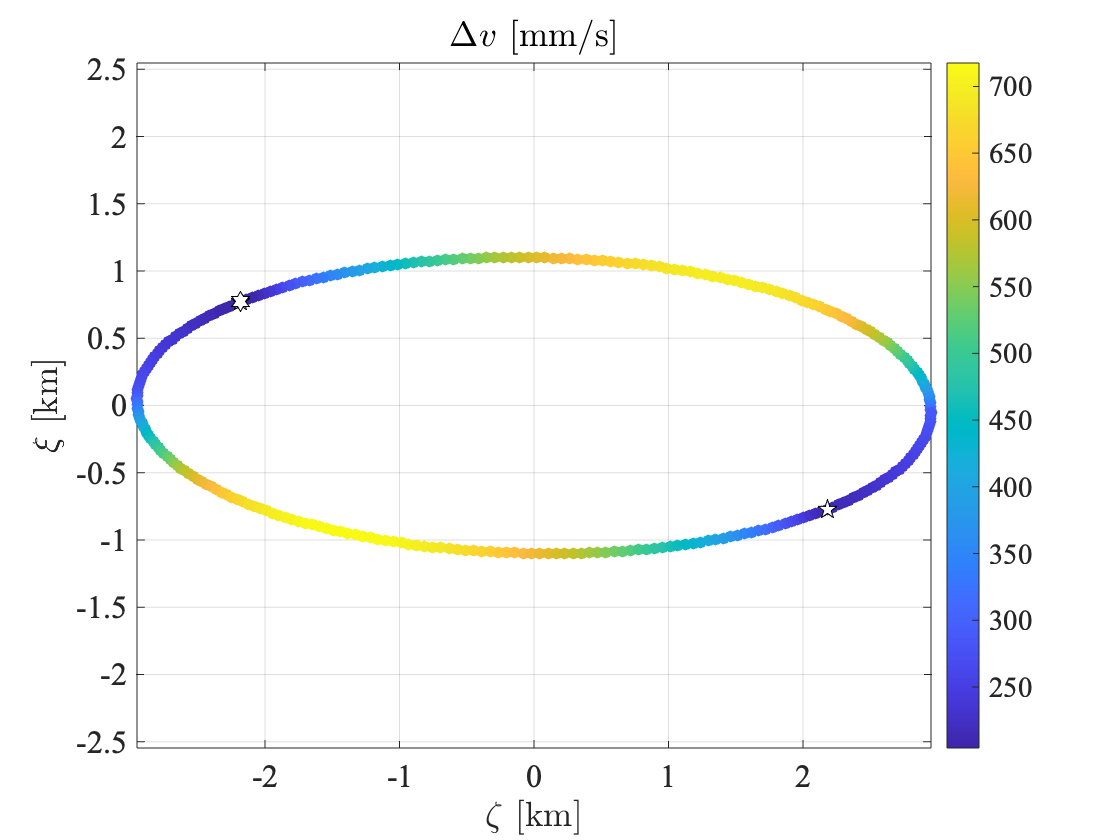}
\caption{Objective function profile on the boundaries of the avoidance region}
\label{EllipseDv}
\end{figure}

As mentioned in Sec. \ref{consLin}, it is not guaranteed that the minor iterations converge to the original non-convex problem's global optimum. Convergence to a local minimum depends on the selected starting point for the minor iteration. However, as the optimal solution lies on the admissible region's elliptical boundary, we can investigate the objective function's behavior on this boundary. This analysis is done by sampling the ellipse with $M$ points and solving $M$ convex optimization problems in which a terminal equality constraint substitutes the keep-out zone one. Figure \ref{EllipseDv} shows the results for $M = 300$, where the pentagram indicates a local minimum  (213.9 mm/s) and the hexagram the global one (204.2 mm/s) analyzed in the previous figures.  For short alert times or a limited duration of the maneuver, it was observed that the objective function has two minima located on opposite sides of the ellipse and corresponding to thrust mainly aligned with either the tangential or the anti-tangential direction. In these cases, the two minima can be identified by starting the minor iterations with two points on the ellipse with opposite coordinates. 

\subsection{Effect of Constraints} \label{sec:constraints}
\begin{figure}[h!]
\subfigure[Different constraints and corresponding solutions indicated with different markers]{
\includegraphics[width=0.47\textwidth]{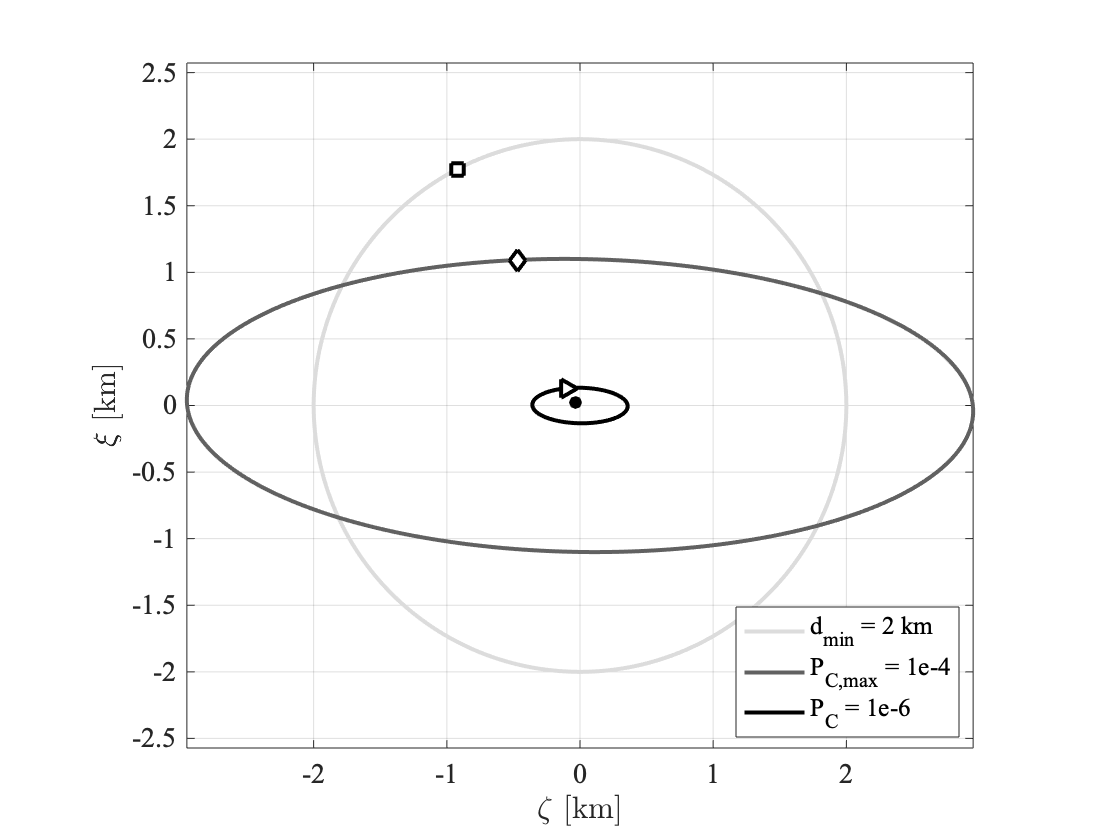}
\label{effectOfObj}}\hfill
\subfigure[Different $\bar{P}_{C,\max}$ contour levels and corresponding solutions indicated with different markers]{
\includegraphics[width=0.47\textwidth]{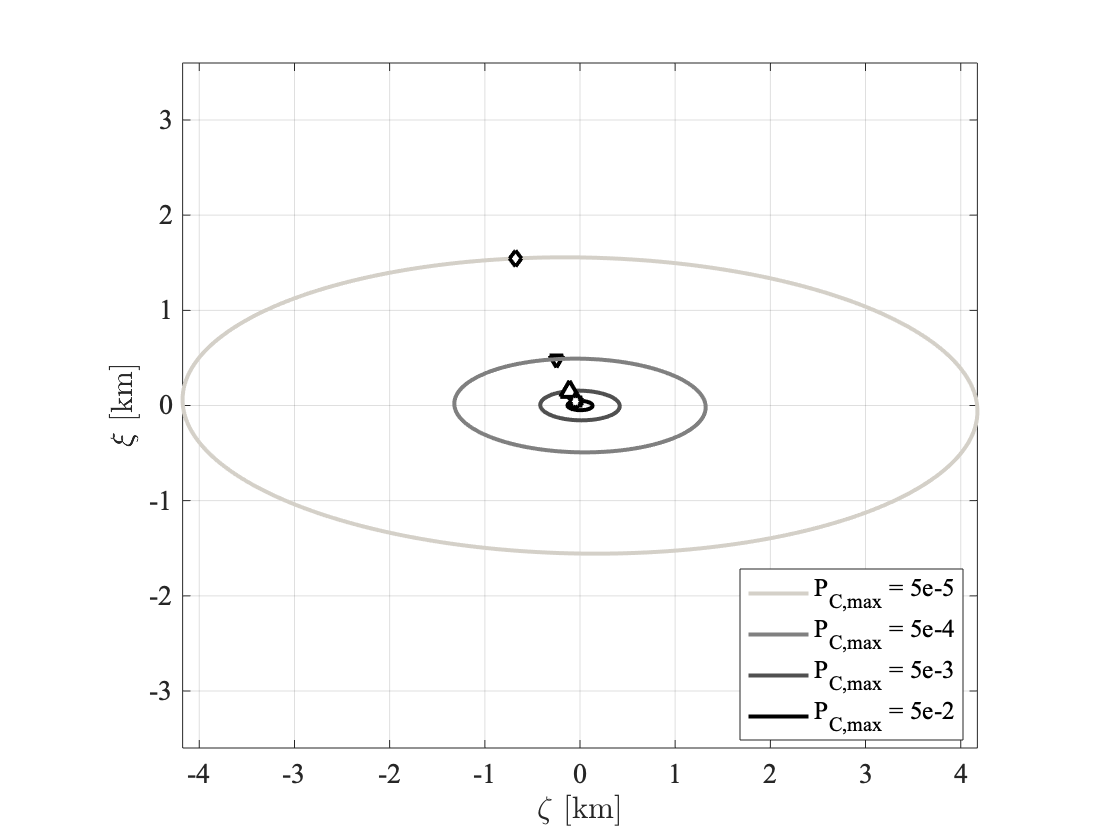}
\label{effectOfObjValue}}
\caption{B-plane analysis according to different constraints}
\label{EffectOfObj}
\end{figure}

The effect of the different constraints on the \gls{cam} design is analyzed, considering impulses on the last two orbits before the conjunction. Figure \ref{effectOfObj} shows the contour lines corresponding to $\bar{d}_{min} = 2$ km, $\bar{P}_{C,\max} = 10^{-4}$ (the case of the previous section), and $\bar{P}_{C} = 10^{-6}$. The collision probability constraint is the less restrictive one, corresponding to a $\Delta v = 28.1$ mm/s and five impulses. The maximum collision probability requires $\Delta v = 288.1$ mm/s and 48 impulses, slightly higher than the previous section's case as it closer to the conjunction. However, these two constraints are characterized by the same elliptical shape on the b-plane with different sizes. The minimum distance constraint is a circle, and the chosen value of $2$ km results in the most constraining one, resulting in $\Delta v = 527.4$ mm/s and 88 impulses. Changing the constraint value makes it possible to identify a line of optimal target points on the b-plane. An example is provided in Fig. \ref{effectOfObjValue} where different values of $\bar{P}_{C,\max}$ are considered, resulting in $\Delta v \in [6.4, 431.1]$ mm/s. A trade-off between safety and propellant consumption is enabled by the rapid calculation of the \gls{cam}.        

\subsection{Time to Closest Approach}\label{sec:time2CA}
\begin{figure}
\subfigure[18-16 orbits to \gls{ca}]{
\includegraphics[width=0.47\textwidth]{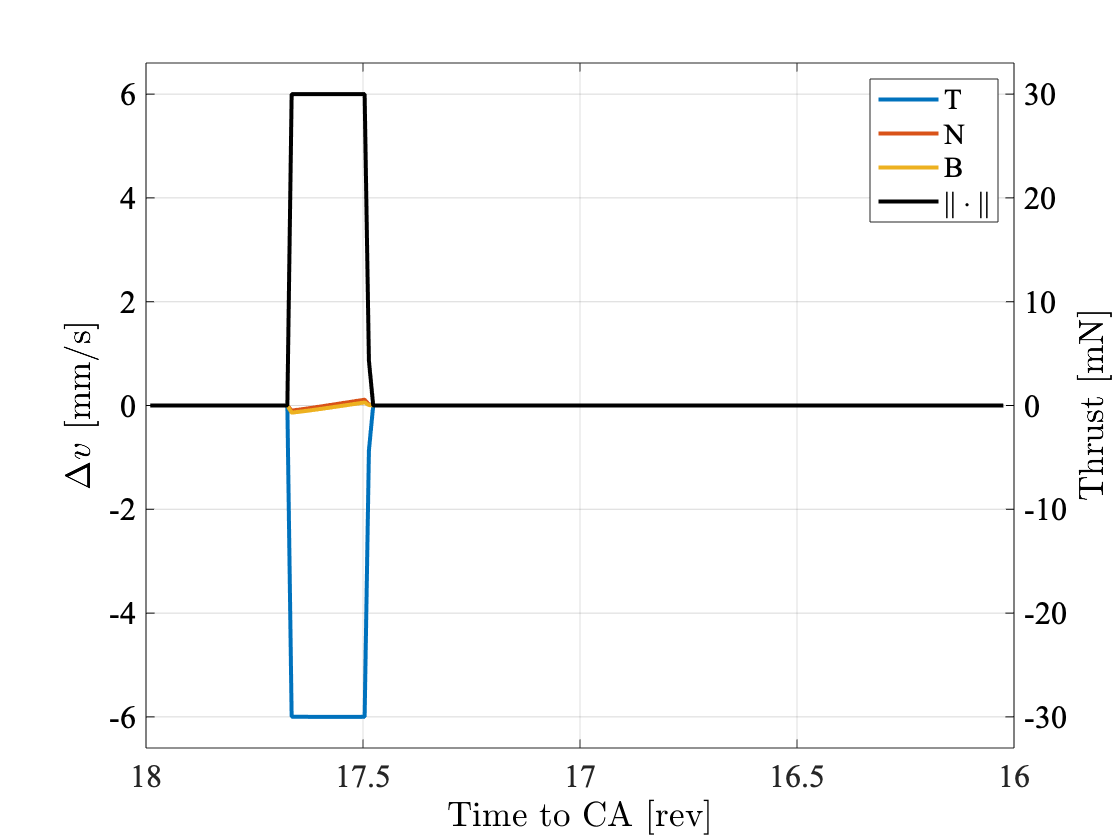}
\label{effecOfTime16}}\hfill
\subfigure[12-10 orbits to \gls{ca}]{
\includegraphics[width=0.47\textwidth]{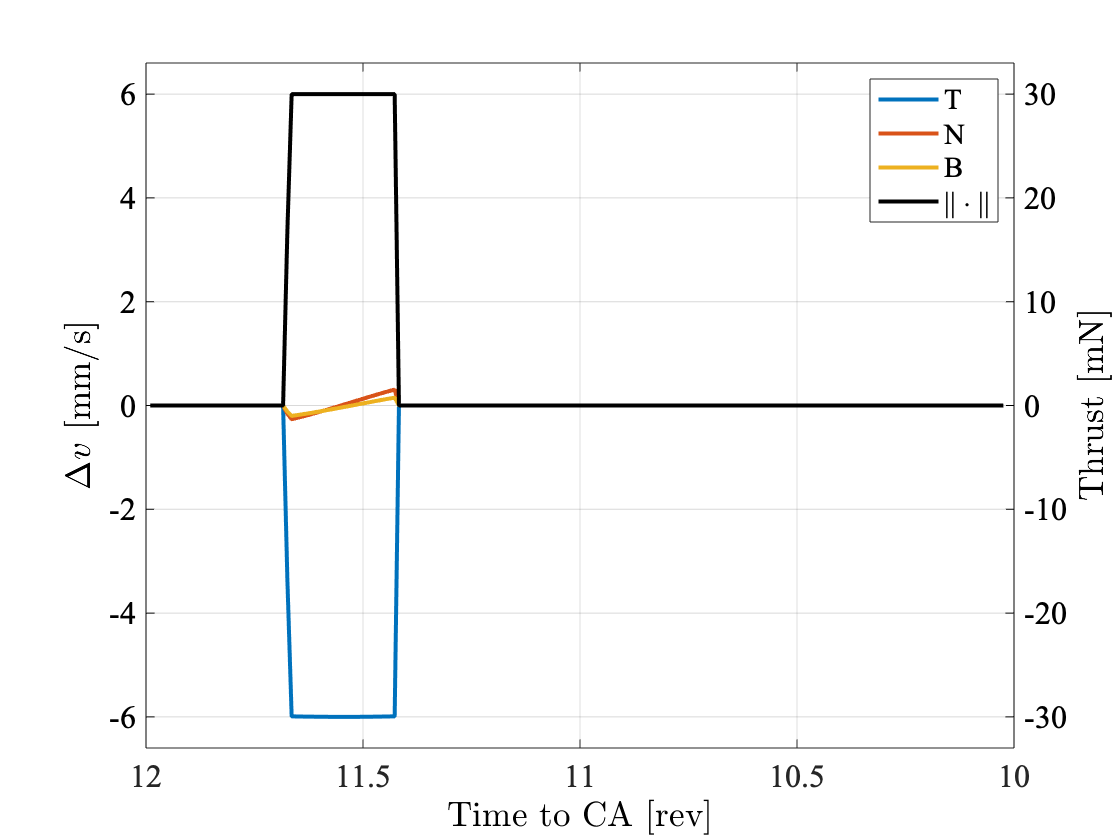}
\label{effecOfTime12}}\hfill
\subfigure[8-6 orbits to \gls{ca}]{
\includegraphics[width=0.47\textwidth]{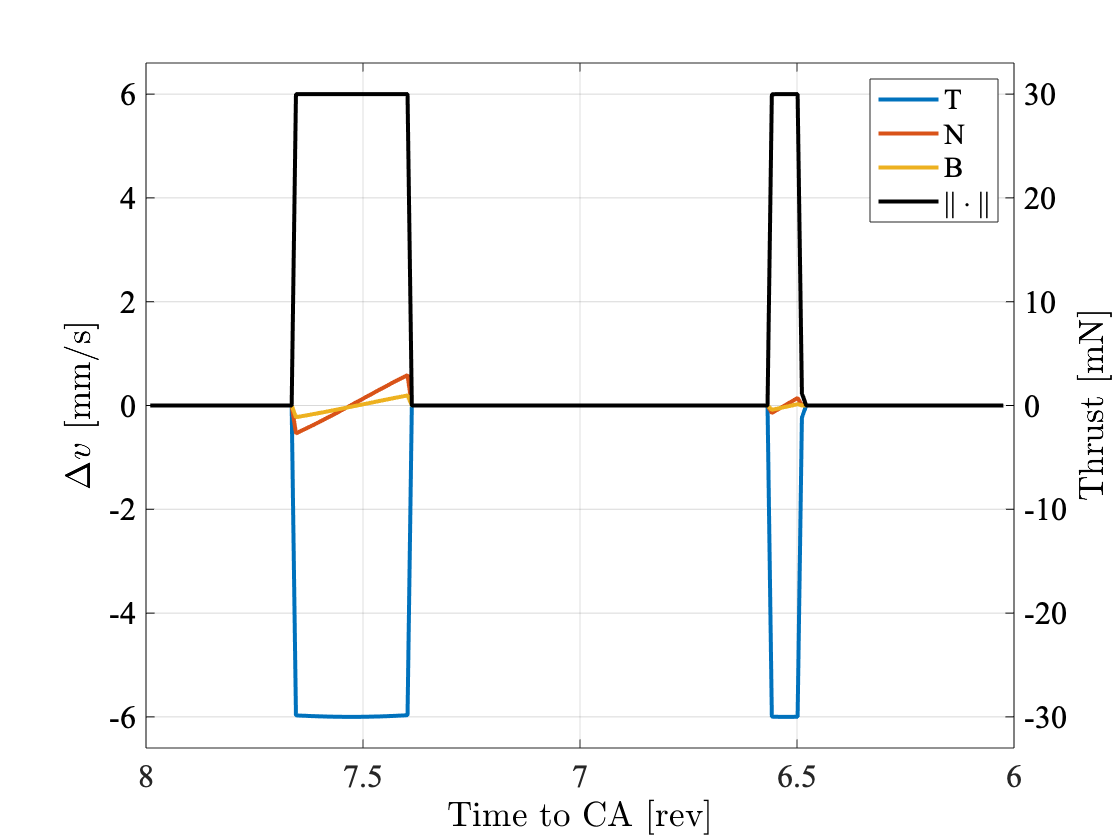}
\label{effecOfTime8}}\hfill
\subfigure[4-2 orbits to \gls{ca}]{
\includegraphics[width=0.47\textwidth]{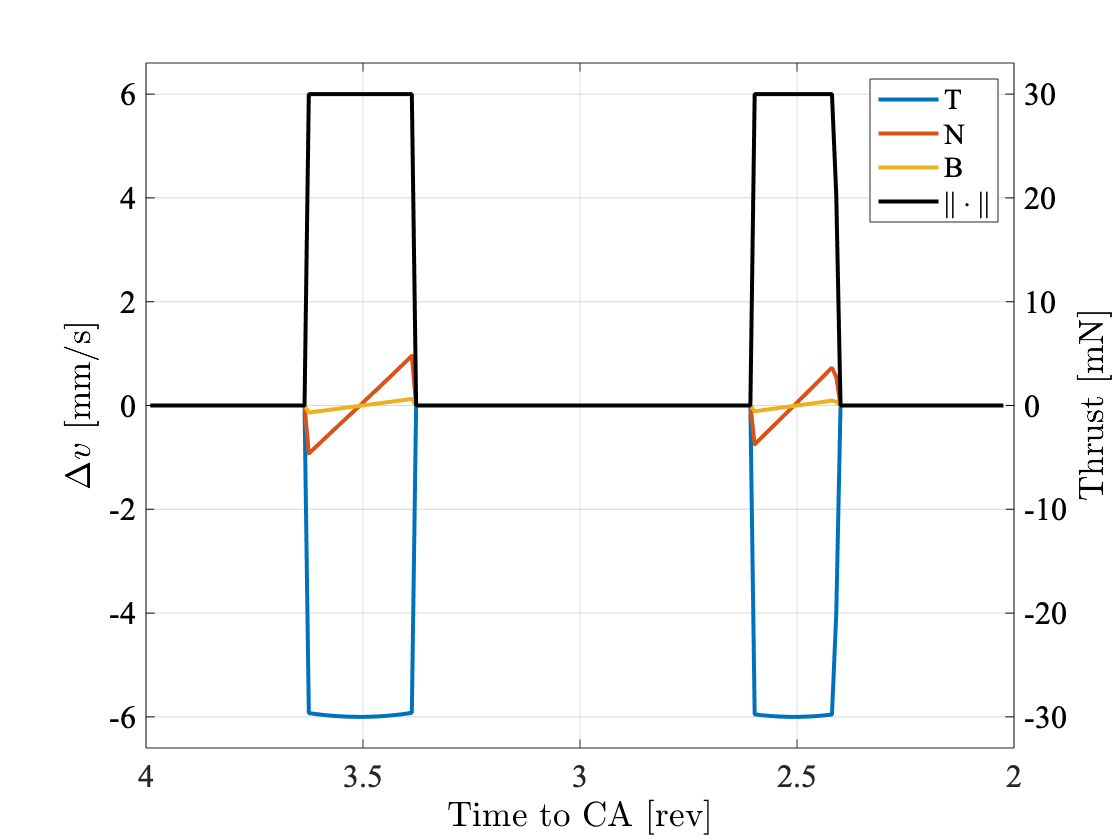}
\label{effecOfTime4}}\hfill
\caption{Impulses profile as a function of time to \gls{ca}}
\label{effectOftime}
\end{figure}

 We analyze alert's time impact on the \gls{cam}. In Fig. \ref{effectOftime} four graphs with alert time decreasing from 18 to 4 orbits are reported, with a maximum of 200 impulses. The $\Delta v$ is 108.9  mm/s for the first case, 153.4 mm/s for the second, 204.2 mm/s for the third and 268.1 mm/s for the last. The optimal maneuver in this case is performed at the first opposition. When the leading time reduces and the $\Delta v$ increases, the optimal maneuver is split into multiple arcs, two in this case. Additionally, the optimal maneuver always has a radial component. Furthermore, a maneuver applied several orbits before the conjunctions has a larger impact on the conjunction geometry, affecting the time to \gls{ca} by more than 2 seconds.   
 
\begin{figure}
\subfigure[Objective function profile on the boundaries of the avoidance region]{
\includegraphics[width=0.47\textwidth]{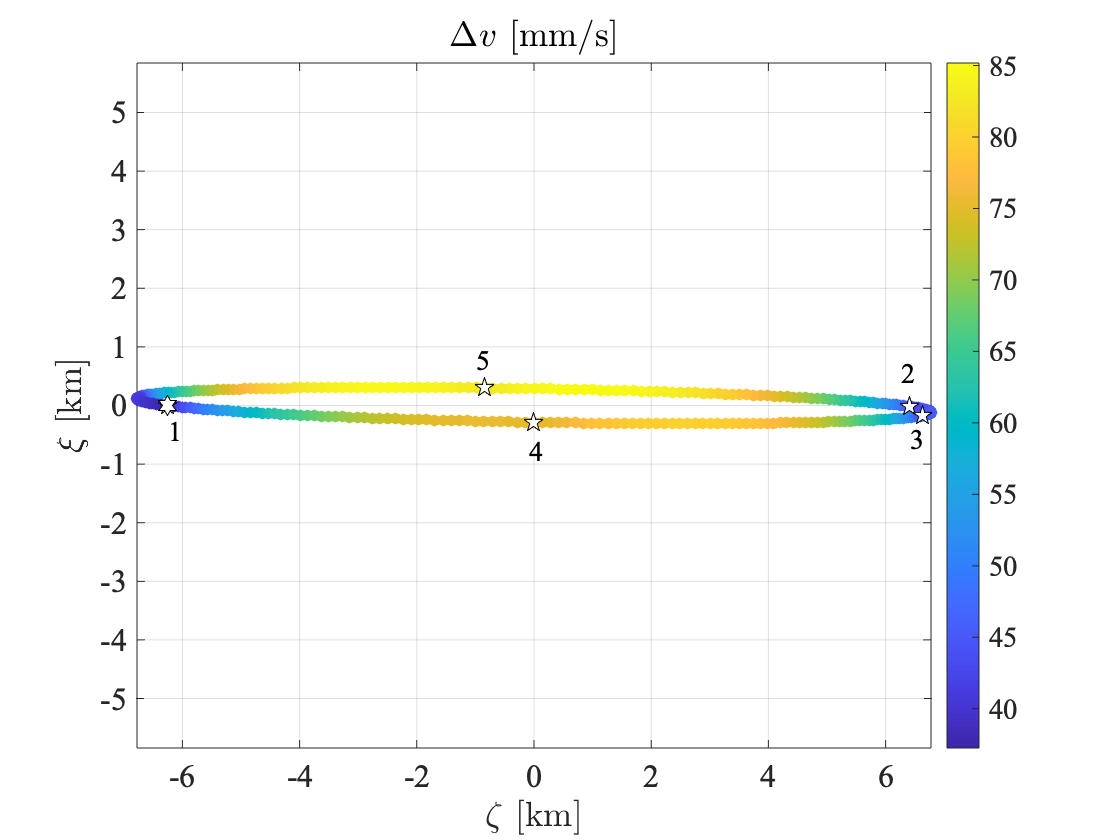}
\label{ellipseDvAnomalous}}\hfill
\subfigure[Solution \#1, global minimum $\Delta v = 37.2$ mm/s]{
\includegraphics[width=0.47\textwidth]{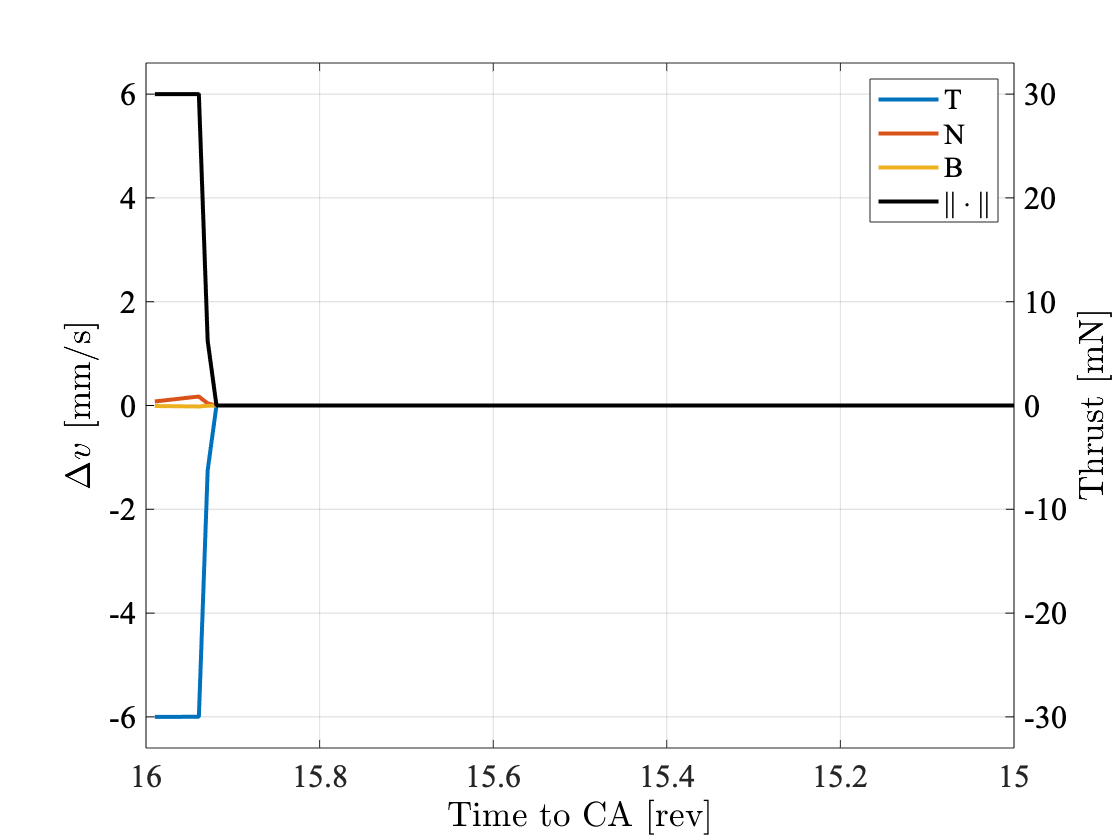}
\label{1}}\hfill
\subfigure[Solution \#2, local minimum $\Delta v = 42.6$ mm/s]{
\includegraphics[width=0.47\textwidth]{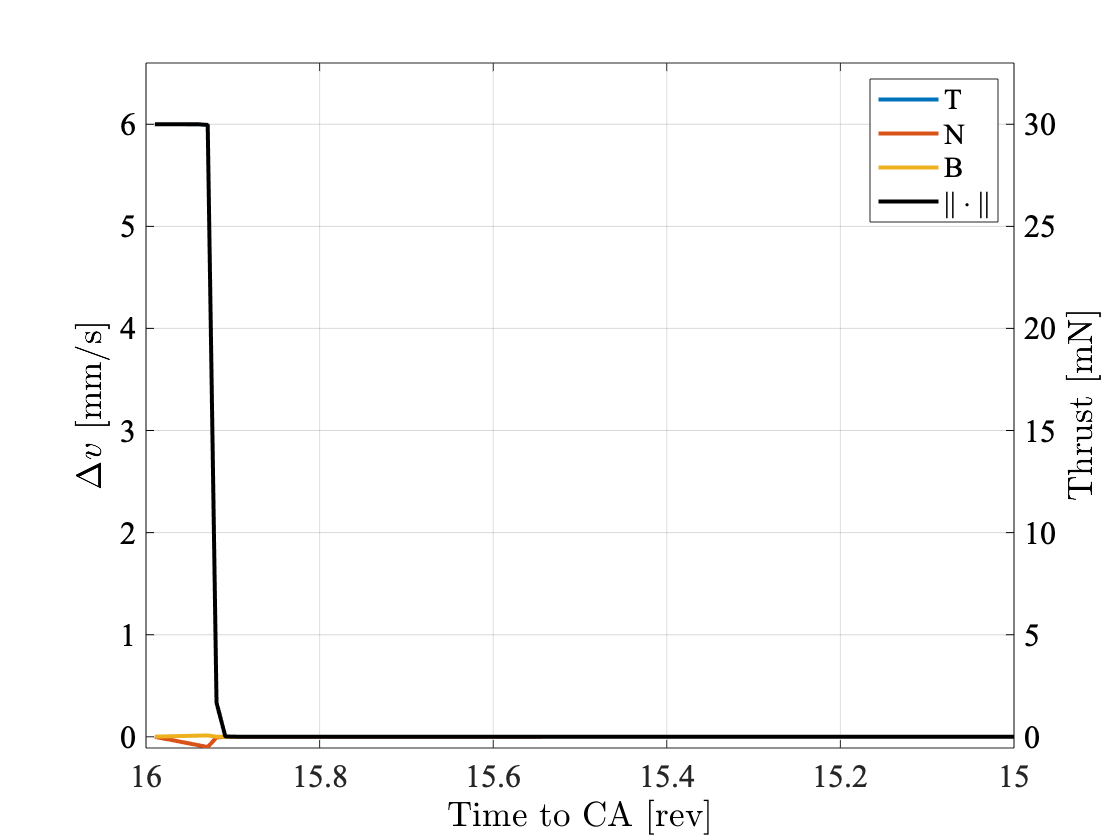}
\label{2}}\hfill
\subfigure[Solution \#3, local minimum $\Delta v = 45.2$ mm/s ]{
\includegraphics[width=0.47\textwidth]{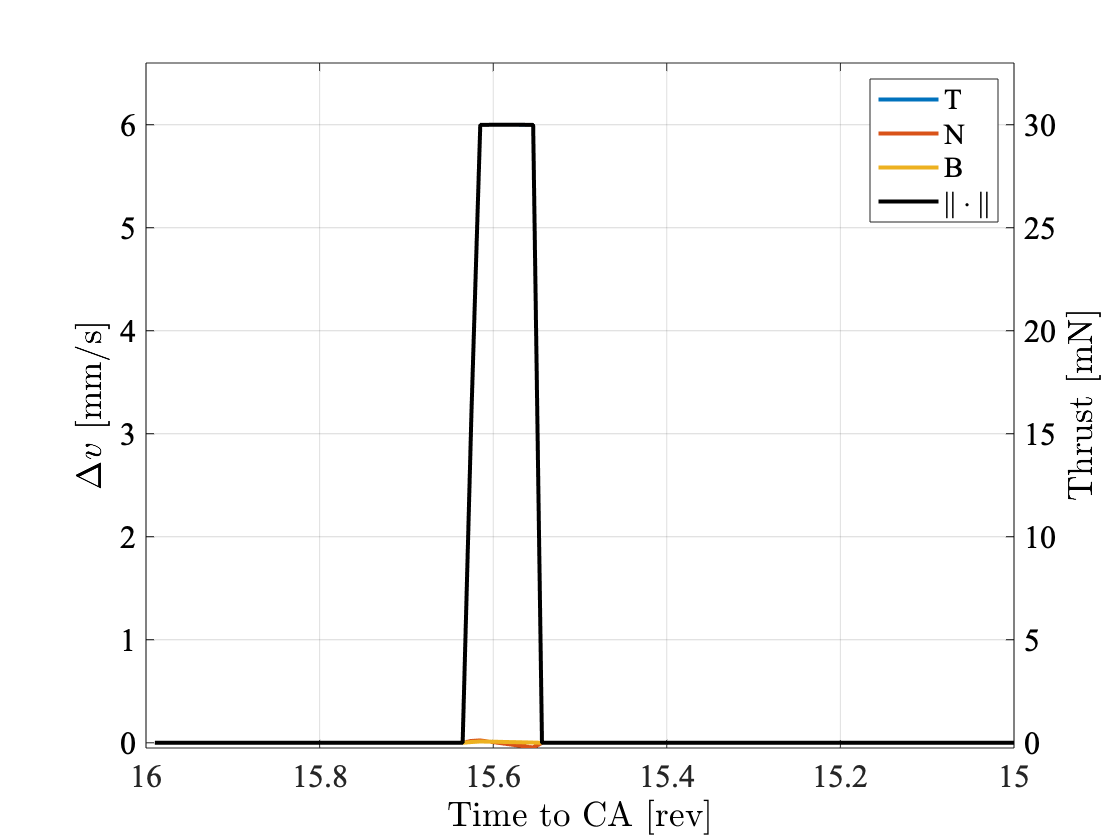}
\label{3}}\hfill
\subfigure[Solution \#4, local minimum $\Delta v = 74.9$ mm/s ]{
\includegraphics[width=0.47\textwidth]{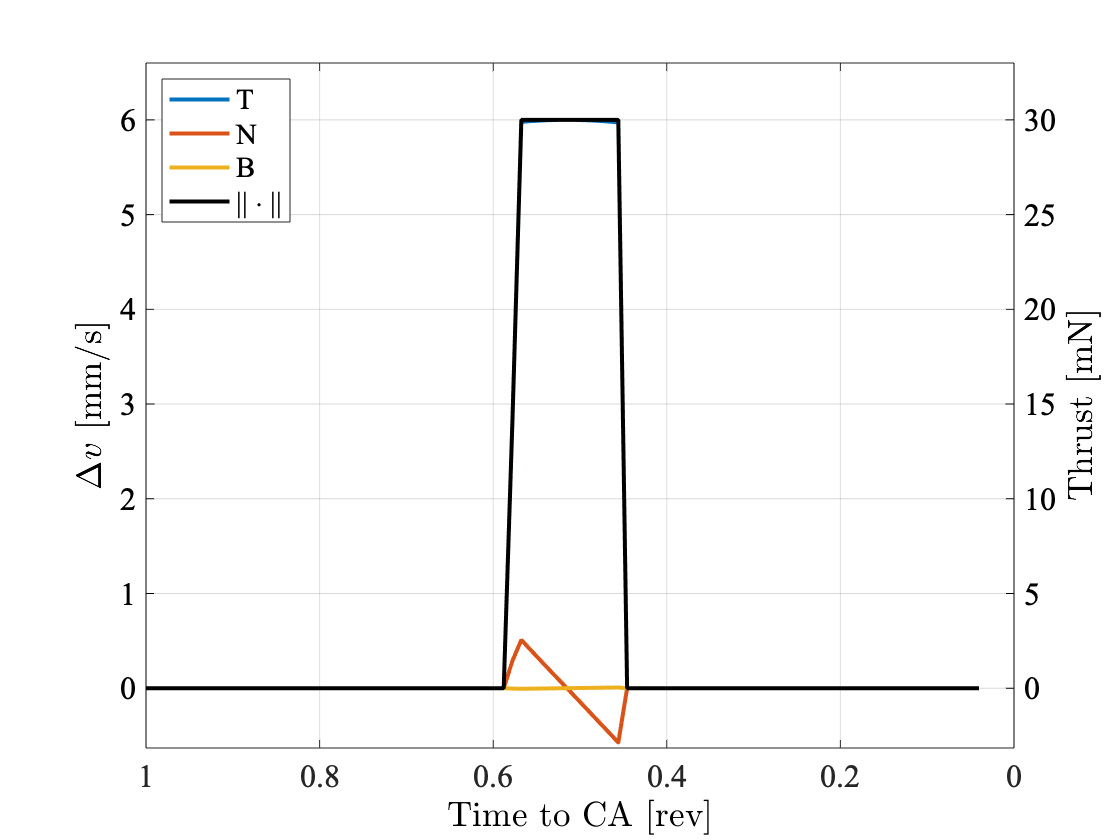}
\label{4}}\hfill
\subfigure[Solution \#5, local minimum $\Delta v = 82.9$ mm/s ]{
\includegraphics[width=0.47\textwidth]{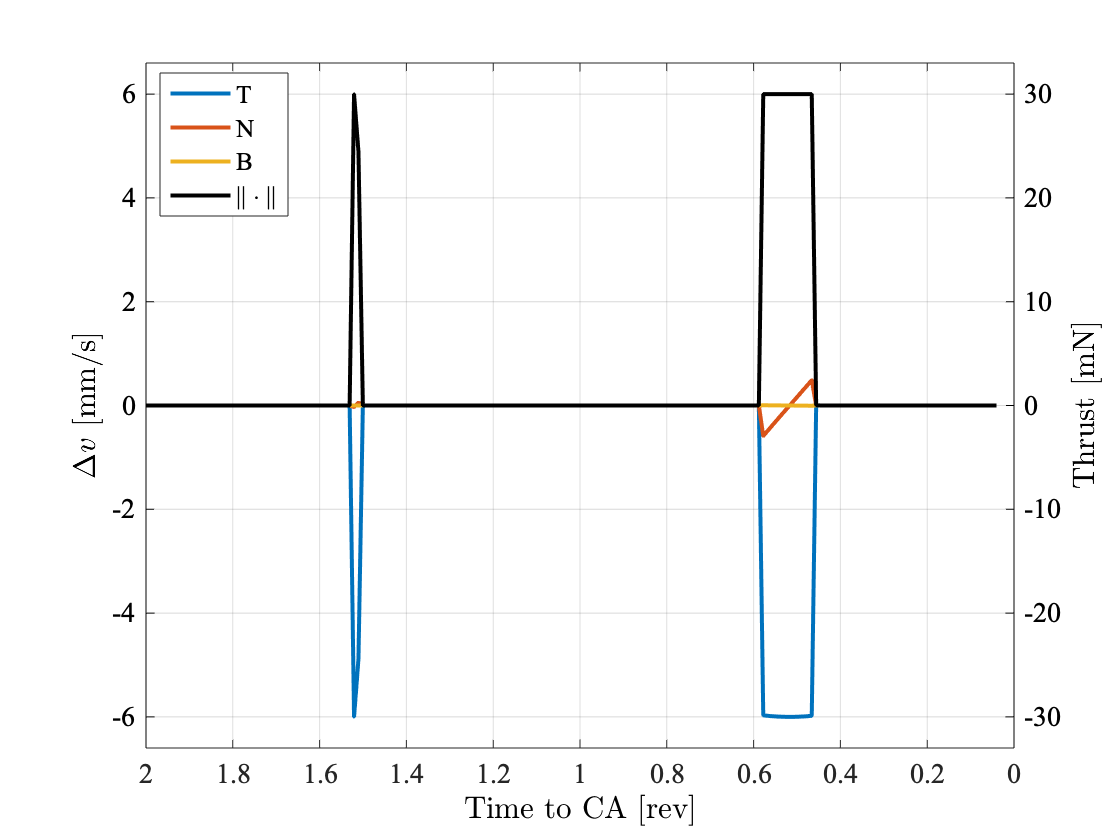}
\label{5}}\hfill
\caption{Example of multiple local minima for long alert times and extended maneuvering windows.}
\label{effectOftimeBis}
\end{figure}

For long alert times, the objective function structure can deviate significantly from the one in Fig. \ref{EllipseDv} due to multiple local minima. An example with five local minima is provided in Fig. \ref{ellipseDvAnomalous} for test case \#10 in the dataset, an alert time of 16 orbits, and up to 1,575 impulses (covering the 16 revolutions entirely). Figures \ref{1}--\ref{5} show that: for the first two minima the maneuver is applied at the earliest time (with thrust mainly either in the tangential or minus tangential direction); the third minimum corresponds to a maneuver 15.5 revolutions before the conjunction; and the spacecraft is maneuvered close to the conjunction at the two local minima with highest $\Delta v$. As a result, when long alert times and extended maneuvering windows are allowed, multiple initial guesses (more than two) are needed to guarantee convergence to the global optimum. The reason for the appearance of multiple local minima and the determination of the guess to compute the global one require further investigation.   

\begin{figure}
\subfigure[Objective function profile on the boundaries of the avoidance region]{
\includegraphics[width=0.47\textwidth]{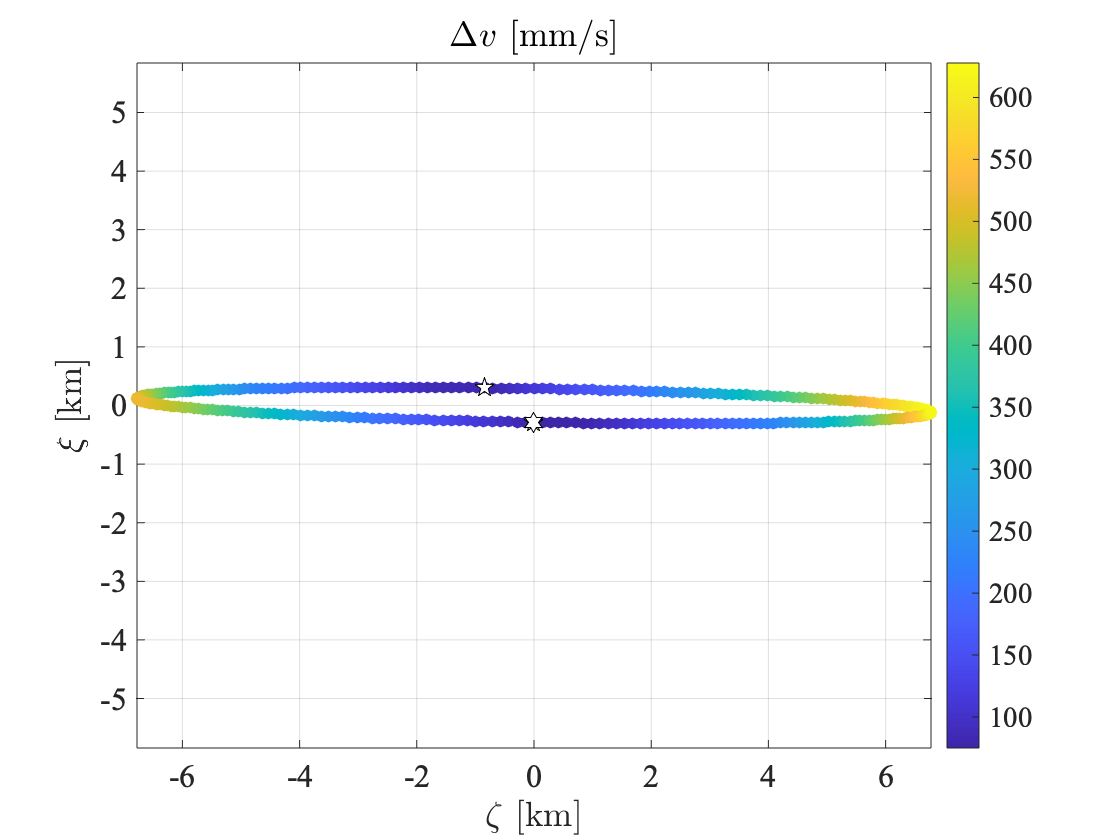}
\label{anomalEllipse}}\hfill
\subfigure[Impulses profile]{
\includegraphics[width=0.47\textwidth]{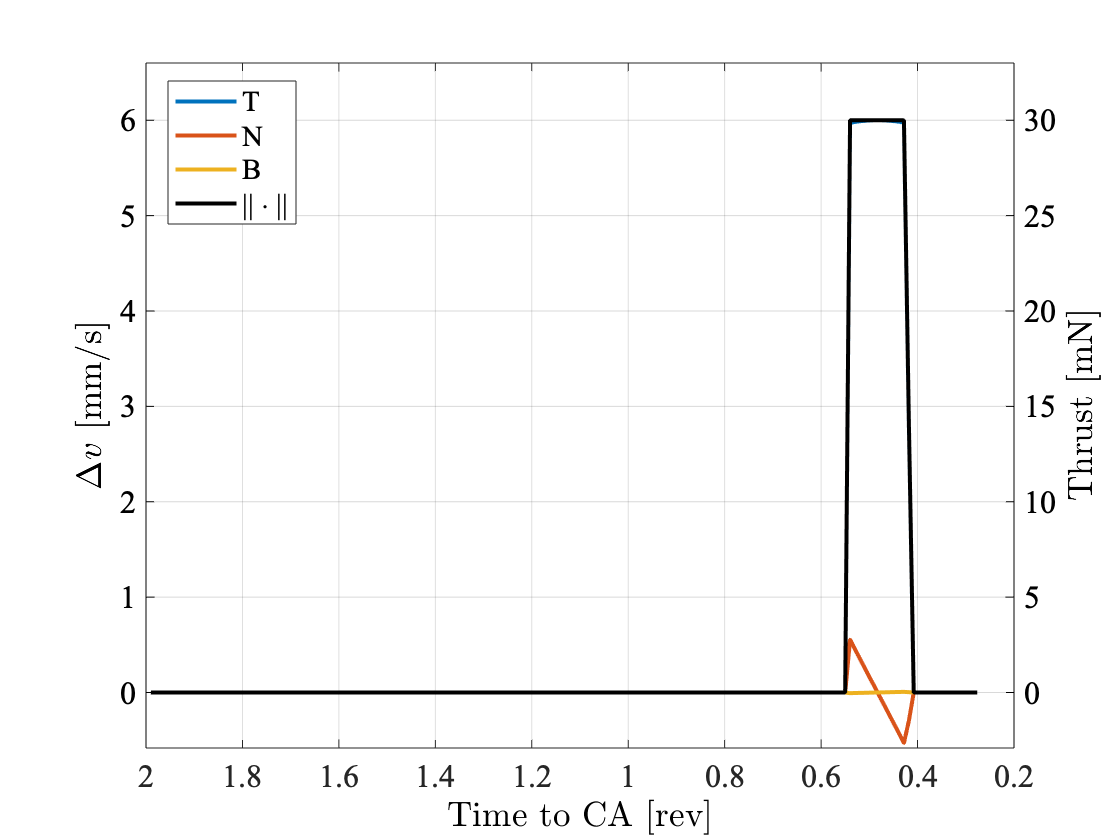}
\label{anomalControl}}\hfill
\caption{Objective function and global optimum solution for case \#10 with short alert time.}
\label{anomalCase}
\end{figure}

If the alert time is reduced to two revolutions and 170 impulses, only two minima remain, as shown in Fig. \ref{anomalEllipse}. The global optimum, indicated with the hexagram, coincides with solution \#4 from Fig. \ref{4}, and is reported for clarity in Fig. \ref{anomalControl}. Thus, depending on the conjunction configuration, there are cases in which the global optimum is obtained by thrusting towards the end of the control window.   

\subsection{Maximum Impulse Magnitude}\label{sec:maxT}
 \begin{figure}[h!]
\subfigure[$\Delta \bar{v} = 200$ mm/s]{
\includegraphics[width=0.47\textwidth]{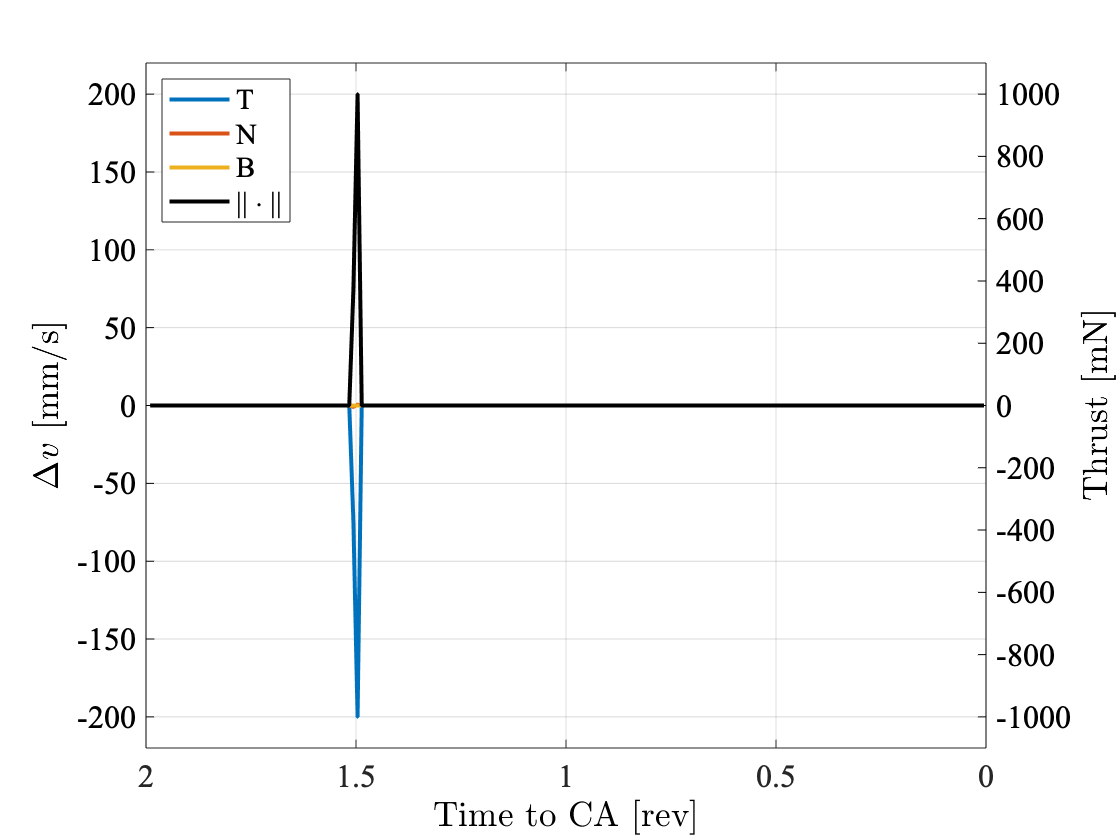}
\label{effectOfThrust10}}\hfill
\subfigure[$\Delta \bar{v} = 10$ mm/s]{
\includegraphics[width=0.47\textwidth]{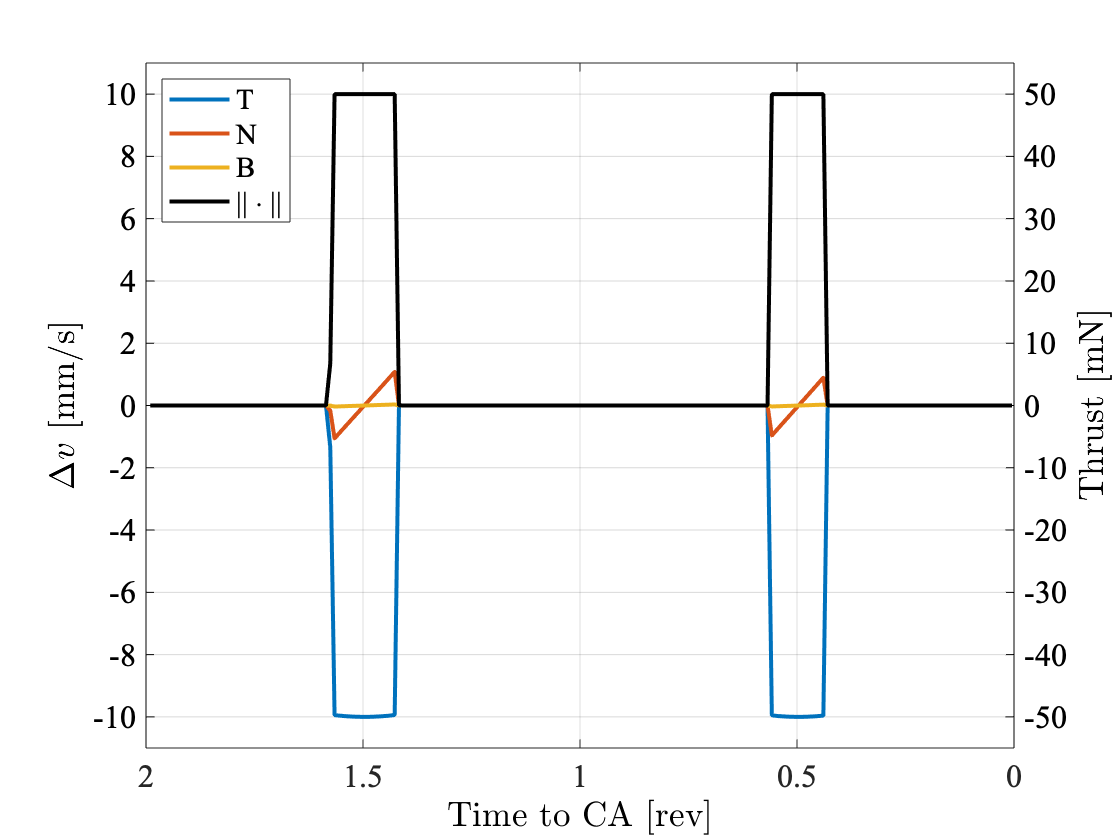}
\label{effectOfThrust1}}\hfill
\subfigure[$\Delta \bar{v} = 3$ mm/s]{
\includegraphics[width=0.47\textwidth]{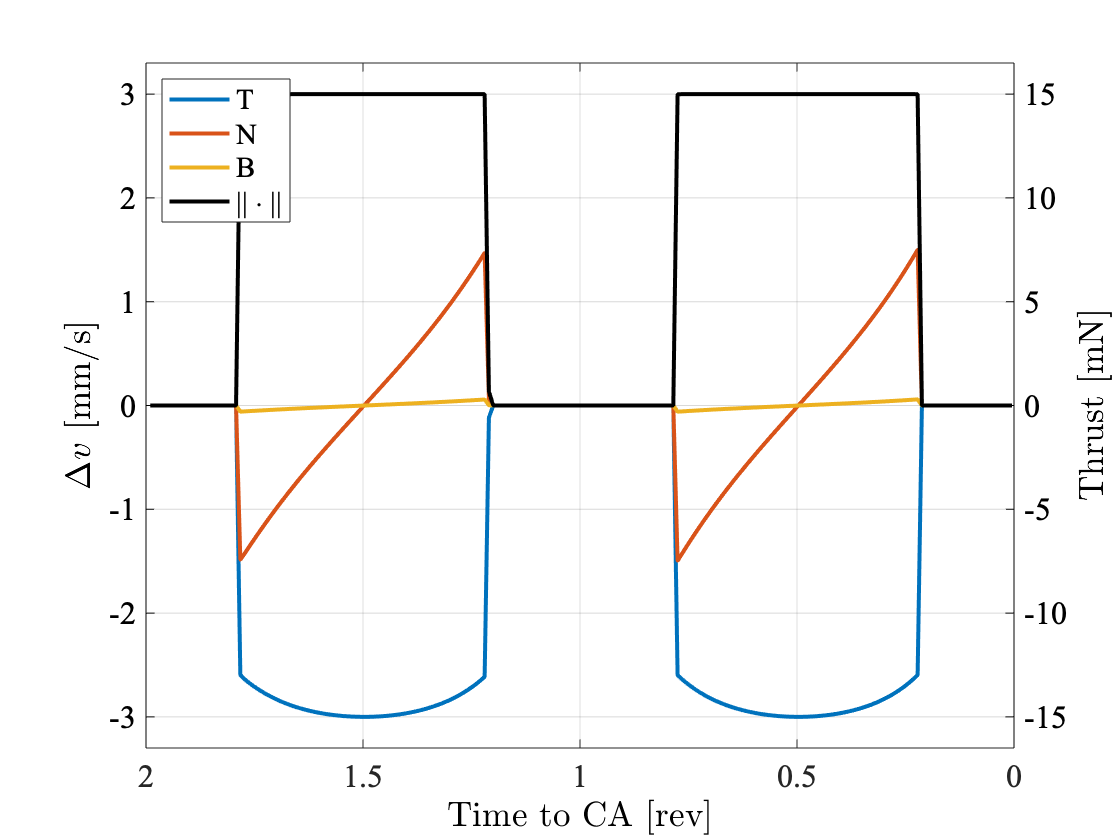}
\label{effectOfThrust01}}\hfill
\subfigure[$\Delta \bar{v} = 2.5$ mm/s]{
\includegraphics[width=0.47\textwidth]{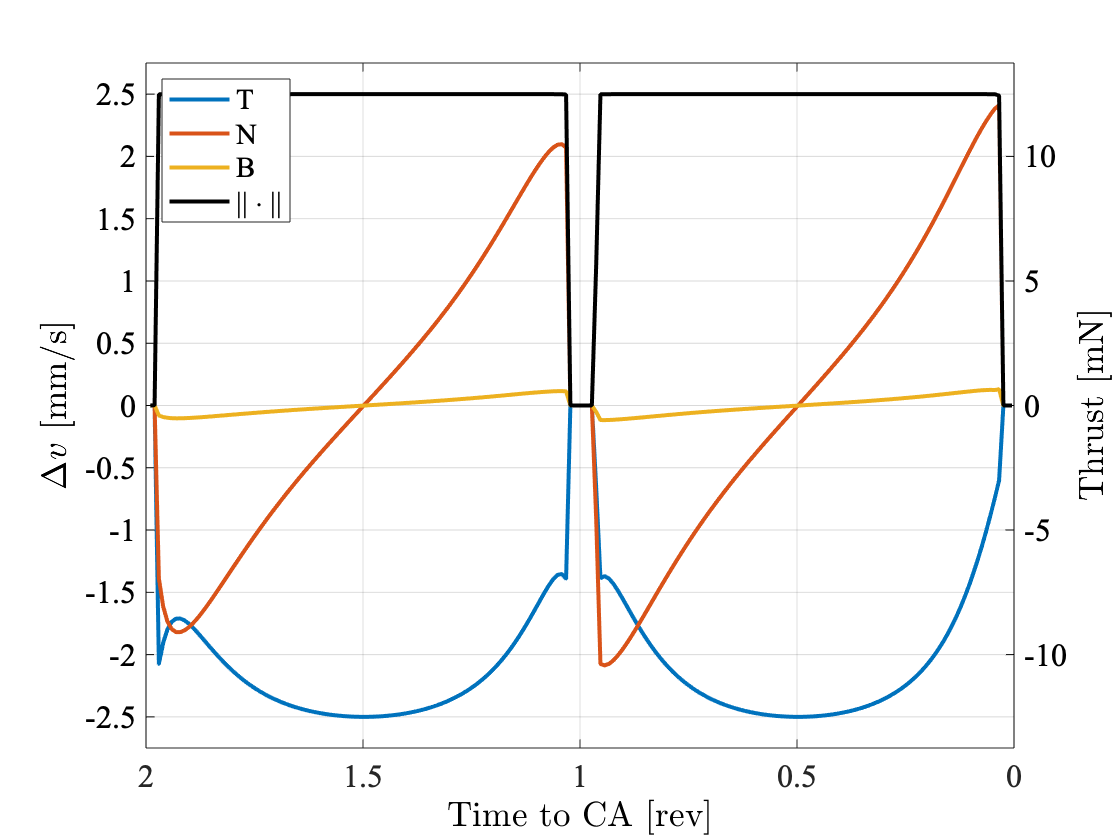}
\label{effectOfThrust005}}\hfill
\caption{Impulses profile as a function of maximum impulse magnitude}
\label{effectOfThrust}
\end{figure}

The effect of varying the impulse magnitude is studied. In Fig.\ref{effectOfThrust} the maximum impulse available is reduced from 200 mm/s to 2.5 mm/s, corresponding to a constant thrust from 1 N to 12.5 mN for a reference 300 kg spacecraft (i.e., going from a chemical to an electric propulsion system). This change has a significant impact on the number of impulses and the efficiency of the maneuver. In the first case, three impulses are sufficient (only one at the maximum thrust value), resulting in a $\Delta v = 275.0$ mm/s. In the latter case, the maneuver requires 191 impulses and a total $\Delta v = 476.1$ mm/s. Closer to the conjunction, the radial component of thrust becomes more relevant, at points even greater than the tangential one, as reported in \cite{Reiter2018}. 

\subsection{Effect of the Eccentricity of the Primary}\label{sec:ecc}

\begin{figure}[h!]
\subfigure[Thrust arcs on primary trajectory, $e = 0.05$ ]{
\includegraphics[width=0.47\textwidth]{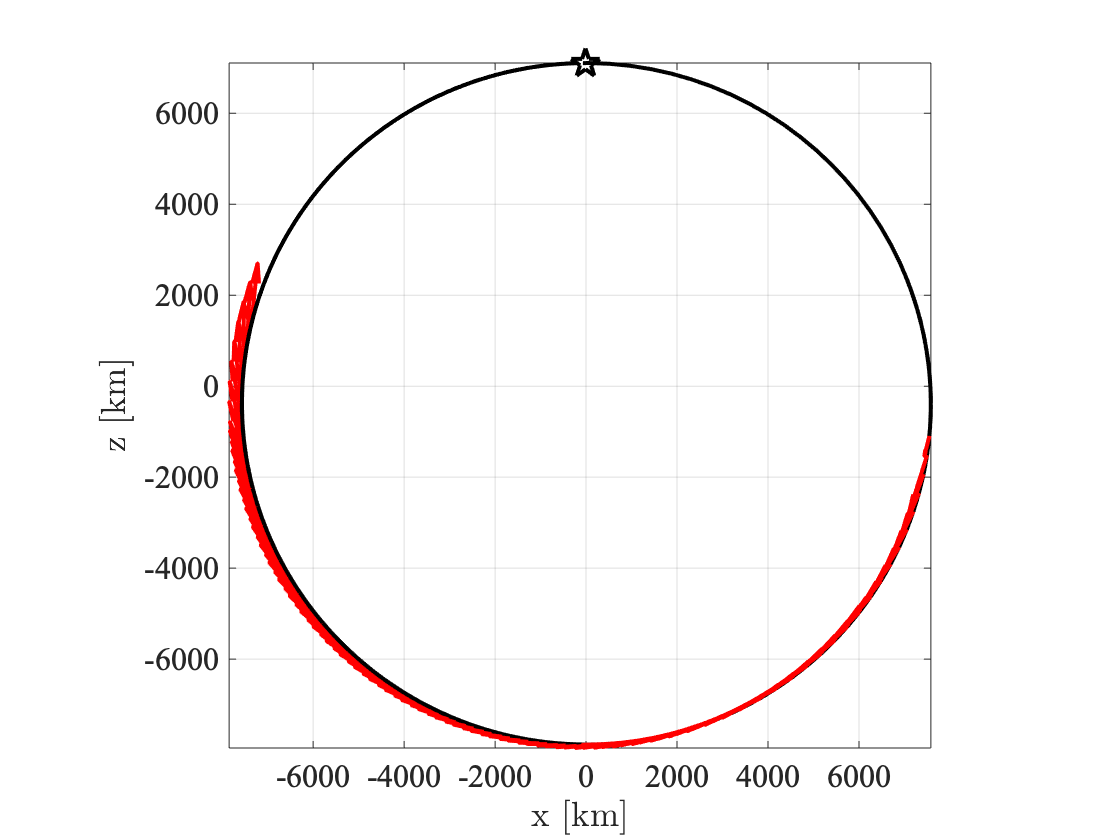}
\label{ecc005Trj}}\hfill
\subfigure[Impulses timing and components, $e = 0.05$]{
\includegraphics[width=0.47\textwidth]{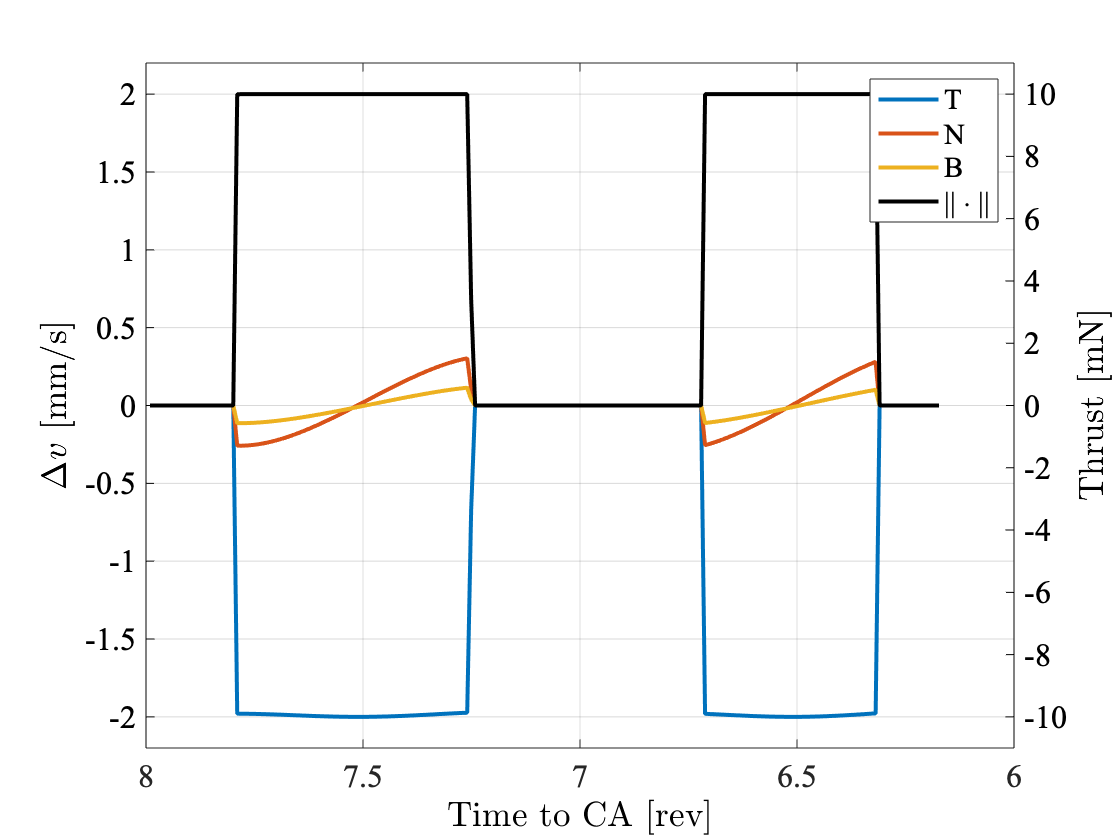}
\label{ecc005Contr}}\hfill
\subfigure[Thrust arcs on primary trajectory, $e = 0.5$]{
\includegraphics[width=0.47\textwidth]{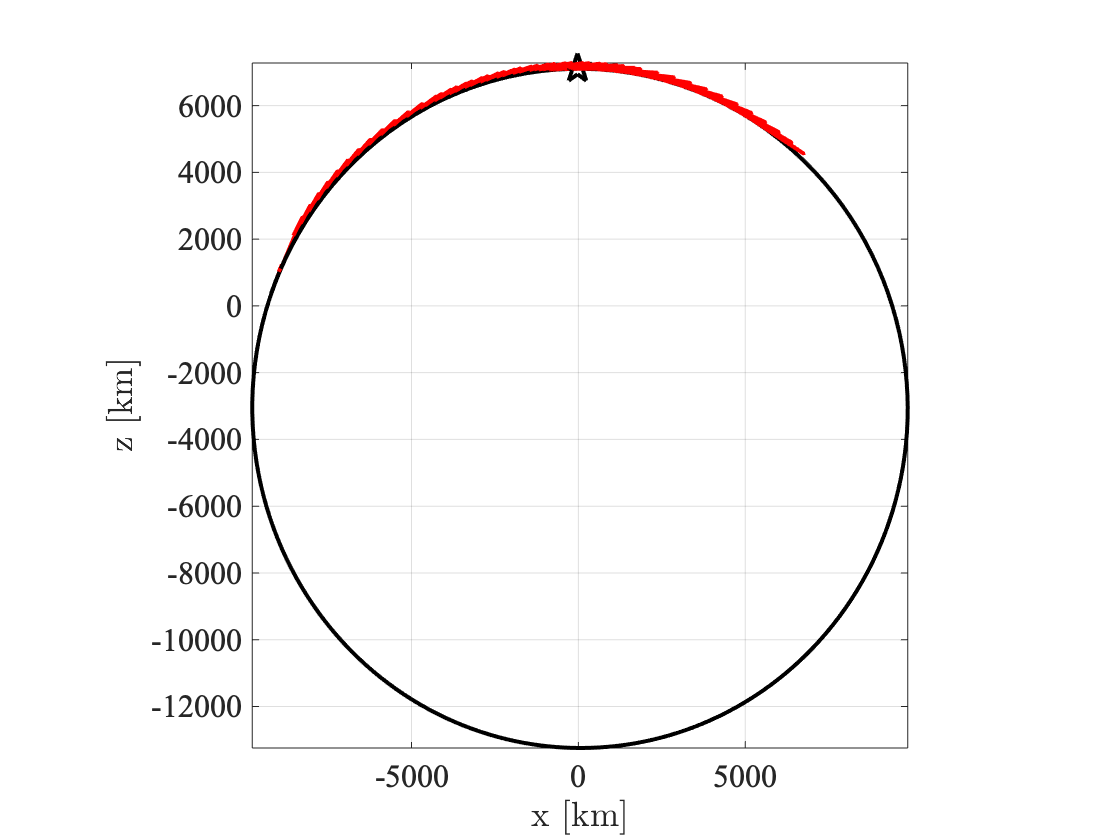}
\label{ecc05Trj}}\hfill
\subfigure[Impulses timing and components, $e = 0.5$]{
\includegraphics[width=0.47\textwidth]{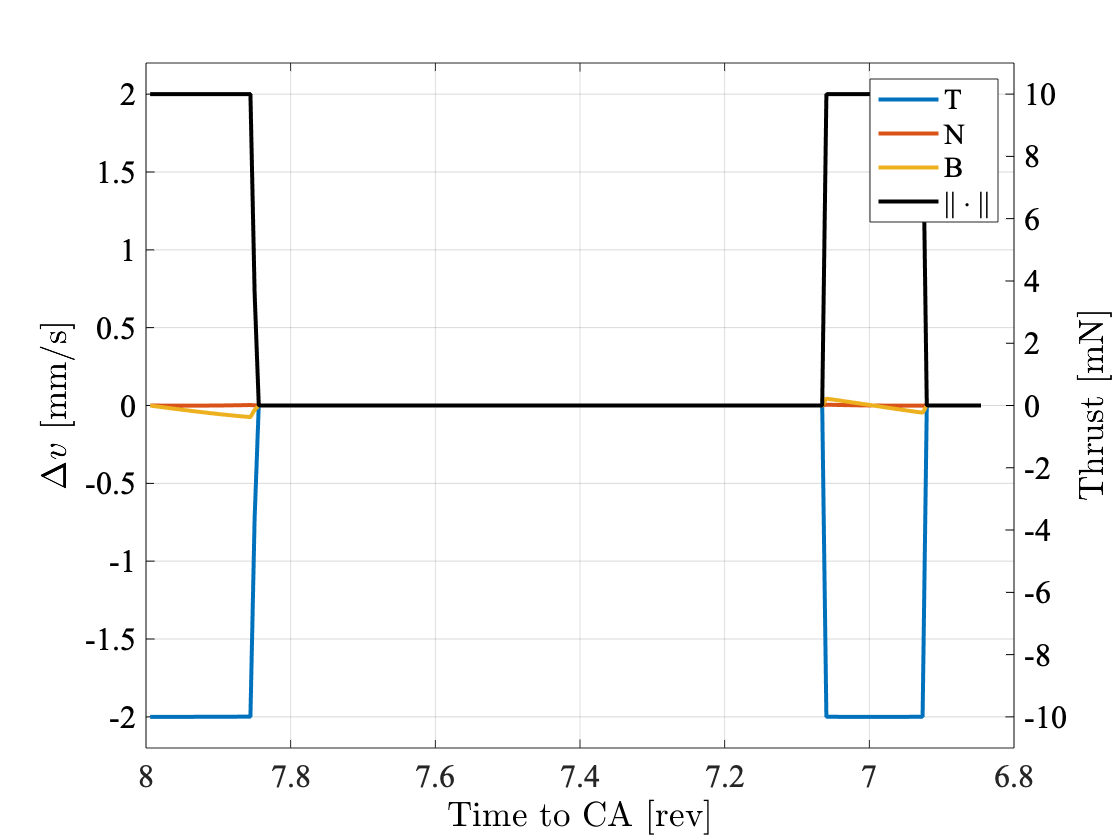}
\label{ecc05Contr}}\hfill
\caption{Effect of eccentricity on the maneuver}
\label{max_ecc}
\end{figure}

The dataset includes highly eccentric cases for the secondary object (with a maximum value of 0.5264), whereas the eccentricity of the primary is always relatively small, with a maximum of 0.0175. To study the suitability of the approach to deal with eccentric cases, we modify the primary's velocity at the conjunction epoch by adding a $\Delta v$ in the tangential direction. As the starting orbit is almost circular, the conjunction occurs close to the perigee of the primary. The maneuver can start 8 revolutions before the conjunction, exploiting up to to 200 impulses separated by one minute. The maximum impulse magnitude is constraned to 2 mm/s, and the constraint on maximum collision probability is applied. Figure \ref{max_ecc} shows the optimal maneuvers for eccentricity $e = 0.05$ and $e= 0.5$. The optimal thrust arc location moves from the opposition (as shown in Fig. \ref{ecc005Trj}) to across the conjunction (as shown in Fig. \ref{ecc05Trj}). This is probably due to an increased maneuver efficiency around the perigee. Correspondingly, the $\Delta v$ drops from 206.7 mm/s to 98.7 mm/s.

\subsection{Full Dataset Simulation}\label{sec:extensiveSim}
The proposed approach is run on 2,170 test cases, using either $P_{C,\max} \le 10^{-4}$, $P_C \le 10^{-6}$, or $\Delta r_{CA} \ge 2$ km as constraints. The maneuvers can be implemented starting from 2 revolutions before the encounter using up to 170 impulses. The maximum impulse magnitude is 6 mm/s, corresponding to a constant thrust of 30 mN in each segment. These simulations are used to assess the properties of the method for a large set of conjunctions, including some considerations on the accuracy of linearizations and the impact of perturbations. In the histograms presented in this section, data in the range $5-95$th percentile are plotted for the sake of readability.  

\begin{figure}[h!]
\subfigure[Distribution of $P_C$]{
\includegraphics[width=0.47\textwidth]{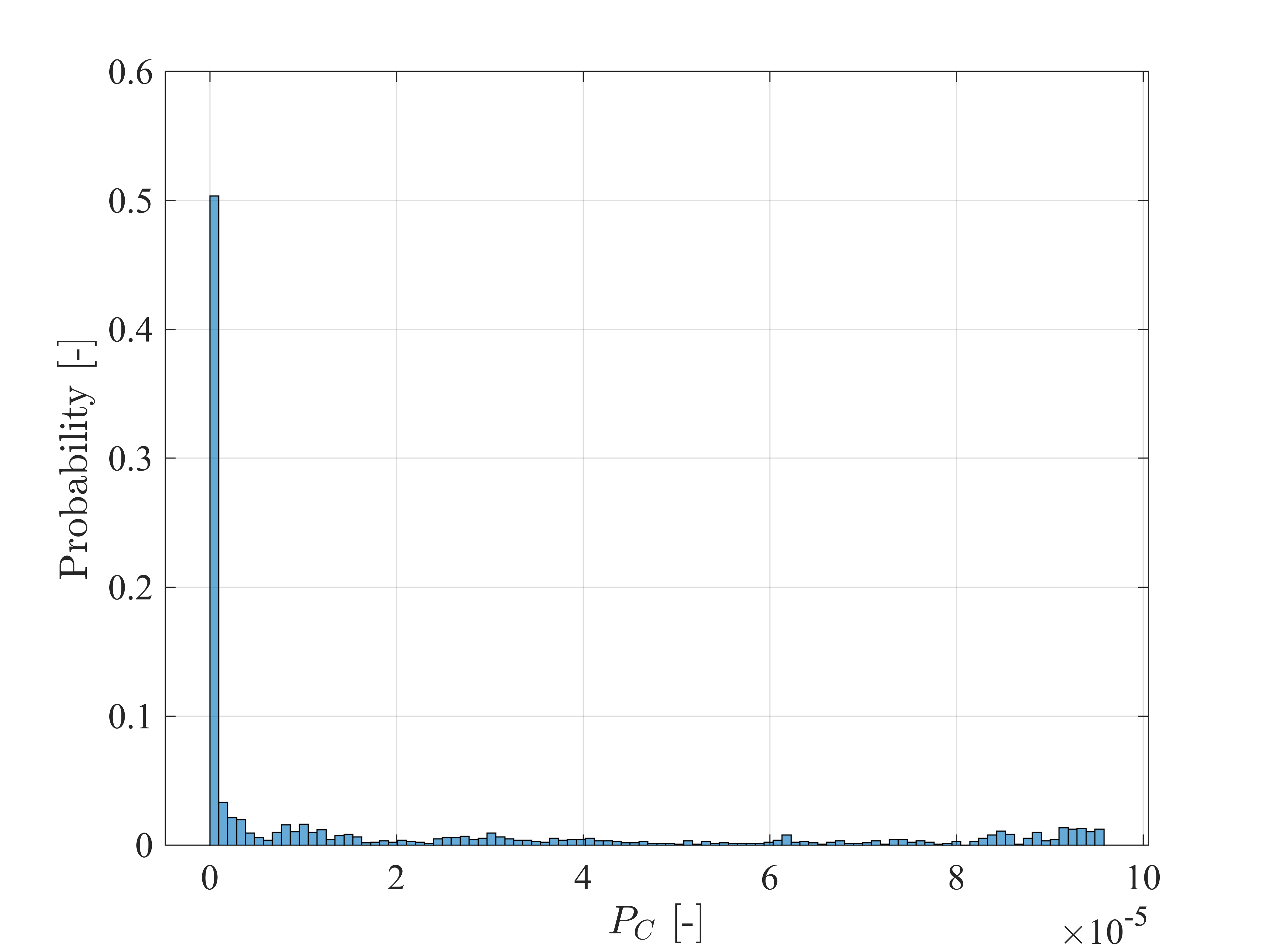}
\label{max_riskPc}}\hfill
\subfigure[Distribution of $\Delta r_{CA}$]{
\includegraphics[width=0.47\textwidth]{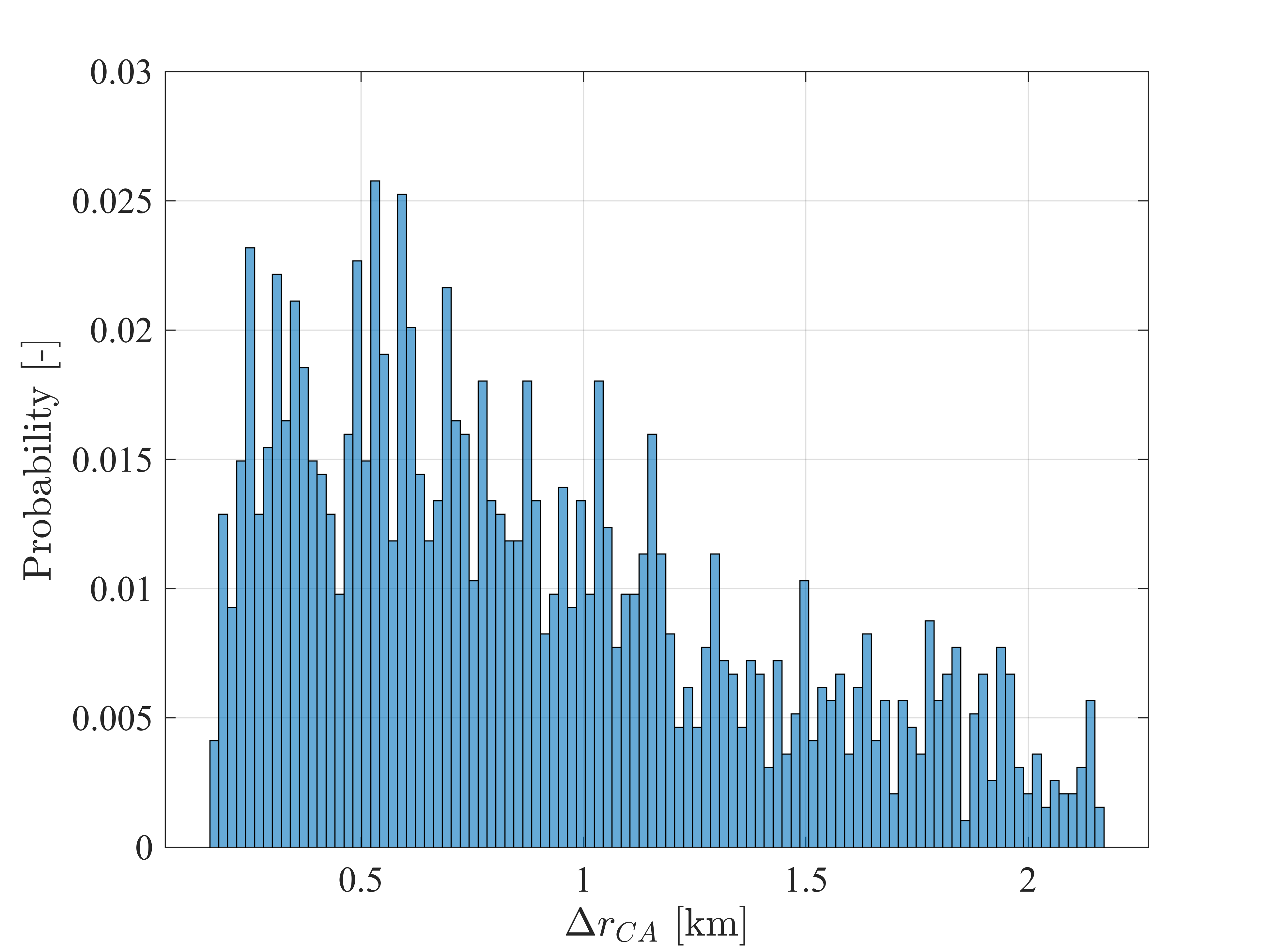}
\label{max_riskminRelDist}}\hfill
\subfigure[Distribution of $\Delta v$]{
\includegraphics[width=0.47\textwidth]{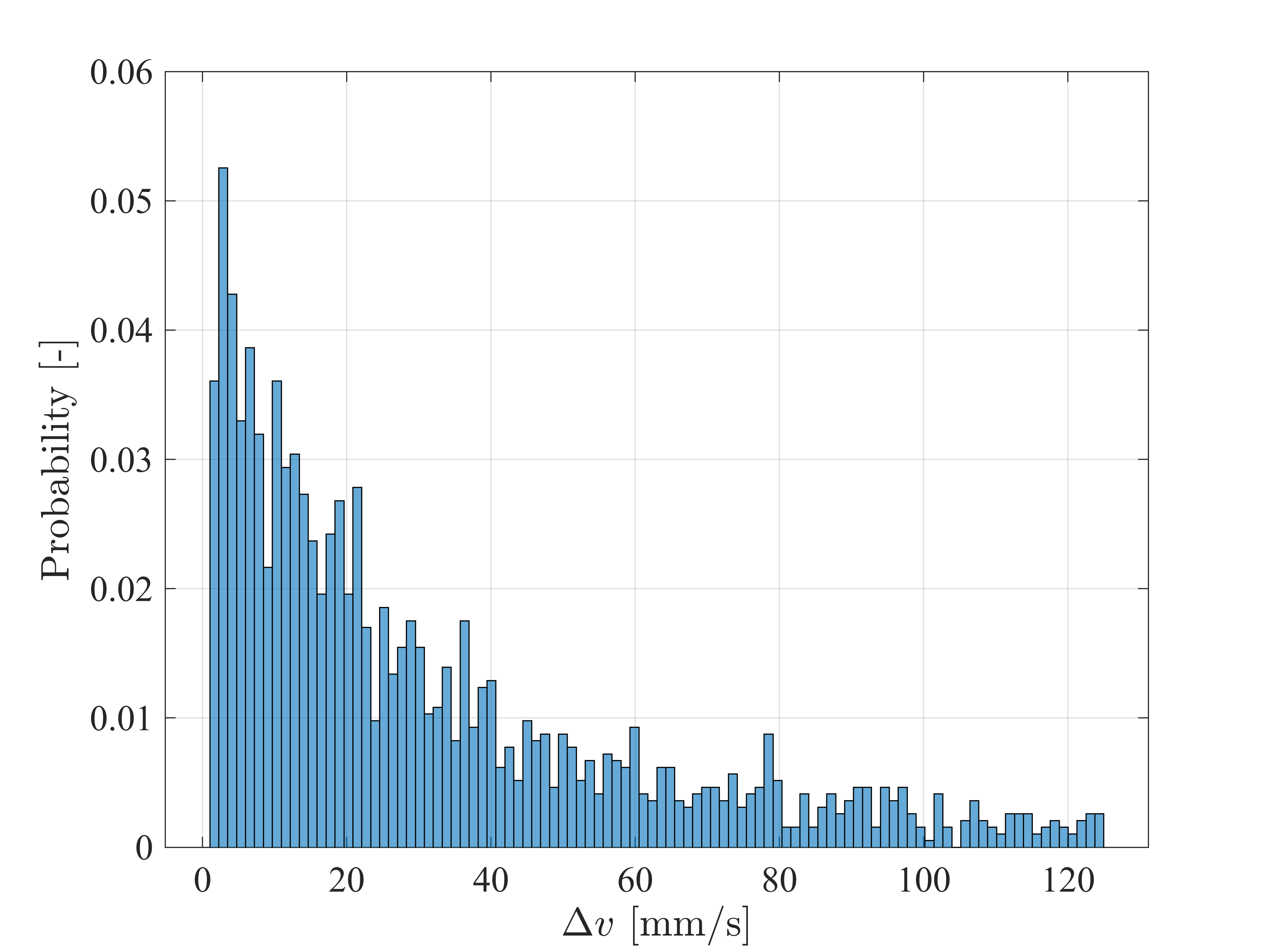}
\label{max_riskdv}}\hfill
\subfigure[Distribution of impulses]{
\includegraphics[width=0.47\textwidth]{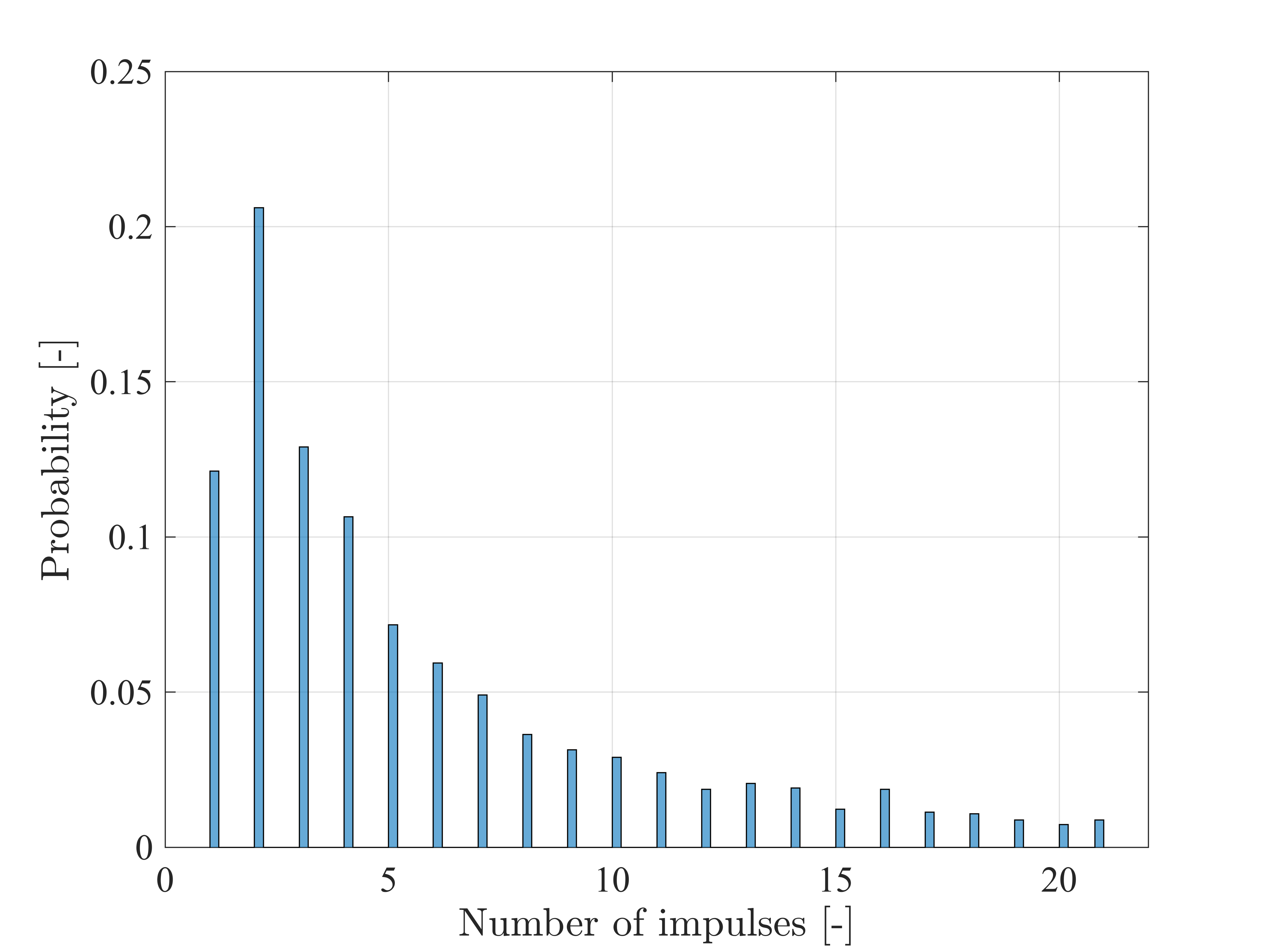}
\label{max_risknImp}}\hfill
\caption{Maximum risk constraint}
\label{max_risk}
\end{figure}

\begin{figure}[h!]
\subfigure[Distribution of $P_{C,\max}$]{
\includegraphics[width=0.47\textwidth]{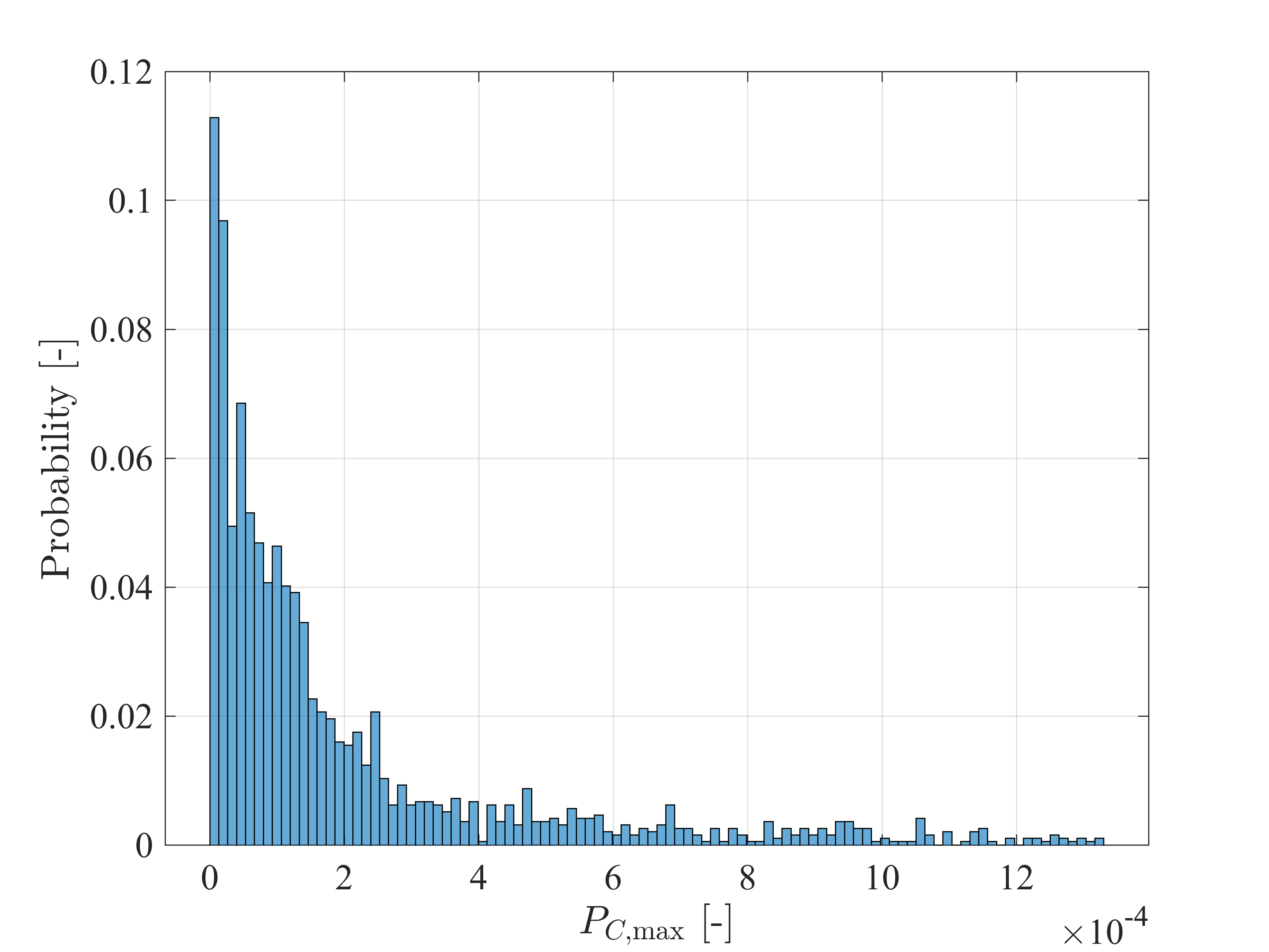}
\label{riskmaxPc}}\hfill
\subfigure[Distribution of $\Delta r_{CA}$]{
\includegraphics[width=0.47\textwidth]{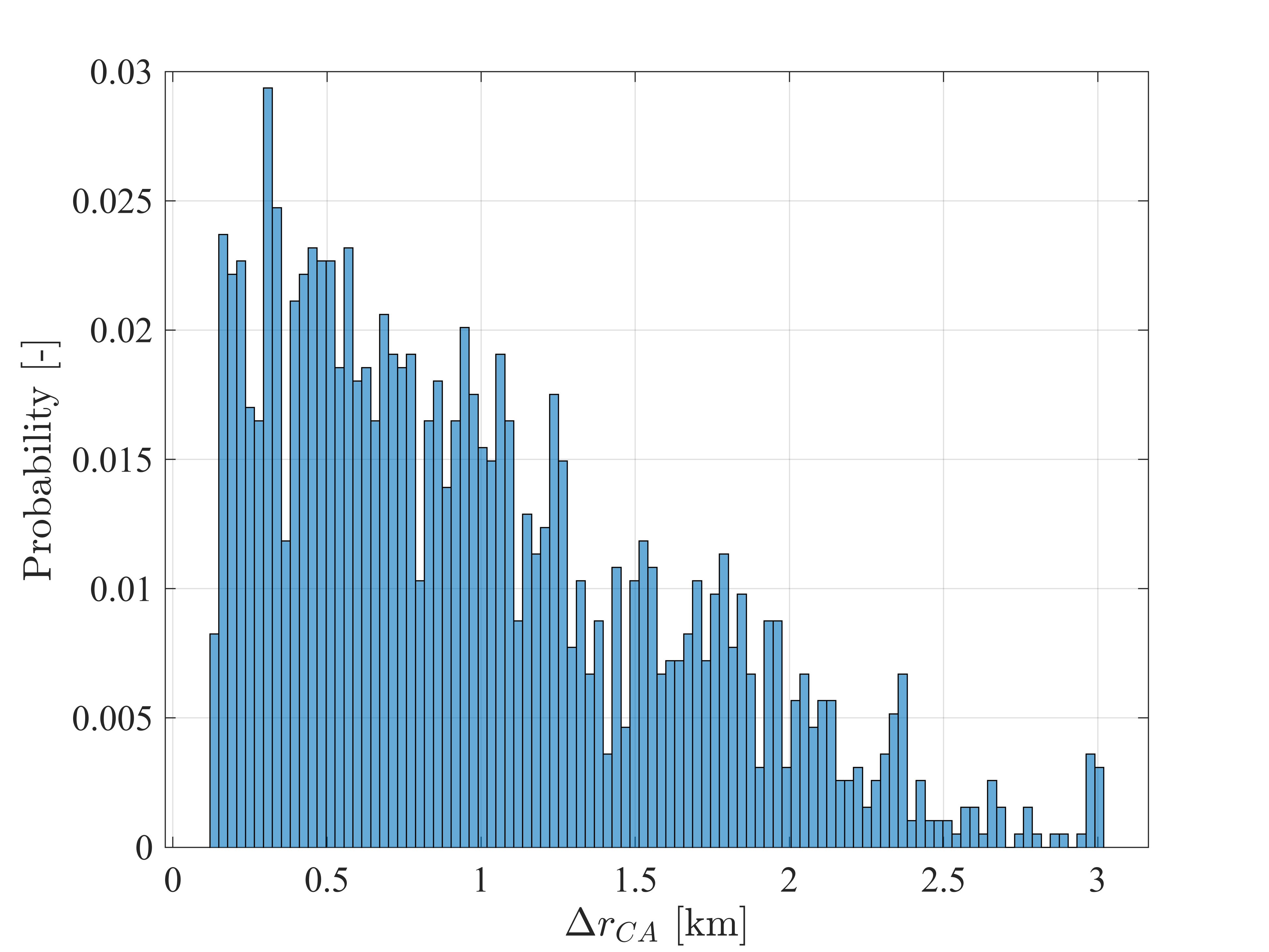}
\label{riskminRelDist}}\hfill
\subfigure[Distribution of $\Delta v$]{
\includegraphics[width=0.47\textwidth]{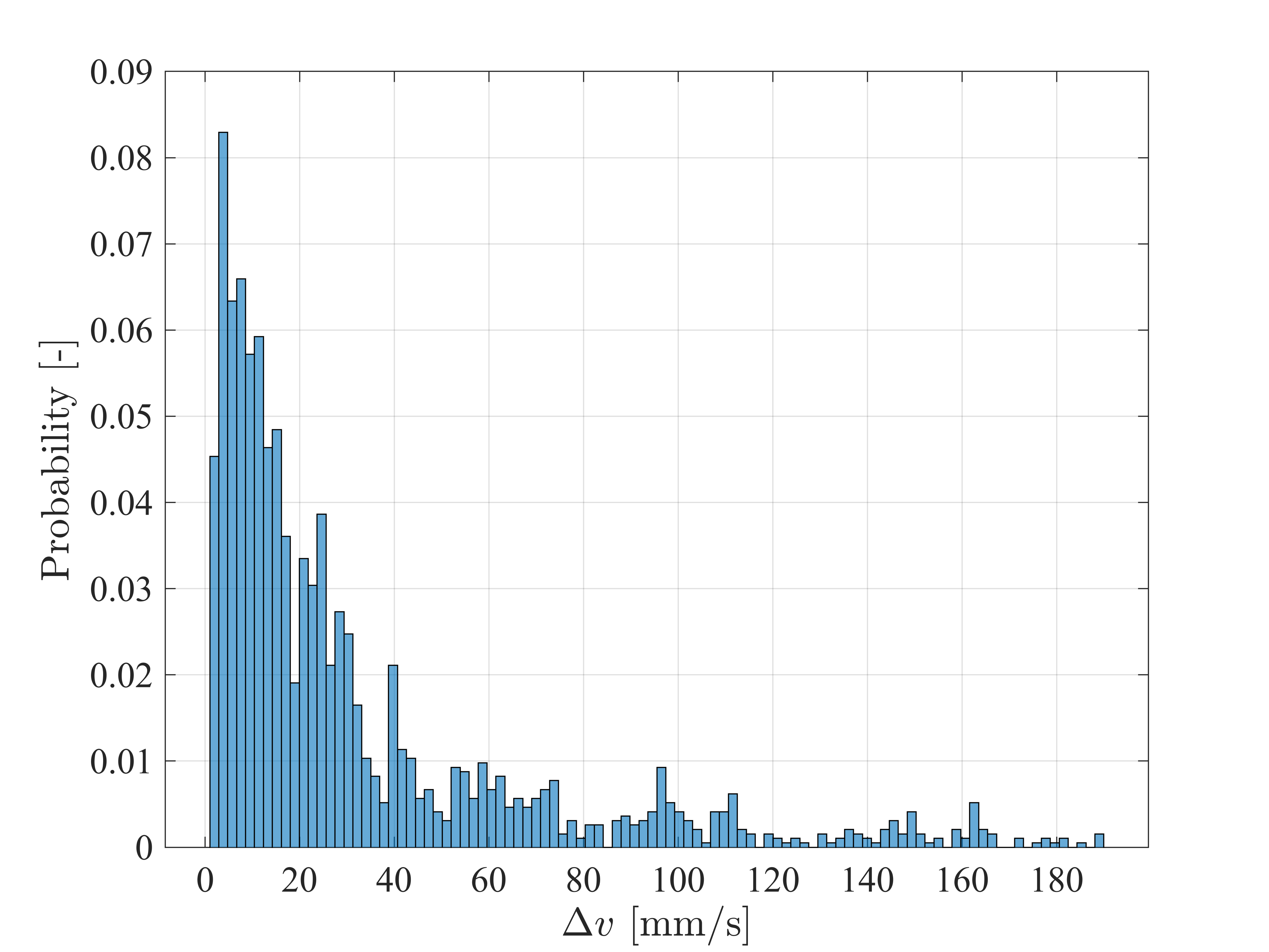}
\label{riskdv}}\hfill
\subfigure[Distribution of impulses]{
\includegraphics[width=0.47\textwidth]{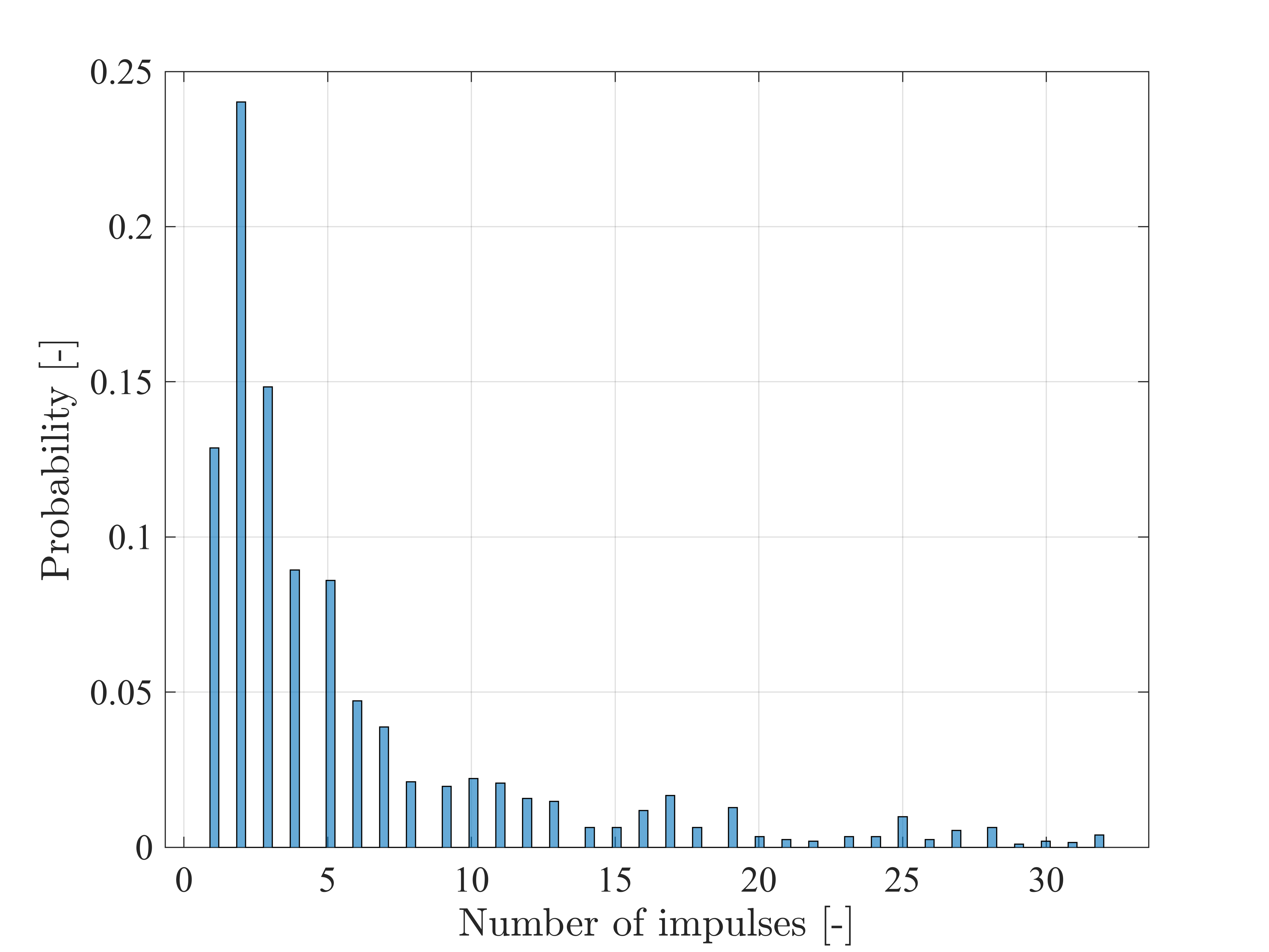}
\label{risknImp}}\hfill
\caption{Risk constraint}
\label{risk}
\end{figure}

\begin{figure}[h!]
\subfigure[Distribution of $P_{C}$]{
\includegraphics[width=0.47\textwidth]{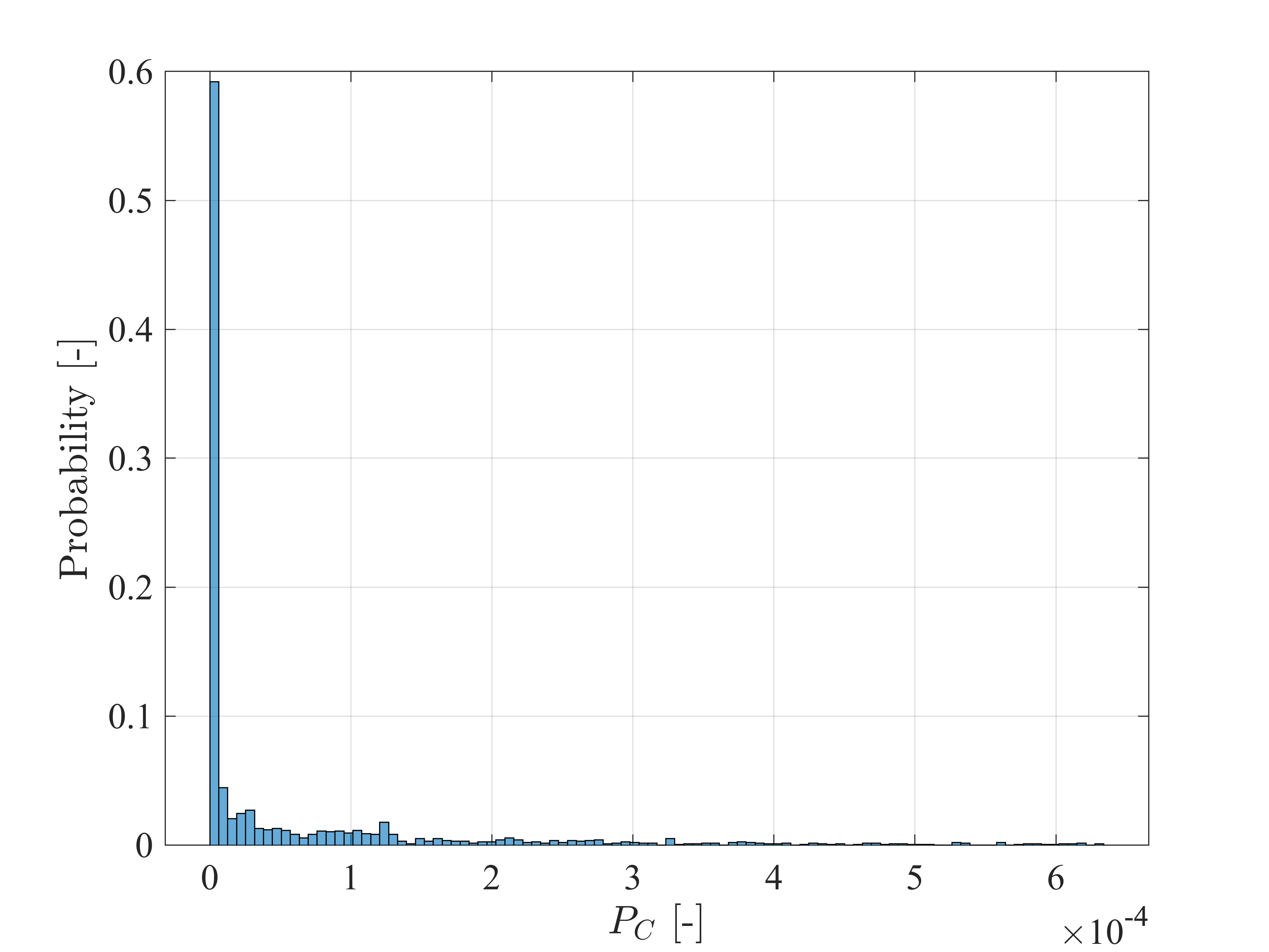}
\label{miss_distancePc}}\hfill
\subfigure[Distribution of $P_{C,\max}$]{
\includegraphics[width=0.47\textwidth]{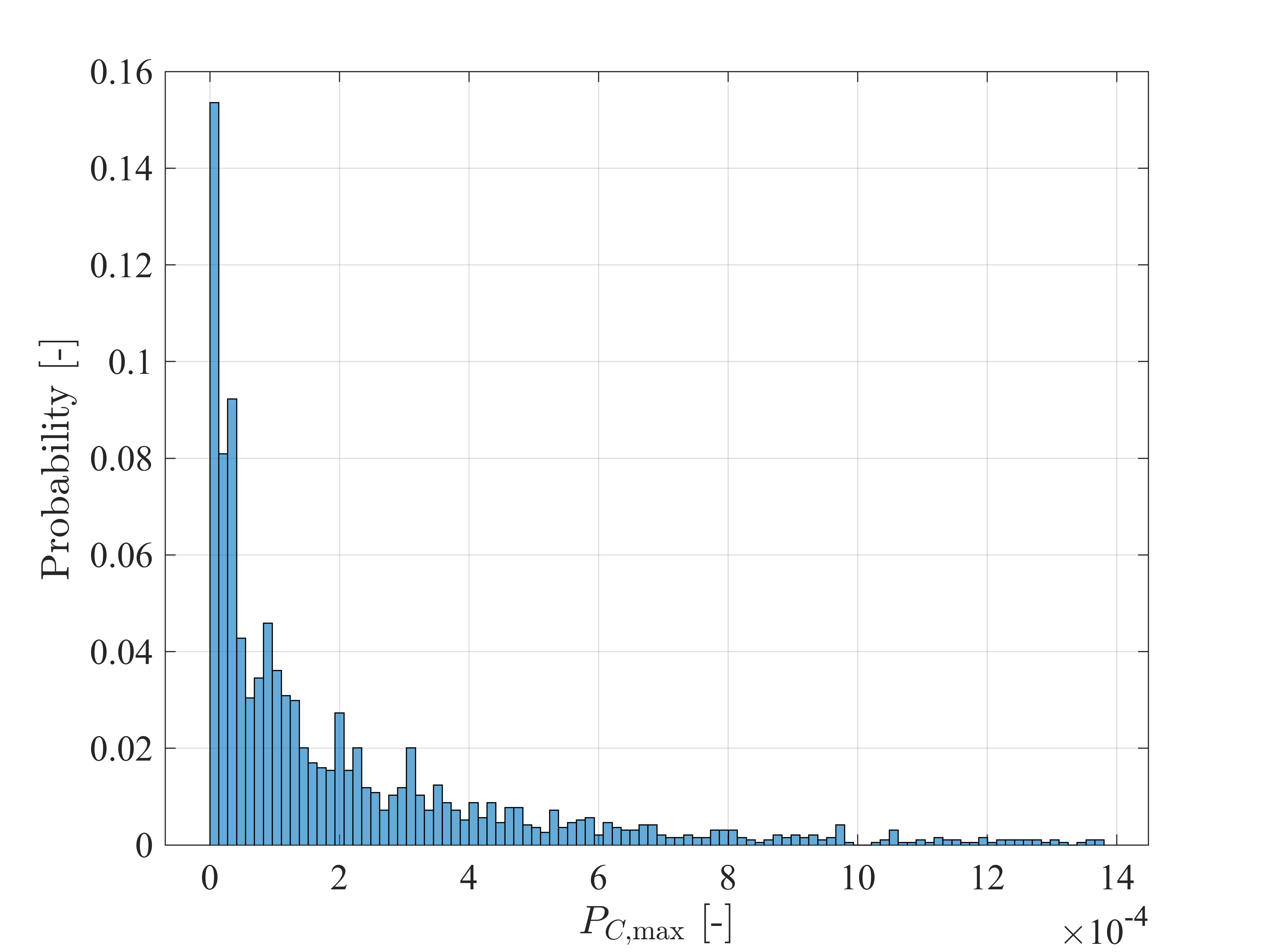}
\label{miss_distancemaxPc}}\hfill
\subfigure[Distribution of $\Delta v$]{
\includegraphics[width=0.47\textwidth]{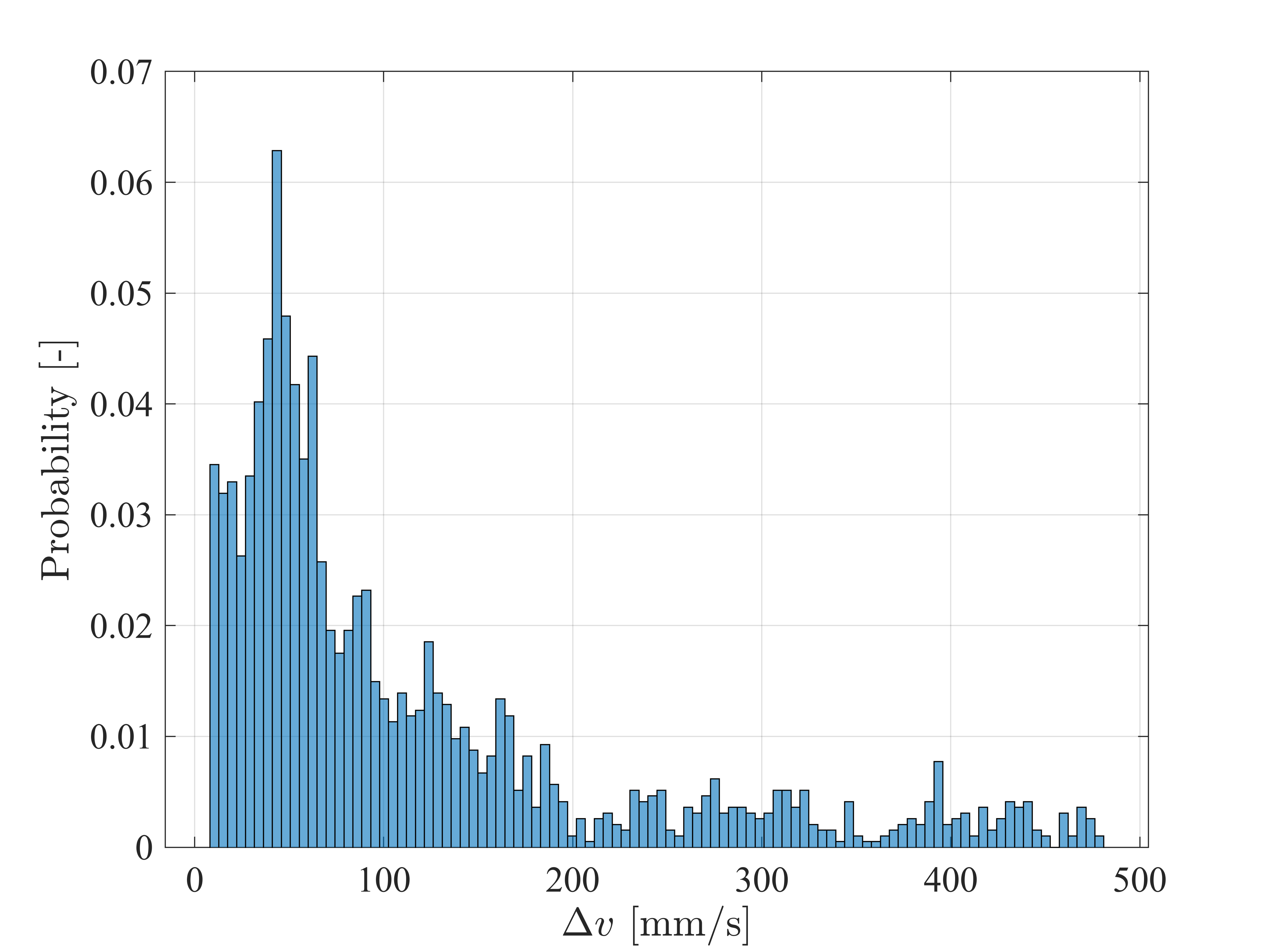}
\label{miss_distancedv}}\hfill
\subfigure[Distribution of impulses]{
\includegraphics[width=0.47\textwidth]{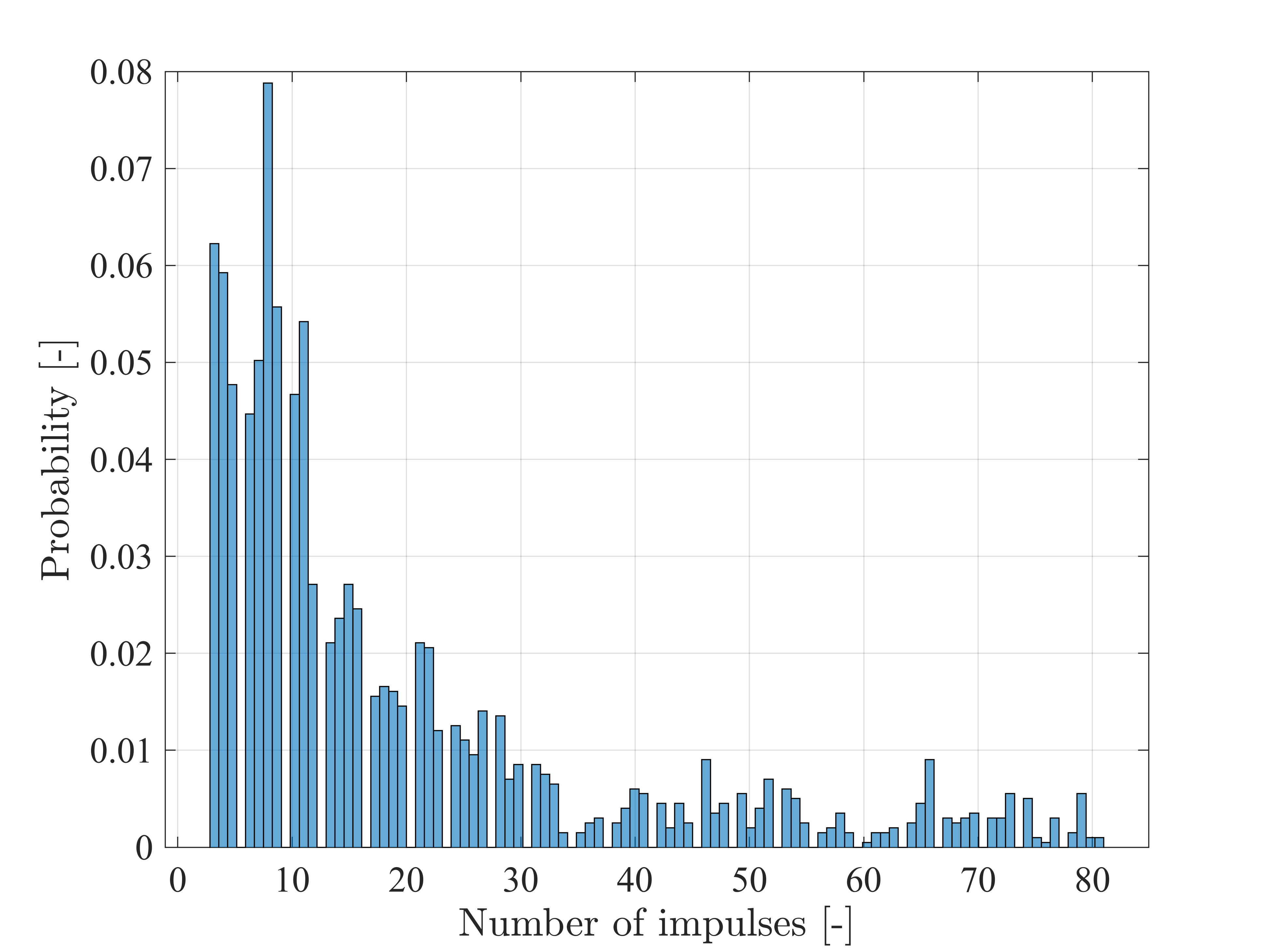}
\label{miss_distancenImp}}\hfill
\caption{Miss distance constraint}
\label{miss_distance}
\end{figure}

Figures \ref{max_risk}--\ref{miss_distance} show the distribution of relevant quantities when the three different constraints are applied. The $P_{C,\max}$ and $P_C$ cases produce similar results in term of median $\Delta v$ (21.2 mm/s versus 17.8 mm/s), relative closest encounter distance (0.7818 km versus 0.8627 km), and number of impulses (4 in both cases). The miss distance constraint is more demanding in terms of median $\Delta v$ (68.9 mm/s) and number of impulses (12), without providing on average either a lower collision risk or a lower maximum collision risk (median values of $1.0202 \times 10^{-4}$ and $1.8363 \times 10^{-6}$, respectively). This indicates that when orbital knowledge statistics are reliable, collision probability constraints should be prioritize over the miss distance one. 

\begin{figure}[h!]
\subfigure[First major iteration]{
\includegraphics[width=0.47\textwidth]{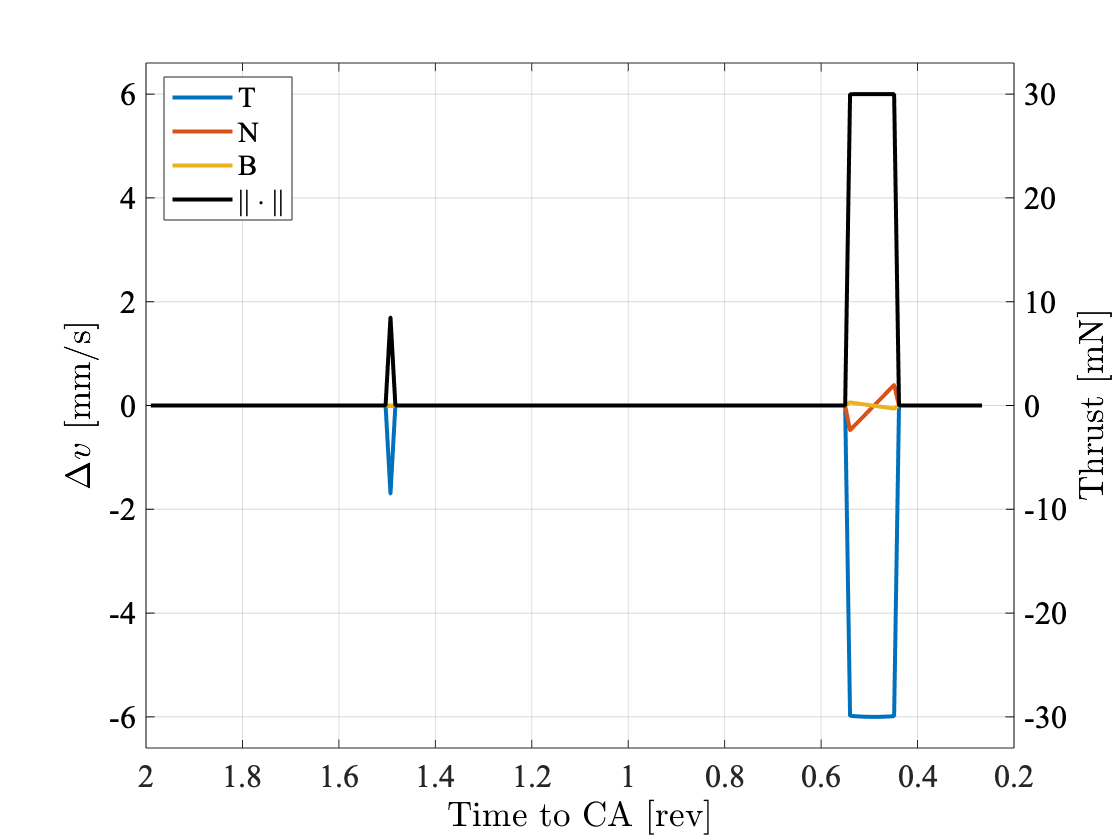}
\label{LinThrust}}\hfill
\subfigure[Fourth major iteration]{
\includegraphics[width=0.47\textwidth]{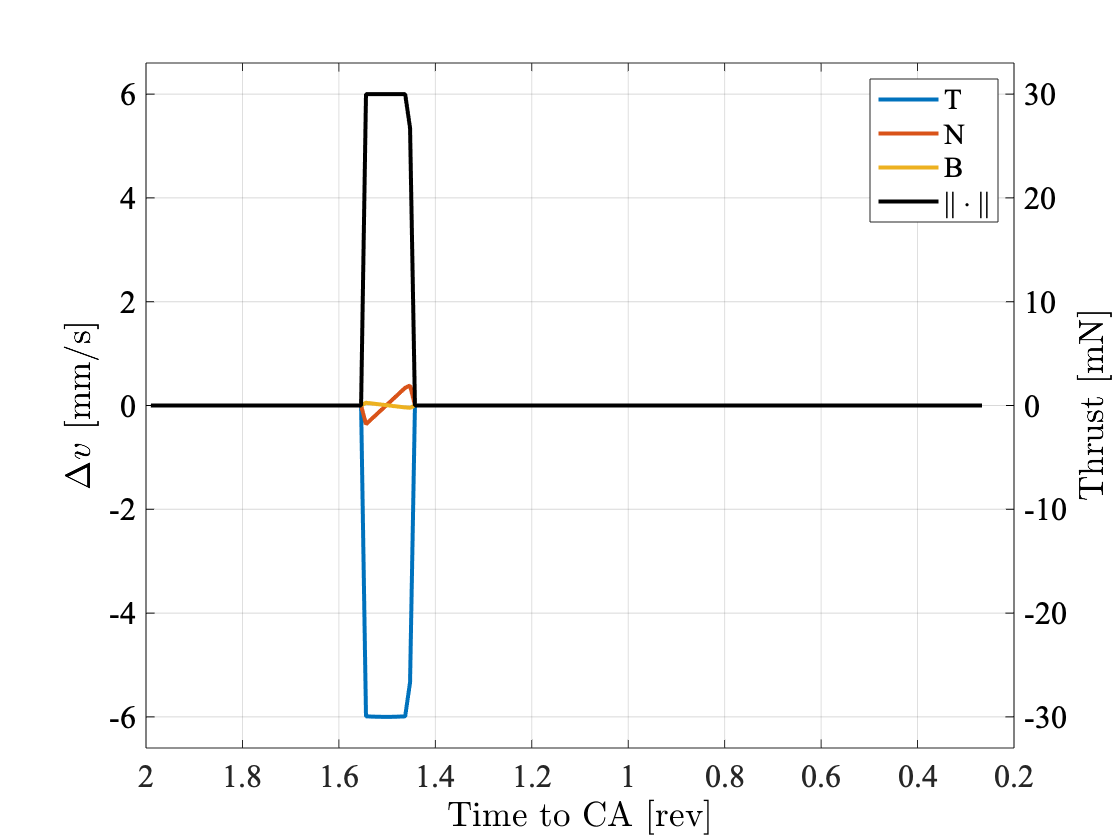}
\label{NonLinThrust}}\hfill
\caption{Effect of nonlinearities on the solution structure for test case \#644}
\label{effectOfNonLin}
\end{figure}

In terms of iterations, the three constraints share a similar behavior, with more than 98\% of the cases requiring only 2 major iterations (with a maximum of 6) and 3-4 minor iterations. In most cases, only one minor iteration is sufficient for the second major iteration, proving the accuracy of the linearized maps. The cases that require more major iterations (i.e., multiple linearizations) are highly sensitive to small variations in the linear maps $\mathbb{A}, \mathbb{B}, \mathbb{C}$. In most cases, this results in a limited number of impulses moving between the end points of a thrust arc, maintaining the total $\Delta v$ almost unchanged. A more interesting case is presented in Fig. \ref{effectOfNonLin} for the case \#644 when the constraint on maximum collision risk is considered. This test case has an encounter velocity of only 94.53 m/s, which is on the boundary of the short-term encounter approximation \cite{Chan2008}. Interestingly, 
while at the first major iteration the maneuver is mainly executed around the last opposition, at the last one, it is entirely performed across the first one. Without the inclusion of the major iterations (i.e. using a single linearization), the solution in Fig. \ref{LinThrust} would have been considered the optimal one. Nevertheless, the difference in the $\Delta v$ is limited with the converged solution consuming 59.3 mm/s, and the first iteration one 61.7 mm/s. Additionally, the constraint violation at the first iteration is also limited, with $P_{C,\max} = 9.863 \times 10^{-5}$. Due to the low encounter velocity, the \gls{cam} changes the time of closest approach by 7 seconds. This is considerably higher than for shorter-term encounters, where the average variation of the time of closest approach is less than half-second.  

\begin{figure}[h!]
\subfigure[Distribution of $\Delta v$ difference]{
\includegraphics[width=0.47\textwidth]{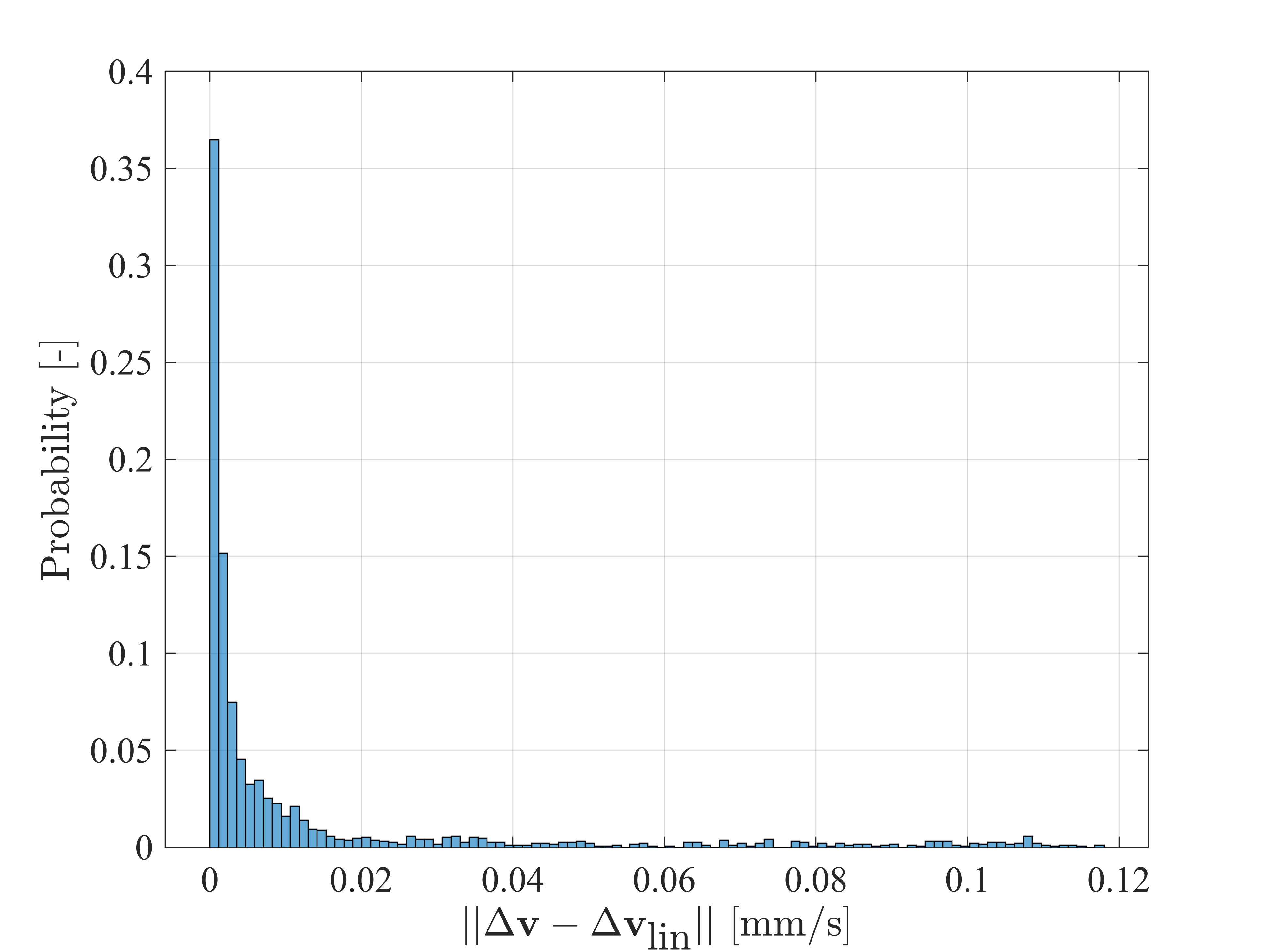}
\label{diffdvLin}}\hfill
\subfigure[Distribution of $\Delta r_{CA}$ difference]{
\includegraphics[width=0.47\textwidth]{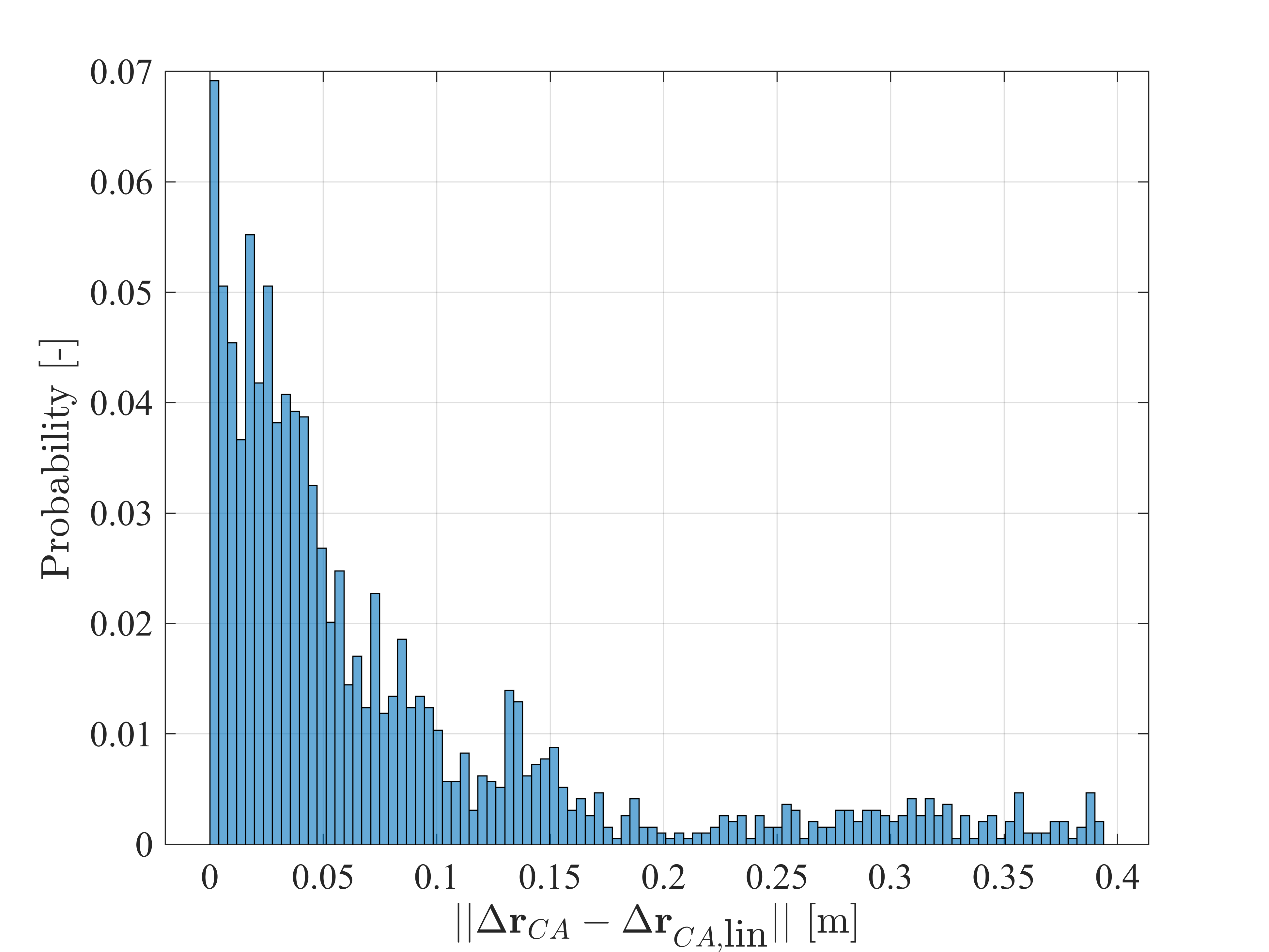}
\label{diffbplaneLin}}\hfill
\caption{Effect of the linearization on \gls{cam} design}
\label{linDifference}
\end{figure}

Finally, to complete this assessment, the miss distance is considered in Fig. \ref{linDifference}. Figure \ref{diffdvLin} shows the difference between the $\Delta v$ computed by Algorithm \ref{algo:SCVX} with respect to the one at the end of the first major iteration, indicated with the subscript {\it lin}. In Fig. \ref{diffbplaneLin} the same analysis is performed for the post-maneuver position on the b-plane. For the vast majority of cases the effect of nonlinearites is negligible, particularly when short alert times are considered, as in this case.    

\begin{figure}[h!]
\subfigure[Distribution of $\Delta v$ difference ]{
\includegraphics[width=0.47\textwidth]{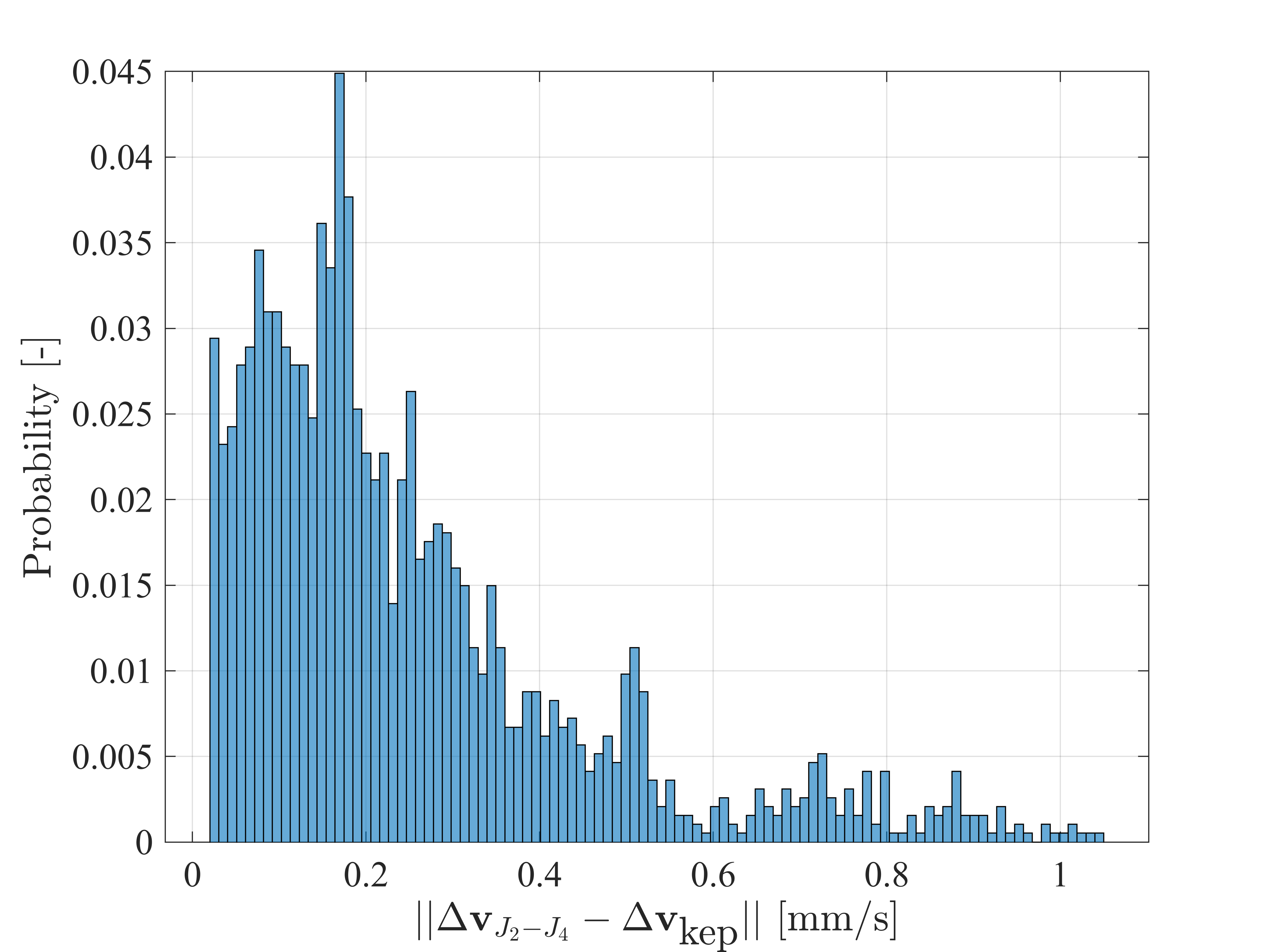}
\label{diffdvKepJ2}}\hfill
\subfigure[Distribution of $\Delta r_{CA}$ difference]{
\includegraphics[width=0.47\textwidth]{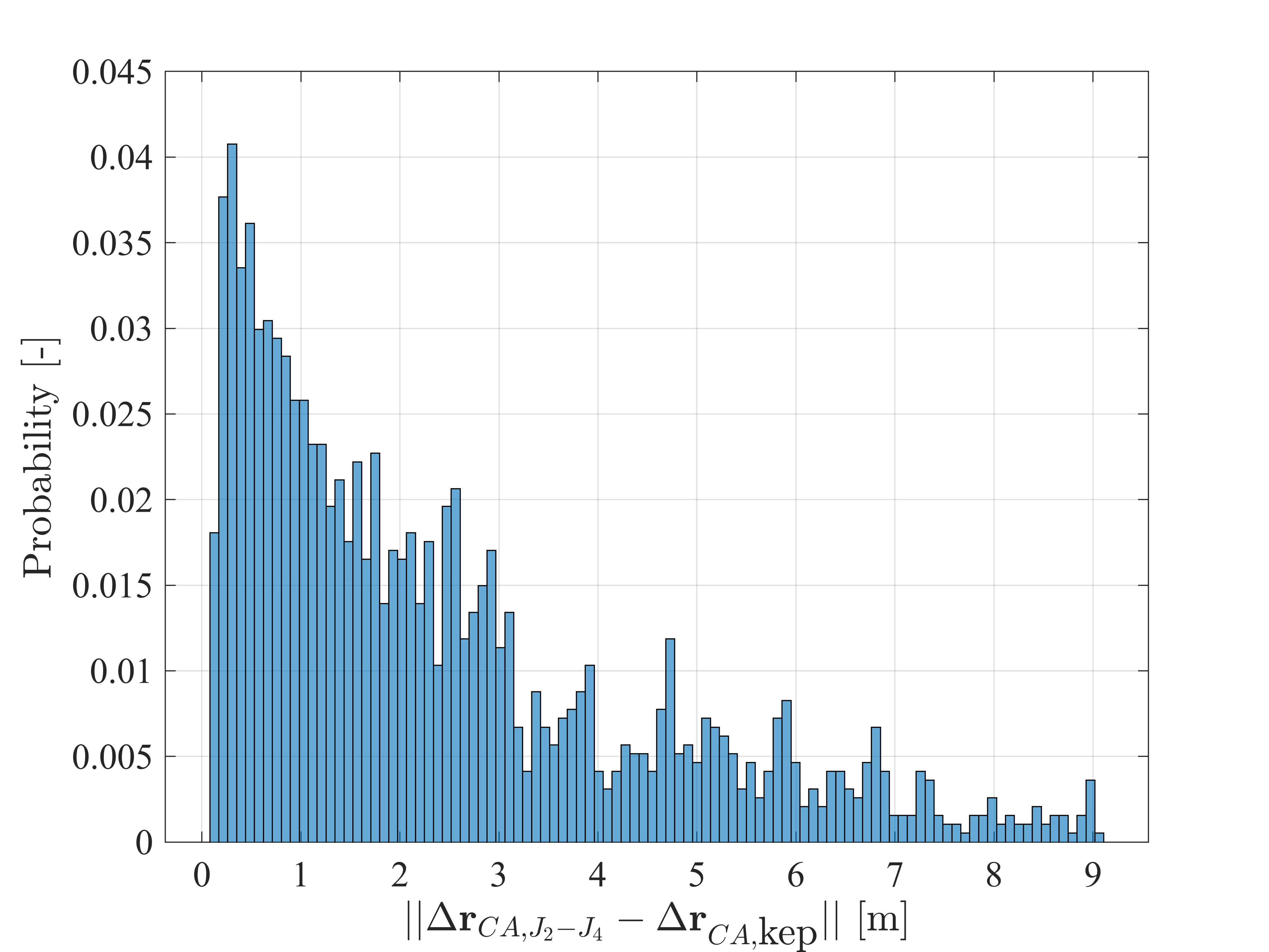}
\label{diffbplaneKepJ2}}\hfill
\caption{Effect of perturbations on \gls{cam} design}
\label{J2Difference}
\end{figure}

The effect of orbital perturbations is studied in Fig. \ref{J2Difference}, in which the results using the dynamics of Eq. \eqref{EoM} are compared against those with a Keplerian model for the miss distance case. As first noted by Patera and Peterson \cite{Patera2003}, the contribution from orbital perturbations are generally small on \gls{cam} design as the changes in the trajectory due to a maneuver are small enough to render the contributions from orbital perturbations negligible with respect to conjunction geometry. In the 90\% of the cases, the $J_2-J_4$ perturbations produce a variation of less than 1 mm/s on the $\Delta v$ and of less than 10 meters the position on the b-plane. 
However, there are cases in which the solution is highly sensitive to the small variations introduced by the perturbations. An example is the test case \#889, for which the solution converged to $\Delta v = 617.1$ mm/s in Keplerian dynamics and to $\Delta v = 1005.6$ mm/s in $J_2-J_4$ dynamics; i.e., to two different minima starting from the same initial guess. 

\begin{figure}[h!]
\centering
\includegraphics[width=0.5\textwidth]{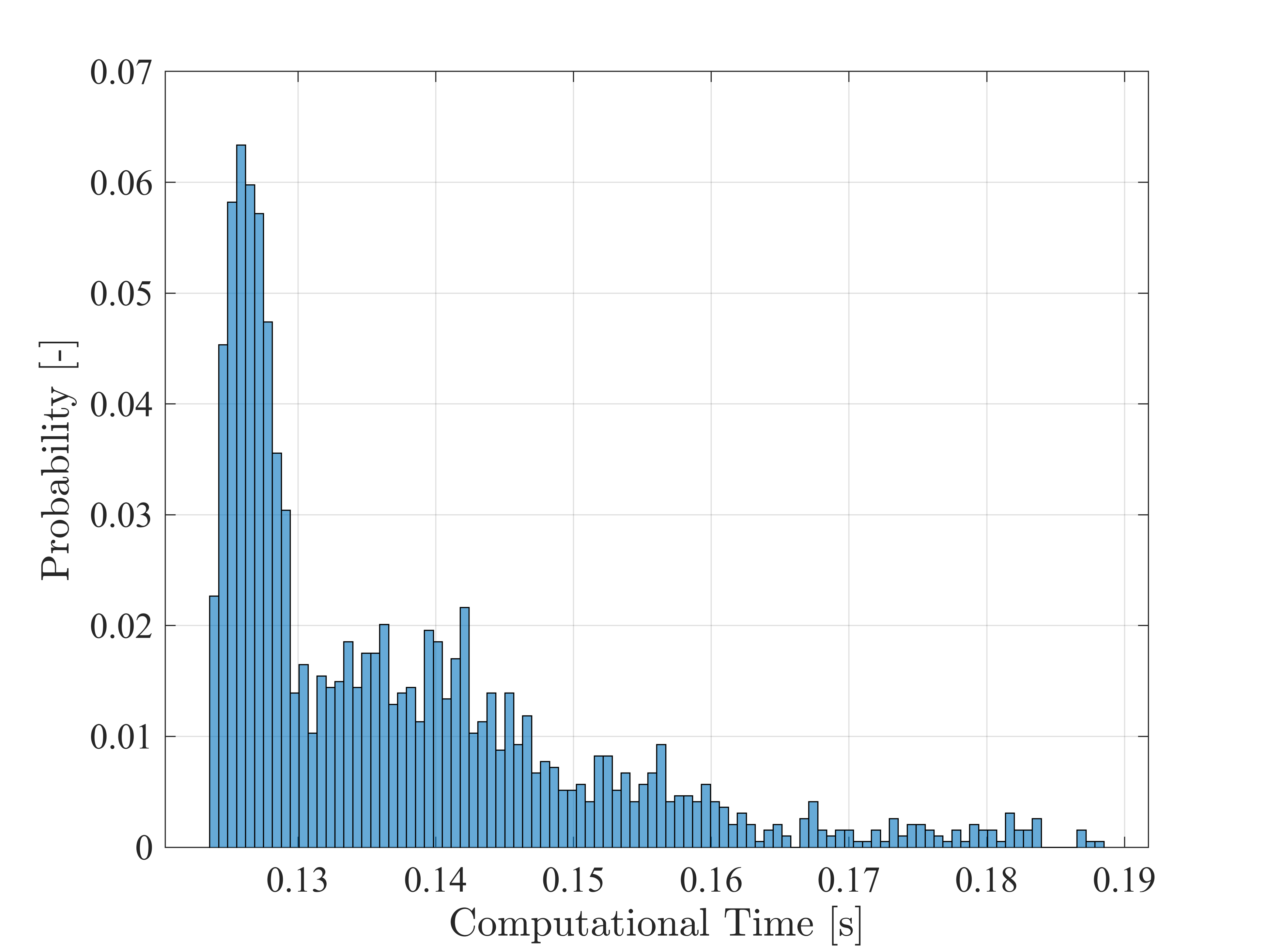}  
\caption{Distribution of computational time}
\label{comTime}
\end{figure}

Figure \ref{comTime} reports the distribution in computational time for the miss distance case, similar results are achieved for the maximum risk and risk ones.  
It is apparent that the convergence properties and efficiency of the algorithm enable the autonomous and quick computation of \gls{cam}s, even potentially onboard. For the latter statement, an assessment of the computational overhead when adapting the codes to be run on an embedded system is needed. However, efficient \gls{socp} solvers tailored for embedded systems usage are already available \cite{Domahidi2013}.

\section{Conclusion}\label{Sec6}

A method based on lossless and successive convexification was proposed for the optimal design of collision avoidance maneuvers (CAMs). The maneuver is modeled as a set of impulses and is thus suitable for handling high- and low-thrust propulsion systems. Maneuvers with minimum propellant consumption were computed meeting constraints either on an estimate of collision risk, maximum collision risk, or miss distance without any prior knowledge on the thrust arc structure and thrust direction. Alert times from few minutes to several orbital periods were considered. The proposed method's convergence properties and efficiency were proven by optimizing CAMs for 2,170 realistic conjunctions, showing that this methodology is promising for future autonomous usage. In the vast majority of the cases linearized Keplerian dynamics were shown to provide accurate results. 

Future efforts will be directed towards handling long-term and multiple encounters and the use of accelerations as decision variables.  

\section*{Acknowledgement}\label{Acknowledgement}
The author acknowledges the work of S\'ebastien Henry for the reconstruction of the conjunction geometries. Laura Pirovano, Harry Holt, and Christian Hofmann's feedback helped improving the quality of the manuscript. Lastly, I am grateful for the reviewers' valuable comments and suggestions. 

\clearpage

\section*{Appendix}\label{sec:simdata}
The data used for the simulations presented in Sec. \ref{sec:MandmIterations}--\ref{sec:maxT} are reported below, those for the extensive simulations are available at \url{github.com/arma1978/conjunction}.

\begin{table}[!ht]
 \centering
\footnotesize{
 \begin{tabular}{rrr}
   \multicolumn{3}{l}{\# Primary}\\
  {\# ECI J2000 Position [km]}  & {\# ECI J2000 Velocity [km/s]} &\\
  $2.33052185175137E+00$ &   $-7.44286282871773E+00$ &\\
  $-1.10370451050201E+03$ &  $-6.13734743652660E-04$ &\\
  $7.10588764299718E+03$ &   $3.95136139293349E-03$ &\\
  \multicolumn{2}{l}{\hspace{1.4cm}\# Covariance matrix RTN [km$^2$]}  & \\
  $9.31700905887535E-05$ & $-2.623398113500550E-04$ & $ 2.360382173935300E-05$  \\
  $-2.623398113500550E-04$ & $1.77796454279511E-02$ & $-9.331225387386501E-05$  \\
  $ 2.360382173935300E-05$ & $-9.331225387386501E-05$ & $1.917372231880040E-05$  \\
   \multicolumn{3}{l}{\# Secondary}\\
  {\# ECI J2000 Position [km]}  & {\# ECI J2000 Velocity [km/s]} &\\
  $2.333465506263321E+00$ &  $7.353740487126315E+00$ &\\
  $-1.103671212478364E+03$ &  $-1.142814049765362E+00$ &\\
  $7.105914958099038E+03$ &  $-1.982472259113771E-01$ &\\
  \multicolumn{2}{l}{\hspace{1.4cm}\# Covariance matrix RTN [km$^2$]}  & \\
  $6.346570910720371E-04$ & $-1.962292216245289E-03$ & $7.077413655227660E-05$  \\
  $-1.962292216245289E-03$ & $8.199899363150306E-01$ & $1.139823810584350E-03$  \\
  $7.077413655227660E-05$ & $1.139823810584350E-03$ & $2.510340829074070E-04$  \\
   \multicolumn{1}{l}{ \# Conjunction details} & \multicolumn{1}{l}{\# Note} & \\
   $R = 29.71$  & m & \\
   $d_{CA}^2 = 8.71655401455392E-01$  & km$^2$ &\\
  $P_C = 1.36040828266536E-01$ & Eq. 5a in \cite{Alfano2007a} & \\
  $P_C = 1.47559666159940E-01$  & Eq. \eqref{eq:approxPC} \\ 
  $P_{C,\max} = 1.92590968666693E-01$ & Eq. \eqref{eq:maxPc} & \\
 \end{tabular}}
 \label{tb:data}
\end{table}

\bibliographystyle{unsrt}  
\bibliography{reference}

\end{document}